%% file: article.tex
\documentclass[preprint,3p,times]{elsarticle}

\usepackage[section]{placeins}
\usepackage{multirow}
\usepackage{makecell}
\usepackage{longtable}
\usepackage{booktabs}
\usepackage{supertabular}
\usepackage[table]{xcolor}

\usepackage{bm}
\usepackage{amssymb}   
\usepackage{amsmath} 

\usepackage{subeqnarray}
\usepackage{cases}  

\usepackage{arydshln}

\usepackage{tikz}
\usepackage{etoolbox}
\usepackage{enumitem}
\newcommand*{\circled}[1]{\lower.7ex\hbox{\tikz\draw (0pt, 0pt)%
    circle (.5em) node {\makebox[1em][c]{\small #1}};}}
\robustify{\circled}
\newcommand*{\blackcircled}[1]{{\lower.7ex\hbox\tikz\draw[fill=black] (0pt, 0pt)%
    circle (.5em) node {\makebox[1em][c]{\small\color{white} #1}};}}
\robustify{\blackcircled}

\usepackage{graphicx}
\usepackage{epstopdf}

\usepackage[bf]{caption2}

\usepackage[ruled,linesnumbered,lined,commentsnumbered]{algorithm2e}

\newtheorem{theorem}{Theorem}
\newtheorem{lemma}{Lemma}

\newtheorem{definition}{Definition}
\newtheorem{corollary}{Corollary}

\newtheorem{example}{Example}

\journal{Elsevier}

\begin{document}

\begin{frontmatter}

\title{Towards building the OP-Mapped WENO schemes: A general 
methodology}


\author[a]{Ruo Li}
\ead{rli@math.pku.edu.cn}

\author[b,c]{Wei Zhong\corref{cor1}}
\ead{zhongwei2016@pku.edu.cn}

\cortext[cor1]{Corresponding author}

\address[a]{CAPT, LMAM and School of Mathematical Sciences, Peking
University, Beijing 100871, China}

\address[b]{School of Mathematical Sciences, Peking University,
Beijing 100871, China}

\address[c]{Northwest Institute of Nuclear Technology, Xi'an 
710024, China}

\input{article_abstract}

\begin{keyword}
Order-preserving mapping \sep OP-Mapped WENO \sep Hyperbolic 
conservation laws


\end{keyword}

\end{frontmatter}

\input{article_introduction}
\input{article_review}
\input{article_GMG_WENO_X}
\input{article_numerics}
\input{article_conclusion}


\bibliographystyle{model1b-shortjournal-num-names}
\bibliography{../refs}

\end{document}

%% file: article_abstract.tex
\begin{abstract}

  A serious and ubiquitous issue in existing mapped WENO schemes is 
  that most of them can hardly preserve high resolutions and in 
  the meantime prevent spurious oscillations on solving hyperbolic 
  conservation laws with long output times. Our goal in this article 
  is to address this widely concerned problem \cite{WENO-PM,WENO-IM,
  WENO-PPM5,WENO-RM260,WENO-MAIMi,MOP-WENO-ACMk}. In our previous 
  work \cite{MOP-WENO-ACMk}, the \textit{order-preserving (OP)} 
  criterion was originally introduced and carefully used to devise a 
  new mapped WENO scheme that performs satisfactorily in long-run 
  simulations, and hence it was indicated that the \textit{OP} 
  criterion plays a critical role in the maintenance of 
  low-dissipation and robustness for the mapped WENO schemes. Thus, 
  in our present work, we firstly define the family of the mapped 
  WENO schemes, whose mappings meet the \textit{OP} criterion, as 
  OP-Mapped WENO. Next, we attentively take a closer look at the 
  mappings of various existing mapped WENO schemes and devise a 
  general formula for them. It helps us to extend the \textit{OP} 
  criterion into the design of the improved mappings. Then, we 
  propose the generalized implementation of obtaining a group of 
  OP-Mapped WENO schemes, named MOP-WENO-X as they are developed 
  from the existing mapped WENO-X schemes, where the notation ``X'' 
  is used to identify the version of the existing mapped WENO 
  scheme. Finally, extensive numerical experiments and comparisons 
  with competing schemes are conducted to demonstrate the enhanced 
  performances of the MOP-WENO-X schemes.

\end{abstract}


%% file: article_introduction.tex
\section{Introduction}
\label{secIntroduction}
The essentially non-oscillatory (ENO) schemes \cite{ENO1987JCP71, 
ENO1987V24, ENO1986, ENO1987JCP83} and the weighted ENO (WENO) 
schemes \cite{ENO-Shu1988, ENO-Shu1989, WENO-LiuXD,WENO-JS, 
WENOoverview} have developed quite successfully in recent decades to 
solve the hyperbolic conservation problems, especially those that 
may generate discontinuities and smooth small-scale structures as 
time evolves in its solution even if the initial condition is 
smooth. The main purpose of this paper is to find a general method 
to introduce the \textit{order-preserving (OP)} mapping proposed in 
our previous work \cite{MOP-WENO-ACMk} for improving the existing 
mapped WENO schemes for the approximation of the hyperbolic 
conservation laws in the form
\begin{equation}
\dfrac{\partial \mathbf{u}}{\partial t} +
\nabla \cdot \mathbf{F}(\mathbf{u}) = 0,
\label{eq:governingEquation}
\end{equation}
where $\mathbf{u} = (u_{1}, \cdots, u_{m}) \in \mathbb{R}^{m}$ is 
the vector of the conserved variables and $\mathbf{F}(\mathbf{u})$ 
is the vector of the Cartesian components of flux.  

Harten et al. \cite{ENO1987JCP71} introduced the ENO schemes. They 
used the smoothest stencil from $r$ possible candidate stencils 
based on the local smoothness to perform a polynomial reconstruction 
such that it yields high-order accuracy in smooth regions but avoids 
spurious oscillations at or near discontinuities. Liu, Osher, and 
Chan \cite{WENO-LiuXD} introduced the first WENO scheme, an improved 
version of the ENO methodology with a cell-averaged approach, by 
using a nonlinear convex combination of all the candidate stencils 
to achieve a higher order of accuracy than the ENO schemes while 
retaining the essential non-oscillatory property at or near 
discontinuities. In other words, it achieves $(r + 1)$th-order of 
accuracy from the $r$th-order ENO schemes \cite{ENO1987JCP71,
ENO1987V24,ENO1986} in the smooth regions while behaves similarly to 
the $r$th-order ENO schemes in regions including discontinuities. In 
\cite{WENO-JS}, Jiang and Shu proposed the classic WENO-JS scheme 
with a new measurement of the smoothness of the numerical solutions 
on substencils (hereafter, denoted by smoothness indicator), by 
using the sum of the normalized squares of the scaled $L_{2}$-norms 
of all the derivatives of $r$ local interpolating polynomials, to 
obtain $(2r - 1)$th-order of accuracy from the $r$th-order ENO 
schemes.

The WENO-JS scheme has become a very popular and quite successful 
methodology for solving compressible flows modeled through 
hyperbolic conservation laws in the form of 
Eq.\eqref{eq:governingEquation}. However, it was less than 
fifth-order for many cases such as at or near critical points of 
order $n_{\mathrm{cp}} = 1$ in the smooth regions. Here, we refer to 
$n_{\mathrm{cp}}$ as the order of the critical point; e.g., 
$n_{\mathrm{cp}} = 1$ corresponds to $f'=0, f'' \neq 0$ and 
$n_{\mathrm{cp}} = 2$ corresponds to $f'=0, f'' = 0, f''' \neq 0$, 
etc. In particular, Henrick et al. \cite{WENO-M} identified that the 
fifth-order WENO-JS scheme fails to yield the optimal convergence 
order at or near critical points where the first derivative vanishes 
but the third derivative does not simultaneously. Then, in the 
same article, they derived the necessary and sufficient conditions 
on the nonlinear weights for optimality of the convergence rate of 
the fifth-order WENO schemes and these conditions were reduced to a 
simpler sufficient condition \cite{WENO-Z} which could be easily 
extended to the $(2r-1)$th-order WENO schemes \cite{WENO-IM}. 
Moreover, also in \cite{WENO-M}, Henrick et al. devised the original 
mapped WENO scheme, named WENO-M hereafter, by constructing a 
mapping function that satisfies the sufficient condition to achieve 
the optimal order of accuracy.

Later, following the idea of incorporating a mapping procedure to 
maintain the nonlinear weights of the convex combination of stencils 
as near as possible to the ideal weights of optimal order accuracy, 
various versions of mapped WENO schemes have been successfully 
proposed. In \cite{WENO-IM}, Feng et al. rewrote the mapping 
function of the WENO-M scheme in a simple and more meaningful form 
and then extended it to a general class of improved mapping 
functions leading to the family of the WENO-IM($k, A$) schemes, 
where $k$ is a positive even integer and $A$ a positive real number. 
It was indicated that by taking $k=2$ and $A=0.1$ in the 
WENO-IM($k,A$) scheme, far better numerical solutions with less 
dissipation and higher resolution could be obtained than that of the 
WENO-M scheme. Unfortunately, the numerical experiments in 
\cite{WENO-RM260} showed that the seventh- and ninth- order 
WENO-IM(2, 0.1) schemes generated evident spurious oscillations near 
discontinuities when the output time is large. In addition, our 
numerical experiments as shown in Fig. \ref{fig:SLP:N3200:PM6}, Fig. 
\ref{fig:BiCWP:N800:PM6} and Fig. \ref{fig:BiCWP:N3200:PM6} of this 
paper indicate that, even for the fifth-order WENO-IM(2, 0.1) 
scheme, the spurious oscillations are also produced when the grid 
number increases or a different initial condition is used. Recently, 
Feng et al. \cite{WENO-PM} pointed out that, when the WENO-M scheme 
was used for solving the problems with discontinuities for long 
output times, its mapping function may amplify the effect from the 
non-smooth stencils leading to a potential loss of accuracy near 
discontinuities. To amend this drawback, a piecewise polynomial 
mapping function with two additional requirements, that is, 
$g'(0) = 0$ and $g'(1)=0$ ($g(x)$ denotes the mapping function), to 
the original criteria in \cite{WENO-M} was proposed. The recommended 
WENO-PM6 scheme \cite{WENO-PM} achieved significantly higher 
resolution than the WENO-M scheme when computing the one-dimensional 
linear advection problem with long output times. However, it may 
generate spurious oscillations near the discontinuities as shown 
in Fig. 8 of \cite{WENO-IM} and Figs. 3-8 of \cite{WENO-RM260}. 

Many other mapped WENO schemes, such as the WENO-PPM$n$
\cite{WENO-PPM5}, WENO-RM($mn0$) \cite{WENO-RM260}, WENO-MAIM$i$
\cite{WENO-MAIMi}, WENO-ACM \cite{WENO-ACM} schemes and et al., have
been successfully developed to enhance the performances of the 
classic WENO-JS scheme in some ways, like achieving optimal 
convergence orders near critical points in smooth regions, having 
lower numerical dissipations, achieving higher resolutions near 
discontinuities, or reducing the computational costs, and we refer 
to the references for more details. However, as mentioned in 
previously published literatures \cite{WENO-IM,WENO-RM260}, most of 
the existing improved mapped WENO schemes could not prevent the 
spurious oscillations near discontinuities, especially for long 
output time simulations. Moreover, when simulating the 
two-dimensional steady problems with strong shock waves, the 
post-shock oscillations, which were systematically studied for WENO 
schemes by Zhang et al. \cite{WENO-ZS}, become very severe in the 
solutions of most of the existing improved mapped WENO schemes 
\cite{WENO-ACM}.

In the previous work of this paper \cite{MOP-WENO-ACMk}, the authors
made a further study of the nonlinear weights of the existing mapped 
WENO schemes, by taking the ones developed in \cite{WENO-PM,WENO-IM,
WENO-MAIMi,MOP-WENO-ACMk} as examples. It was found that the order 
of the nonlinear weights for the substencils of the same global 
stencil has been changed at many points in the mapping process of 
all these considered mapped WENO schemes. The order-change of the 
nonlinear weights is caused by weights increasing of non-smooth 
substencils and weights decreasing of smooth substencils. 
It was revealed that this is the essential cause of the potential 
loss of accuracy of the WENO-M scheme and the spurious oscillation 
generation of the existing improved mapped WENO schemes, through 
theoretical analysis and extensive numerical tests. In the same 
article, the definition of the \textit{order-preserving (OP)} 
mapping was given and suggested as an additional criterion in the 
design of the mapping function. Then a new mapped WENO scheme with 
its mapping function satisfying the additional criterion was 
proposed. Extensive numerical experiments showed that this scheme 
can achieve the optimal convergence order of accuracy even at 
critical points. It also can decrease the numerical dissipations and 
obtain high resolution but does not generate spurious oscillation 
near discontinuities even if the output time is large. Moreover, it 
was observed clearly that it exhibits a significant advantage in 
reducing the post-shock oscillations when calculating the steady 
problems with strong shock waves in two dimension.

In this article, the idea of introducing the 
\textit{OP} criterion into the design of the mapping 
functions proposed in \cite{MOP-WENO-ACMk} is extended to various 
existing mapped WENO schemes. First of all, we give a common name of 
\textit{OP-Mapped WENO} to the family of the mapped WENO schemes 
whose mappings are \textit{OP}. A general formula for the 
mapping functions of various existing mapped WENO schemes is 
obtained that allows the extension of the \textit{OP} criterion to 
all existing mapped WENO schemes. The notation MOP-WENO-X is used to 
denote the improved mapped WENO scheme considering the \textit{OP} 
criterion based on the existing WENO-X scheme. A new function named 
\textbf{minDist} is defined (see Definition \ref{definition:minDist} 
in subsection \ref{subsec:newMapping} below). A general algorithm to 
construct \textit{OP} mappings through the existing mapping 
functions by using the \textbf{minDist} function is proposed. 
Extensive numerical tests are conducted to demonstrate the 
performances of the MOP-WENO-X schemes: (1) a series of accuracy 
tests shows the capacity of the MOP-WENO-X schemes to achieve the 
optimal convergence order in smooth regions with first-order 
critical points and their advantages in long output time simulations 
of the problems with very high-order critical points; (2) the 
one-dimensional linear advection equation with two kinds of initial 
conditions for long output times are then presented to demonstrate 
that the MOP-WENO-X schemes can obtain high resolution and meanwhile 
avoid spurious oscillation near discontinuities; (3) some benchmark 
tests of steady problems with strong shock waves modeled via the 
two-dimensional Euler equations are computed, it is clear that the 
MOP-WENO-X schemes also enjoy a significant advantage in reducing 
the post-shock oscillations.

The remainder of this paper is organized as follows. In Section 
\ref{secMappedWENO}, we briefly review the preliminaries to 
understand the finite volume method and the procedures of the WENO-JS
\cite{WENO-JS}, WENO-M \cite{WENO-M} and other versions of mapped 
WENO schemes. Section \ref{MOP-WENO-X} presents a general method to 
introduce the \textit{OP} mapping for improving the existing 
mapped WENO schemes. Some numerical results are provided in Section 
\ref{NumericalTests} to illustrate the advantages of the proposed 
WENO schemes. Finally, concluding remarks are given in Section 
\ref{secConclusions} to close this paper.


%% file: article_review.tex
\section{Brief review of the WENO schemes}
\label{secMappedWENO}

For simplicity of presentation but without loss of generality, we 
denote our discussion to the following one-dimensional linear 
hyperbolic conservation equation
\begin{equation}
\dfrac{\partial u}{\partial t}+\dfrac{\partial f(u)}{\partial x}=0,
\quad x_{l} < x < x_{r}, t > 0,
\label{eq:1D-hyperbolicLaw}
\end{equation}
with the initial condition $u(x,0) = u_{0}(x)$. We only confine our 
attention to the uniform meshes in this paper and the WENO method 
with non-uniform meshes can refer to \cite{non-uniform-WENO-01,
non-uniform-WENO-02}. Throughout this paper, we assume that the 
given domain $[x_{l}, x_{r}]$ is discretized into the set of uniform 
cells $I_{j} := [x_{j-1/2}, x_{j+1/2}], j = 1,\cdots,N$ with the 
cell size $\Delta x = \frac{x_{r} - x_{l}}{N}$. The associated cell 
centers and cell boundaries are denoted by $x_{j}=x_{l} + (j - 1/2)
\Delta x$ and $x_{j \pm 1/2} = x_{j} \pm \Delta x/2$, respectively. 
The notation $\bar{u}(x_{j}, t)=\dfrac{1}{\Delta x}
\int_{x_{j-1/2}}^{x_{j+1/2}}u(\xi,t)\mathrm{d}\xi$ indicates the 
cell average of $I_{j}$. The one-dimensional linear hyperbolic 
conservation equation in Eq.\eqref{eq:1D-hyperbolicLaw} can be 
approximated by a system of ordinary differential equations, 
yielding the semi-discrete finite volume form:
\begin{equation}
\begin{array}{l}
\dfrac{\mathrm{d}\bar{u}_{j}(t)}{\mathrm{d}t}
\approx \mathcal{L}(u_{j}), \\
\mathcal{L}(u_{j}) = -\dfrac{1}{\Delta x}\Big( \hat{f}_{j+1/2} - 
\hat{f}_{j-1/2} \Big),
\end{array}
\label{eq:discretizedFunction}
\end{equation}
where $\bar{u}_{j}(t)$ is the numerical approximation to the cell 
average $\bar{u}(x_{j}, t)$, and the numerical flux $\hat{f}_{j \pm
1/2}$ is a replacement of the physical flux function $f(u)$ at the 
cell boundaries $x_{j \pm 1/2}$ and it is defined by $\hat{f}_{j \pm
1/2} = \hat{f}(u_{j \pm 1/2}^{-}, u_{j \pm 1/2}^{+})$. The notations 
$u_{j \pm 1/2}^{\pm}$ refer to the limits of $u$ and their values of 
$u_{j \pm 1/2}^{\pm}$ can be obtained by the technique of 
reconstruction, like the WENO reconstruction procedures narrated 
later. In this paper, we choose the global Lax-Friedrichs flux
\begin{equation*}
\hat{f}(a,b) = \frac{1}{2}\big[f(a) + f(b) -\alpha(b - a)\big], 
\end{equation*}
where $\alpha = \max_{u} \lvert f'(u) \rvert$ is a constant and the 
maximum is taken over the whole range of $u$.

\subsection{The WENO-JS reconstruction}
Firstly, we review the process of the classic fifth-order WENO-JS 
reconstruction \cite{WENO-JS}. For brevity, we describe only the 
reconstruction procedure of the left-biased $u_{j + 1/2}^{-}$, and 
the right-biased one $u_{j + 1/2}^{+}$ can trivially be computed by 
mirror symmetry with respect to the location $x_{j + 1/2}$ of 
$u_{j + 1/2}^{-}$. We drop the subscript ``-'' below just for 
simplicity of notation. 

To construct the values of $u_{j + 1/2}$ from known cell average 
values $\overline{u}_{j}$, a 5-point global stencil 
$S^{5}=\{I_{j-2}, I_{j-1}, I_{j}, I_{j+1},$\\$I_{j+2}\}$ is used in 
the fifth-order WENO-JS scheme. It is subdivided into three 3-point 
substencils $S_{s} = \{I_{j+s-2}, I_{j+s-1}, I_{j+s}\}$ with 
$s=0, 1, 2$. It is known that the third-order approximations of 
$u(x_{j+1/2}, t)$ associated with these substencils are explicitly 
given by
\begin{equation}
\begin{array}{l}
\begin{aligned}
&u_{j+1/2}^{0} = \dfrac{1}{6}(2\bar{u}_{j-2} - 7\bar{u}_{j-1}
+ 11\bar{u}_{j}), \\
&u_{j+1/2}^{1} = \dfrac{1}{6}(-\bar{u}_{j-1} + 5\bar{u}_{j}
+ 2\bar{u}_{j+1}), \\
&u_{j+1/2}^{2} = \dfrac{1}{6}(2\bar{u}_{j} + 5\bar{u}_{j+1}
- \bar{u}_{j+2}).
\end{aligned}
\end{array}
\label{eq:approx_ENO}
\end{equation}

Then the $u_{j + 1/2}$ of global stencil $S^{5}$ is computed by a 
weighted average of those third-order approximations of substencils, 
taking the form
\begin{equation}
u_{j + 1/2} = \sum\limits_{s = 0}^{2}\omega_{s}u_{j + 1/2}^{s}.
\label{eq:approx_WENO}
\end{equation}
The nonlinear weights $\omega_{s}$ in the classic WENO-JS scheme are 
defined as
\begin{equation} 
\omega_{s}^{\mathrm{JS}} = \dfrac{\alpha_{s}^{\mathrm{JS}}}{\sum_{l =
 0}^{2} \alpha_{l}^{\mathrm{JS}}}, \alpha_{s}^{\mathrm{JS}} = \dfrac{
 d_{s}}{(\epsilon + \beta_{s})^{2}}, \quad s = 0,1,2,
\label{eq:weights:WENO-JS}
\end{equation} 
where $d_{0}, d_{1}, d_{2}$ are called the ideal weights of 
$\omega_{s}$ since they generate the central upstream fifth-order 
scheme for the global stencil $S^{5}$. It is known that 
$d_{0} = 0.1, d_{1} = 0.6, d_{2} = 0.3$ and in smooth regions we can 
get $\sum\limits_{s=0}^{2} d_{s} u^{s}_{j+1/2} = u(x_{j + 1/2}, t) + 
O(\Delta x^{5})$. $\epsilon$ is a small positive number introduced 
to prevent the denominator from becoming zero. The parameters 
$\beta_{s}$ are the smoothness indicators for the third-order 
approximations $u_{j+1/2}^{s}$ and their explicit formulas can be 
obtained from \cite{WENO-JS}, taking the form
\begin{equation*}
\begin{array}{l}
\begin{aligned}
\beta_{0} &= \dfrac{13}{12}\big(\bar{u}_{j - 2} - 2\bar{u}_{j - 1} + 
\bar{u}_{j} \big)^{2} + \dfrac{1}{4}\big( \bar{u}_{j - 2} - 4\bar{u}_
{j - 1} + 3\bar{u}_{j} \big)^{2}, \\
\beta_{1} &= \dfrac{13}{12}\big(\bar{u}_{j - 1} - 2\bar{u}_{j} + \bar
{u}_{j + 1} \big)^{2} + \dfrac{1}{4}\big( \bar{u}_{j - 1} - \bar{u}_{
j + 1} \big)^{2}, \\
\beta_{2} &= \dfrac{13}{12}\big(\bar{u}_{j} - 2\bar{u}_{j + 1} + \bar
{u}_{j + 2} \big)^{2} + \dfrac{1}{4}\big( 3\bar{u}_{j} - 4\bar{u}_{j 
+ 1} + \bar{u}_{j + 2} \big)^{2}.
\end{aligned}
\end{array}
\end{equation*}

In general, the fifth-order WENO-JS scheme is able to recover the
optimal convergence rate of accuracy in smooth regions. However, 
when at or near critical points where the first derivative vanishes 
but the third derivative does not simultaneously, it loses accuracy 
and its convergence rate of accuracy decreases to third-order or 
even less. We refer to \cite{WENO-M} for more details.

\subsection{The mapped WENO reconstructions}
\label{subsecWENO-M}
To address the issue of the WENO-JS scheme mentioned above, Henrick 
et al. \cite{WENO-M} made a systematic truncation error analysis of 
Eq.\eqref{eq:discretizedFunction} in its corresponding finite 
difference version by using the Taylor series expansions of the Eq.
\eqref{eq:approx_ENO}, and hence they derived the necessary and 
sufficient conditions on the weights for the fifth-order WENO scheme 
to achieve the formal fifth-order of convergence at smooth regions 
of the solution, taking the form 
\begin{equation}
\sum\limits_{s=0}^{2}(\omega_{s}^{\pm} - d_{s})=O(\Delta x^{6}), 
\quad
\sum\limits_{s=0}^{2}A_{s}(\omega_{s}^{+} - \omega_{s}^{-})=
O(\Delta x^{3}),
\quad  
\omega_{s}^{\pm} - d_{s} = O(\Delta x^{2}),
\label{eq:nec-suf_cond}
\end{equation}
where the superscripts ``+'' and ``-'' on $\omega_{s}$ correspond to 
their use in either $u_{j+1/2}^{s}$ and $u_{j-1/2}^{s}$ stencils 
respectively. Since the first equation in Eq.\eqref{eq:nec-suf_cond}
always holds due to the normalization, a simpler sufficient 
condition for the fifth-order convergence is given as \cite{WENO-Z}
\begin{equation}
\omega^{\pm}_{s} - d_{s} = O(\Delta x^{3}), \quad s = 0, 1, 2.
\label{eq:suf_cond}
\end{equation}
The conditions Eq.\eqref{eq:nec-suf_cond} or Eq.\eqref{eq:suf_cond}
may not hold in the case of smooth extrema or at critical points 
when the fifth-order WENO-JS scheme is used. An innovative idea of 
fixing this deficiency, originally proposed by Henrick in 
\cite{WENO-M}, is to design a mapping function to make $\omega_{s}$ 
approximating the ideal weights $d_{s}$ at critical points to the 
required third order $O(\Delta x^{3})$. The first mapping function devised by Henrick et al. in \cite{WENO-M} is given as
\begin{equation}
\big( g^{\mathrm{M}} \big)_{s}(\omega) = \dfrac{ \omega \big( d_{s} +
d_{s}^2 - 3d_{s}\omega + \omega^{2} \big) }{ d_{s}^{2} + (1 - 2d_
{s})\omega }, \quad \quad s = 0, 1, 2.
\label{mapFuncWENO-M}
\end{equation}
We can verify that $\big( g^{\mathrm{M}}\big)_{s}(\omega)$ meets the 
conditions in Eq.\eqref{eq:suf_cond} as it is a non-decreasing 
monotone function on $[0, 1]$ with finite slopes and satisfies the 
following properties.
\begin{lemma} 
The mapping function $\big( g^{\mathrm{M}} \big)_{s}(\omega)$
defined by Eq.(\ref{mapFuncWENO-M}) satisfies: \\

C1. $0 \leq \big(g^{\mathrm{M}}\big)_{s}(\omega)\leq 1, \big(
g^{\mathrm{M}}\big)_{s}(0)=0, \big(g^{\mathrm{M}}\big)_{s}(1) = 1$;

C2. $\big( g^{\mathrm{M}} \big)_{s}(d_{s}) = d_{s}$;

C3. $\big( g^{\mathrm{M}} \big)_{s}'(d_{s}) = \big(g^{\mathrm{M}}
\big)_{s}''(d_{s}) = 0$. 
\label{lemmaWENO-Mproperties}
\end{lemma}

Following Henrick's idea, a great many improved mapping functions 
were successfully proposed \cite{WENO-IM,WENO-PM,WENO-PPM5,WENO-RM260,WENO-MAIMi,WENO-ACM,MOP-WENO-ACMk}. To clarify our major concern 
and provide convenience to readers but for brevity in the 
description, we only state some mapping functions in the following 
context, and we refer to references for properties similar to Lemma 
\ref{lemmaWENO-Mproperties} and more details of these mapping 
functions.

$\blacksquare$ WENO-IM$(k, A)$ \cite{WENO-IM}
\begin{equation}
\big( g^{\mathrm{IM}} \big)_{s}(\omega; k, A) = d_{s} + \dfrac{\big( 
\omega - d_{s} \big)^{k + 1}A}{\big( \omega - d_{s} \big)^{k}A + 
\omega(1 - \omega)}, \quad A > 0, k = 2n, n \in \mathbb{N}^{+}.
\label{mappingFunctionWENO-IM}
\end{equation}

$\blacksquare$ WENO-PM$k$ \cite{WENO-PM}
\begin{equation}
\big( g^{\mathrm{PM}} \big)_{s}(\omega) = c_{1}(\omega - d_{s})^{k+1}
(\omega + c_{2}) + d_{s}, \quad \quad k \geq 2,
\label{mappingFunctionWENO-PM}
\end{equation}
where $c_{1}, c_{2}$ are constants with specified parameters $k$ and 
$d_{s}$, taking the following forms
\begin{equation*}
\begin{array}{ll}
c_{1} = \left\{
\begin{array}{ll}
\begin{aligned}
&(-1)^{k}\dfrac{k+1}{d_{s}^{k+1}}, & 0 \leq \omega \leq d_{s}, \\
&-\dfrac{k+1}{(1-d_{s})^{k+1}},    & d_{s} < \omega \leq 1,
\end{aligned}
\end{array}\right.
&
c_{2} = \left\{
\begin{array}{ll}
\begin{aligned}
&\dfrac{d_{s}}{k+1},        &  0 \leq \omega \leq d_{s}, \\
&\dfrac{d_{s}-(k+2)}{k+1},  & d_{s} < \omega \leq 1. 
\end{aligned}
\end{array}\right.
\end{array}
\end{equation*}

$\blacksquare$ WENO-PPM$n$ \cite{WENO-PPM5}
\begin{equation}
\big( g_{s}^{\mathrm{PPM}n} \big)_{s}(\omega) = \left\{
\begin{array}{ll}
\begin{aligned}
&\big( g_{s,\mathrm{L}}^{\mathrm{PPM}n} \big)_{s}(\omega), & \omega 
\in [0, d_{s}]\\
&\big( g_{s,\mathrm{R}}^{\mathrm{PPM}n} \big)_{s}(\omega), & \omega 
\in (d_{s}, 1],
\end{aligned}
\end{array}
\right. 
\label{eq:mappingFunction_WENO-PPMn}
\end{equation}
and for $n = 5$,
\begin{equation}
\big( g_{s,\mathrm{L}}^{\mathrm{PPM}5} \big)_{s}(\omega) = d_{s}
\big( 1 + (a - 1)^{5} \big), \quad 
\big( g_{s,\mathrm{R}}^{\mathrm{PPM}5} \big)_{s}(\omega) = d_{s} + 
b^{4}\big( \omega - d_{s} \big)^{5}.
\label{eq:mappingFunction_WENO-PPM5}
\end{equation}
where $a = \omega/d_{s},b = 1/(d_{s} - 1)$.

$\blacksquare$ WENO-RM$(mn0)$ \cite{WENO-RM260}
\begin{equation}
\big( g^{\mathrm{RM}} \big)_{s}(\omega) = d_{s} + \dfrac{(\omega - 
d_{s})^{n + 1}}{a_{0} + a_{1}\omega + \cdots + a_{m + 1}
\omega^{m + 1}}, \quad m \leq n \leq 8,
\label{eq:mappingFunction_WENO-RM}
\end{equation}
where
\begin{equation}\left\{
\begin{array}{l}
\begin{aligned}
& a_{i} = C_{n + 1}^{i} (-d_{s})^{n - i}, \quad i = 0,1,\cdots,m, \\
& a_{m + 1} = (1 - d_{s})^{n} - \sum\limits_{i = 0}^{m}a_{i}.
\end{aligned}
\end{array}
\right.
\label{eq:coeff_WENO-RM}
\end{equation}
And $m=2, n=6$ is recommended in \cite{WENO-RM260}, then
\begin{equation}
\big( g^{\mathrm{RM}} \big)_{s}(\omega) = d_{s} + \dfrac{(\omega - 
d_{s})^{7}}{a_{0} + a_{1}\omega + a_{2}\omega^{2} + a_{3}\omega^{3}}
, \quad \omega \in [0,1]
\label{eq:mappingFunction_WENO-RM260}
\end{equation}
where
\begin{equation}
a_{0} = d_{s}^{6}, \quad 
a_{1} = -7d_{s}^{5},  \quad 
a_{2} = 21d_{s}^{4}, \quad 
a_{3} = (1 - d_{s})^{6} - \sum\limits_{i = 0}^{2}a_{i}.
\label{eq:coeffOf_WENO-RM260}
\end{equation}

$\blacksquare$ WENO-MAIM1 \cite{WENO-MAIMi}
\begin{equation}
\big( g^{\mathrm{MAIM}1} \big)_{s}\big( \omega \big) = d_{s} + 
\dfrac{A \bigg( \frac{1 + (-1)^{k}}{2} + \frac{1 + (-1)^{k + 1}}{2} 
\cdot \mathrm{sgm}\big( \omega - d_{s}, \delta, 1, k \big) \bigg) 
\cdot (\omega - d_{s})^{k + 1}}{A \bigg( \frac{1 + (-1)^{k}}{2} + 
\frac{1 + (-1)^{k + 1}}{2} \cdot \mathrm{sgm}\big( \omega - d_{s}, 
\delta, 1, k \big) \bigg) \cdot (\omega - d_{s})^{k} + 
\omega^{\frac{d_{s}}{m_{s}\omega + \epsilon_{\mathrm{A}}}} (1 - 
\omega)^{\frac{1 - d_{s}}{m_{s}(1 - \omega)+\epsilon_{\mathrm{A}}}}},
\label{eq:mapFuncMAIM1}
\end{equation}
with 
\begin{equation}
\mathrm{sgm}\big( x, \delta, B, k \big) = \left\{ 
\begin{array}{ll}
\begin{aligned}
&\dfrac{x}{|x|}, & |x| \geq \delta, \\ 
&\dfrac{x}{\Big(B\big( \delta ^2 - x^2 \big)\Big)^{k + 3} + \ |x|}, 
& |x| < \delta.
\end{aligned}
\end{array} \right.
\label{eq:sgm}
\end{equation} 
In Eq.\eqref{eq:mapFuncMAIM1}, $k\in\mathbb{N}^{+}, A>0, \delta >0$ 
with $\delta \rightarrow 0$, $\epsilon_{\mathrm{A}}$ is a very small 
positive number to prevent the denominator from becoming zero, and 
$m_{s}\in\Big[\frac{\alpha_{s}}{k+1}, M \big)$ with $M$ being a 
finite positive constant real number and $\alpha_{s}$ a positive 
constant that only depends on $s$ in the fifth-order WENO-MAIM1 
scheme. In Eq.\eqref{eq:sgm}, the positive parameter $B$ is a scale 
transformation factor introduced to adjust the shape of the mapping 
function and it is set to be $1$ in WENO-MAIM1 while to be other 
values in the following WENO-ACM schemes.

$\blacksquare$ WENO-ACM \cite{WENO-ACM}
\begin{equation}
\big( g^{\mathrm{ACM}} \big)_{s}(\omega) = \left\{
\begin{array}{ll}
\begin{aligned}
&\dfrac{d_{s}}{2}\mathrm{sgm}(\omega - \mathrm{CFS}_{s}, \delta_{s}, 
B, k) + \dfrac{d_{s}}{2}, & \omega \leq d_{s}, \\
&\dfrac{1-d_{s}}{2}\mathrm{sgm}(\omega -\overline{\mathrm{CFS}}_{s}, 
\delta_{s}, B, k) + \dfrac{1 + d_{s}}{2}, & \omega > d_{s},
\end{aligned}
\end{array}
\right.
\label{eq:mapFuncACM}
\end{equation}
where $\mathrm{CFS}_{s} \in (0,d_{s})$, $\overline{\mathrm{CFS}}_{s} 
=1- \frac{1-d_{s}}{d_{s}} \times \mathrm{CFS}_{s}$ with 
$\overline{\mathrm{CFS}}_{s} \in (d_{s}, 1)$, and $\delta_{s} < \min
\Big\{\mathrm{CFS}_{s}, d_{s} - \mathrm{CFS}_{s}, (1 - d_{s})\Big(1 
- \frac{\mathrm{CFS}_{s}}{d_{s}}\Big) , \frac{1-d_{s}}{d_{s}}
\mathrm{CFS}_{s}\Big\}$.

$\blacksquare$ MIP-WENO-ACM$k$ \cite{MOP-WENO-ACMk}
\begin{equation}
\big( g^{\mathrm{MIP-ACM}k} \big)_{s}(\omega) = \left\{
\begin{array}{ll}
k_{s} \omega, & \omega \in [0, \mathrm{CFS}_{s}), \\
d_{s}, & \omega \in [\mathrm{CFS}_{s},\overline{\mathrm{CFS}}_{s}],\\
1 - k_{s} (1 - \omega), & \omega \in (\overline{\mathrm{CFS}}_{s},1],
\end{array}
\right.
\label{eq:mappingFunctionMIP-ACMk}
\end{equation}
where $\mathrm{CFS}_{s}\in (0,d_{s})$, $\overline{\mathrm{CFS}}_{s} =
1- \frac{1-d_{s}}{d_{s}}\times\mathrm{CFS}_{s}$ with 
$\overline{\mathrm{CFS}}_{s} \in (d_{s}, 1)$, and $k_{s} \in 
\Big[0, \frac{d_{s}}{\mathrm{CFS}_{s}}\Big]$. 

By using the mapping function $\big(g^{\mathrm{X}}\big)_{s}(\omega)$
, where the superscript ``X'' corresponds to ``M'', ``PM6'', or 
``IM'', etc., the nonlinear weights of the associated WENO-X scheme 
are defined as
\begin{equation*}
\omega_{s}^{\mathrm{X}} = \dfrac{\alpha _{s}^{\mathrm{X}}}{\sum_{l = 
0}^{2} \alpha _{l}^{\mathrm{X}}}, \alpha_{s}^{\mathrm{X}} = \big( g^{
\mathrm{X}} \big)_{s}(\omega^{\mathrm{JS}}_{s}), \quad s = 0,1,2,
\end{equation*}
where $\omega_{s}^{\mathrm{JS}}$ are calculated by 
Eq.(\ref{eq:weights:WENO-JS}).

In references, it has been analyzed and proved in detail that the
WENO-X schemes can retain the optimal order of accuracy in smooth 
regions even at or near critical points.


%% file: article_GMG_WENO_X.tex
\section{A general method to introduce order-preserving mapping for 
mapped WENO schemes}
\label{MOP-WENO-X}

\subsection{The OP-Mapped WENO}
Before giving Definition \ref{def:OP-Mapped_WENO} below, to maintain 
coherence and for the readers' convenience, we state the 
definition of \textit{order-preserving/non-order-preserving} mapping 
and \textit{OP/non-OP} point proposed in \cite{MOP-WENO-ACMk}.

\begin{definition}(order-preserving/non-order-preserving mapping) 
Suppose that $\big( g^{\mathrm{X}} \big)_{s}(\omega), s = 0,\cdots,
r-1$ is a monotone increasing piecewise mapping function of the 
$(2r-1)$th-order mapped WENO-X scheme. If for $\forall m, n \in \{0, 
\cdots,r-1\}$, when $\omega_{m} > \omega_{n}$, we have
\begin{equation}
\big( g^{\mathrm{X}}\big)_{m}(\omega_{m}) \geq \big( 
g^{\mathrm{X}} \big)_{n}(\omega_{n}).
\label{def:order-preserving_mappng}
\end{equation}
and when $\omega_{m} = \omega_{n}$, we have $\big( g^{\mathrm{X}}
\big)_{m}(\omega_{m}) = \big( g^{\mathrm{X}} \big)_{n}(\omega_{n})$, 
then we say the set of mapping functions \Big\{$\big( g^{\mathrm{X}}
\big)_{s}(\omega), s = 0,\cdots,r-1$\Big\} is 
\textbf{order-preserving (OP)}. Otherwise, we say the set of mapping 
functions \Big\{$\big( g^{\mathrm{X}}\big)_{s}(\omega), s = 0,\cdots,
r-1\Big\}$ is \textbf{non-order-preserving (non-OP)}.
\label{def:OPM}
\end{definition}

\begin{definition}(OP/non-OP point) Let $S^{2r-1}$ denote the 
$(2r-1)$-point global stencil centered around $x_{j}$. Assume that 
$S^{2r-1}$ is subdivided into $r$-point substencils $\{S_{0},\cdots,
S_{r-1}\}$ and $\omega_{s}$ are the nonlinear weights corresponding 
to the substencils $S_{s}$ with $s=0,\cdots,r-1$, which are used as 
the independent variables by the mapping function. Suppose that 
$\big( g^{\mathrm{X}}\big)_{s}(\omega), s=0,\cdots,r-1$ is the 
mapping function of the mapped WENO-X scheme, then we say that a 
\textbf{non-OP} mapping process occurs at $x_{j}$, if $\exists m, n 
\in \{0,\cdots,r-1\}$, s.t.
\begin{equation}\left\{
\begin{array}{ll}
\begin{aligned}
&\big(\omega_{m} - \omega_{n}\big)\bigg(\big(g^{\mathrm{X}}\big)_{m}
(\omega_{m}) - \big(g^{\mathrm{X}}\big)_{n}(\omega_{n})\bigg) < 0, 
&\mathrm{if} \quad \omega_{m} \neq \omega_{n},\\
&\big(g^{\mathrm{X}}\big)_{m}(\omega_{m}) \neq \big(g^{\mathrm{X}}
\big)_{n}(\omega_{n}), &\mathrm{if} \quad \omega_{m}=\omega_{n}.
\end{aligned}
\end{array}\right.
\end{equation}
And we say $x_{j}$ is a \textbf{non-OP point}. Otherwise, we say 
$x_{j}$ is an \textbf{OP point}.
\label{def:MaRe}
\end{definition}

\begin{definition}(OP-Mapped WENO)
The family of the mapped WENO schemes with \textit{OP} mappings is 
collectively referred to as \textbf{OP-Mapped WENO} in our study.
\label{def:OP-Mapped_WENO}
\end{definition}

\subsection{A general formula for the existing mapping functions}
We rewrite the mapping function of the WENO-X scheme, that is, 
$\big( g^{\mathrm{X}} \big)_{s}(\omega), s=0,1,\cdots,r-1$, to be a 
general formula, given as
\begin{equation}
g^{\mathrm{X}}\big(\omega; m_{\mathrm{P}},P_{s,1},\cdots,
P_{s,m_{\mathrm{P}}}\big) = \big( g^{\mathrm{X}} \big)_{s}(\omega),
\label{eq:generalFormulation}
\end{equation}
where $m_{\mathrm{P}}$ is the number of the parameters related with 
$s$ indicating the substencil, and $P_{s,1},\cdots,P_{s,m_{
\mathrm{P}}}$ are these parameters. Clearly, we have $m_{\mathrm{P}} 
= 0$ for the WENO-JS scheme and $m_{\mathrm{P}} \geq 1$ for other 
mapped WENO schemes. In Table \ref{table:GF_parameters}, taking $9$ 
different WENO schemes as examples, we have presented their 
parameters of $m_{\mathrm{P}}$ and $P_{s,1},\cdots,P_{s,m_{\mathrm{P}
}}$. Let $n_{\mathrm{X}}$ denote the order of the specified critical 
point, namely $\omega = d_{s}$, of the mapping function of the 
WENO-X scheme, that is, $\big( g^{\mathrm{X}} \big)'_{s}(d_{s}) = 
\cdots = \big( g^{\mathrm{X}} \big)^{(n_{\mathrm{X}})}_{s}(d_{s})=0, 
\big( g^{\mathrm{X}} \big)^{(n_{\mathrm{X}} + 1)}_{s}(d_{s}) \neq 0$
. To simplify the description of Theorem \ref{theorem:MOP-mapping} 
below, we present $n_{\mathrm{X}}$ of the WENO-X scheme in the sixth 
column of Table \ref{table:GF_parameters}. 

\begin{table}[ht]
\footnotesize
\centering
\caption{The parameters $m_{\mathrm{P}}$ and $P_{s,1},\cdots,P_{s,m_{
\mathrm{P}}}$ for the WENO-JS scheme and some existing mapped WENO 
schemes whose mapping functions are \textit{non-OP}.}
\label{table:GF_parameters}
\begin{tabular*}{\hsize}
{@{}@{\extracolsep{\fill}}lllllll@{}}
\toprule
No.&Scheme, WENO-X & $m_{\mathrm{P}}$ & $P_{s,1},\cdots,
P_{s,m_{\mathrm{P}}}$ &   Parameters  & $n_{\mathrm{X}}$ & Ref. \\
\hline
1   & WENO-JS        & $0$              & None             &  None
& None  & See \cite{WENO-JS}    \\
\hdashline
2   & WENO-M         & $1$              & $P_{s,1} = d_{s}$ & None 
& 2 &   See \cite{WENO-M}       \\
3   & WENO-IM($k,A$) & $1$ & $P_{s,1} = d_{s}$ & $k=2.0, A = 0.1$  
& $k$  & See \cite{WENO-IM}    \\ 
4   & WENO-PM$k$     & $1$          & $P_{s,1} = d_{s}$    & $k = 6$
& $k$  & See \cite{WENO-PM}     \\
5   & WENO-PPM$n$    & $1$         & $P_{s,1} = d_{s}$    & $n = 5$  
& $4$  & See \cite{WENO-PPM5} \\
6   & WENO-RM($mn0$) & $1$ & $P_{s,1} = d_{s}$    & $m = 2, n = 6$ 
& $3, 4$  & See \cite{WENO-RM260} \\
7   & WENO-MAIM$1$   & $2$   & $P_{s,1} = d_{s}, P_{s,2} = m_{s}$
& $k = 10, A = 1.0\mathrm{e-}6, m_{s} = 0.06$ 
& $k, k + 1$ & See \cite{WENO-MAIMi} \\
8  & WENO-ACM     & $2$ & $P_{s,1}=d_{s}, P_{s,2}=\mathrm{CFS}_{s}$ 
& $A = 20, k = 2, \mu = 1\mathrm{e-}6, \mathrm{CFS}_{s} = d_{s}/10$ 
& $\infty$ & See \cite{WENO-ACM} \\
9  & MIP-WENO-ACM$k$ & $3$& $P_{s,1} = d_{s}, P_{s,2} = 
\mathrm{CFS}_{s}, P_{s,3} = k_{s}$ & $k_{s} = 0.0, 
\mathrm{CFS}_{s} = d_{s}/10$  & $\infty$& See \cite{MOP-WENO-ACMk} \\
\bottomrule
\end{tabular*}
\end{table}

\begin{lemma}
For the WENO-X scheme shown in Table \ref{table:GF_parameters}, the 
mapping function $\big( g^{\mathrm{X}} \big)_{s}(\omega), s=0,1,
\cdots,r-1$ is monotonically increasing over $[0, 1]$.
\label{lemma:general_mapping}
\end{lemma}
\textbf{Proof.}
See the corresponding references given in the last column of Table 
\ref{table:GF_parameters}.
$\hfill\square$ \\

\subsection{The new mapping functions}\label{subsec:newMapping}
Firstly, we give the \textbf{minDist} function by the following 
definition.

\begin{definition}({\rm{\textbf{minDist}}} function)
Define the {\rm{\textbf{minDist}}} function as follows
\begin{subnumcases}{}
\bm{\mathrm{minDist}}\big(x_{0},\cdots,x_{r-1};d_{0},\cdots,
d_{r-1};\omega\big) = x_{k^{*}}, 
\label{eq:minDist:01} \\
k^{*} = \min\Bigg(\mathrm{IndexOf}\bigg(\min\Big\{\lvert \omega-d_{0}
\rvert, \lvert \omega-d_{1}\rvert, \cdots, \lvert \omega-d_{r-1}
\rvert\Big\}\bigg)\Bigg),
\label{eq:minDist:02}
\end{subnumcases}
where $d_{s}, s=0,\cdots,r-1$ is the optimal weight, $\omega$ is the 
nonlinear weight being the independent variable of the mapping 
function, and the function $\mathrm{IndexOf}(\cdot)$ returns a set 
of the subscripts of ``$\cdot$'', that is, if $\min\Big\{\lvert 
\omega-d_{0}\rvert, \lvert \omega-d_{1}\rvert, \cdots, \lvert \omega
- d_{r-1}\rvert\Big\} = \lvert \omega - d_{m_{1}} \rvert = \lvert 
\omega - d_{m_{2}} \rvert = \cdots = \lvert \omega - d_{m_{M}} 
\rvert$, then
\begin{equation}
\mathrm{IndexOf}\bigg(\min\Big\{\lvert \omega-d_{0}
\rvert, \lvert \omega-d_{1}\rvert, \cdots, \lvert \omega-d_{r-1}
\rvert\Big\}\bigg) = \Big\{m_{1},m_{2},\cdots,m_{M}\Big\}.
\label{eq:IndexOf}
\end{equation}
\label{definition:minDist}
\end{definition}

Let $\mathcal{D} = \Big\{d_{0},d_{1},\cdots,d_{r-1}\Big\}$ be an 
array of all the ideal weights of the $(2r-1)$th-order WENO schemes. 
We build a new array by sorting the elements of $\mathcal{D}$ in 
ascending order, that is, $\mathcal{\widetilde{D}} = 
\Big\{\widetilde{d}_{0}, \widetilde{d}_{1},\cdots,\widetilde{d}_{r-1}
\Big\}$. In other words, the arrays $\mathcal{D}$ and 
$\mathcal{\widetilde{D}}$ have the same elements with different 
arrangements, and the elements of $\mathcal{\widetilde{D}}$ satisfy
\begin{equation}
0 < \widetilde{d}_{0} < \widetilde{d}_{1} < \cdots < 
\widetilde{d}_{r-1} < 1.
\label{eq:tilde_ds}
\end{equation}

\begin{definition}
Let $\mathcal{G} = \Big\{ \big( g^{\mathrm{X}} \big)_{0}(\omega), 
\big( g^{\mathrm{X}} \big)_{1}(\omega), \cdots, \big( 
g^{\mathrm{X}} \big)_{r-1}(\omega) \Big\}$ be an array of all the 
mapping functions of the $(2r-1)$th-order mapped WENO-X scheme. We 
define a new array by sorting the elements of $\mathcal{G}$ in a 
new order, that is, $\mathcal{\widetilde{G}}=\Big\{\widetilde{
\big(g^{\mathrm{X}}\big)}_{0}(\omega),\widetilde{\big(g^{\mathrm{X}}
\big)}_{1}(\omega), \cdots, \widetilde{\big(g^{\mathrm{X}}\big)}_{
r-1}(\omega) \Big\}$, where $\widetilde{\big(g^{\mathrm{X}}\big)}_{s}
(\omega)$ is the mapping function associated with $\widetilde{d}_{s}$
.
\label{def:newMappingSet}
\end{definition}

\begin{lemma}
Denote $\widetilde{d}_{-1} = 0, \widetilde{d}_{r} = 1$. Let 
$\mathring{d}_{-1} = \widetilde{d}_{-1}, \mathring{d}_{0} = 
\frac{\widetilde{d}_{0} + \widetilde{d}_{1}}{2}, \cdots,
\mathring{d}_{r-2}=\frac{\widetilde{d}_{r-2}+\widetilde{d}_{r-1}}{2}
, \mathring{d}_{r-1} = \widetilde{d}_{r}$. For $\forall i = 0,1,
\cdots,r-1$, if $\omega \in (\mathring{d}_{i-1}, \mathring{d}_{i}]$, 
then 
\begin{equation*}
\min\Bigg(\mathrm{IndexOf}\bigg(\min\Big\{\lvert \omega-
\widetilde{d}_{0}\rvert, \lvert \omega-\widetilde{d}_{1}\rvert, 
\cdots, \lvert \omega-\widetilde{d}_{r-1}\rvert\Big\}\bigg)\Bigg)= i.
\end{equation*}
\label{lemma:Omega_i:minDist}
\end{lemma}
\textbf{Proof.}\\

(1) We first prove the cases of $i=1,\cdots,r-2$. When 
$\widetilde{d}_{i} \leq \omega \leq \frac{\widetilde{d}_{i}+
\widetilde{d}_{i+1}}{2}$, as Eq.\eqref{eq:tilde_ds} holds, we get
\begin{equation}\left\{
\begin{array}{l}
0 \leq \omega - \widetilde{d}_{i} \leq \widetilde{d}_{i+1} - \omega 
< \cdots < \widetilde{d}_{r-1} - \omega,\\
0 \leq \omega - \widetilde{d}_{i} < \omega -\widetilde{d}_{i-1}
< \cdots < \omega - \widetilde{d}_{0}.
\end{array}
\right.
\label{eq:prf:Lemma:minDist:01}
\end{equation}
Similarly, when $\frac{\widetilde{d}_{i-1}+\widetilde{d}_{i}}{2} < 
\omega < \widetilde{d}_{i}$, we get
\begin{equation}\left\{
\begin{array}{l}
0 < \widetilde{d}_{i} - \omega < \omega - \widetilde{d}_{i-1} < 
\cdots < \omega - \widetilde{d}_{0},\\
0 < \widetilde{d}_{i} - \omega < \widetilde{d}_{i+1} - \omega < 
\cdots < \widetilde{d}_{r-1} - \omega.
\end{array}
\right.
\label{eq:prf:Lemma:minDist:02}
\end{equation}
Then, according to Eq.\eqref{eq:prf:Lemma:minDist:01}\eqref{eq:prf:Lemma:minDist:02}, we obtain
\begin{equation}
\min\Big\{\lvert \omega - \widetilde{d}_{0} \rvert, \cdots, \lvert 
\omega - \widetilde{d}_{i-1} \rvert, \lvert \omega - 
\widetilde{d}_{i} \rvert,\lvert \omega - 
\widetilde{d}_{i+1} \rvert,\cdots,\lvert \omega - \widetilde{d}_{r-1}
\rvert\Big\} = \lvert \omega - \widetilde{d}_{i} \rvert = \lvert 
\omega - \widetilde{d}_{i+1} \rvert, 
\quad i = 1,\cdots,r-2,
\label{eq:prf:Lemma:minDist:03}
\end{equation}
where the last equality holds if and only if $\omega - 
\widetilde{d}_{i} = \widetilde{d}_{i+1} - \omega$.

(2) For the case of $i=0$, we know that $\omega \in 
(\mathring{d}_{-1},\mathring{d}_{0}] =\Big(0,\frac{\widetilde{d}_{0}+
\widetilde{d}_{1}}{2}\Big]$. When $\widetilde{d}_{0} \leq \omega \leq
\frac{\widetilde{d}_{0}+\widetilde{d}_{1}}{2}$, we have
\begin{equation}
0 \leq \omega - \widetilde{d}_{0} \leq \widetilde{d}_{1} - \omega < 
\cdots < \widetilde{d}_{r-1} - \omega.
\label{eq:prf:Lemma:minDist:04}
\end{equation}
And when $0 < \omega < \widetilde{d}_{0}$, we have
\begin{equation}
0 < \widetilde{d}_{0} - \omega < \widetilde{d}_{1} - \omega < 
\cdots < \widetilde{d}_{r-1} - \omega.
\label{eq:prf:Lemma:minDist:05}
\end{equation}
Then, according to Eq.\eqref{eq:prf:Lemma:minDist:04}\eqref{eq:prf:Lemma:minDist:05}, we 
obtain
\begin{equation}
\min\Big\{\lvert \omega - \widetilde{d}_{0} \rvert, \cdots, \lvert 
\omega - \widetilde{d}_{i-1} \rvert, \lvert \omega - 
\widetilde{d}_{i} \rvert,\lvert \omega - 
\widetilde{d}_{i+1} \rvert,\cdots,\lvert \omega - \widetilde{d}_{r-1}
\rvert\Big\} = \lvert \omega - \widetilde{d}_{0} \rvert = \lvert 
\omega - \widetilde{d}_{1} \rvert,
\label{eq:prf:Lemma:minDist:06}
\end{equation}
where the last equality holds if and only if $\omega - 
\widetilde{d}_{0} = \widetilde{d}_{1} - \omega$.

(3) As the proof of the case of $i=r-1$ is very similar to that of 
the case $i=0$, we do not state it here for simplicity. And we can 
get that, if $\omega \in (\mathring{d}_{r-2}, \mathring{d}_{r-1}]$, 
then
\begin{equation}
\min\Big\{\lvert \omega - \widetilde{d}_{0} \rvert, \cdots, \lvert 
\omega - \widetilde{d}_{i-1} \rvert, \lvert \omega - 
\widetilde{d}_{i} \rvert,\lvert \omega - 
\widetilde{d}_{i+1} \rvert,\cdots,\lvert \omega - \widetilde{d}_{r-1}
\rvert\Big\} = \lvert \omega - \widetilde{d}_{r-1} \rvert.
\label{eq:prf:Lemma:minDist:07}
\end{equation}

(4) Thus, according to Eq.\eqref{definition:minDist} and
Eq.\eqref{eq:prf:Lemma:minDist:03}\eqref{eq:prf:Lemma:minDist:06}\eqref{eq:prf:Lemma:minDist:07}, we obtain
\begin{equation*}
\min\Bigg(\mathrm{IndexOf}\bigg(\min\Big\{\lvert \omega - 
\widetilde{d}_{0} \rvert, \cdots, \lvert \omega - \widetilde{d}_{i-1}
\rvert, \lvert \omega - \widetilde{d}_{i} \rvert,\lvert \omega - 
\widetilde{d}_{i+1} \rvert,\cdots,\lvert \omega - \widetilde{d}_{r-1}
\rvert\Big\}\bigg)\Bigg)=i.
\end{equation*}

Up to now, we have finished the proof of Lemma 
\ref{lemma:Omega_i:minDist}.
$\hfill\square$ \\

For simplicity of description and according to Lemma 
\ref{lemma:Omega_i:minDist}, we introduce intervals $\Omega_{i}$ 
defined as follows
\begin{equation}
\Omega_{i}=\Big\{\omega \lvert \mathrm{\textbf{minDist}}(
\widetilde{d}_{0}, \widetilde{d}_{1}, \cdots, \widetilde{d}_{r-1}; 
\widetilde{d}_{0}, \widetilde{d}_{1}, \cdots, \widetilde{d}_{r-1}; 
\omega) = \widetilde{d}_{i} \Big\} = (\mathring{d}_{i-1}, 
\mathring{d}_{i}], \quad i = 0, 1, \cdots, r-1. 
\label{def:Omega_i}
\end{equation}
If $\omega \in \Omega = (0,1]$, it is trivial to verify that: (1) 
$\Omega = \Omega_{0} \bigcup \Omega_{1} \bigcup \cdots \bigcup 
\Omega_{r-1}$; (2) for $\forall i,j = 0,1,\cdots,r-1$ and $i\neq j$, 
$\Omega_{i}\bigcap\Omega_{j} = \varnothing$.

\begin{figure}[ht]
\centering
  \includegraphics[height=0.285\textwidth]
  {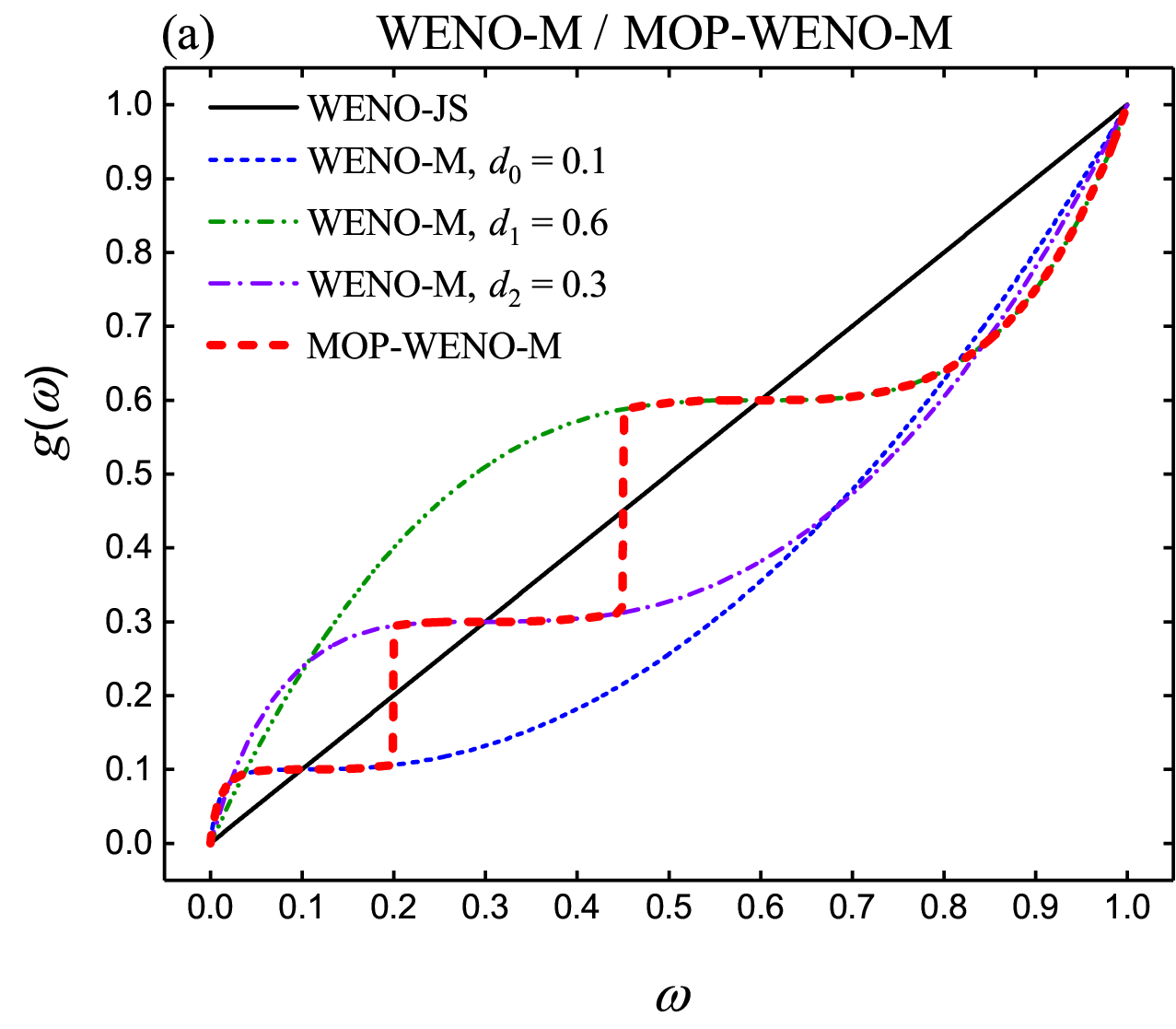}
  \includegraphics[height=0.285\textwidth]
  {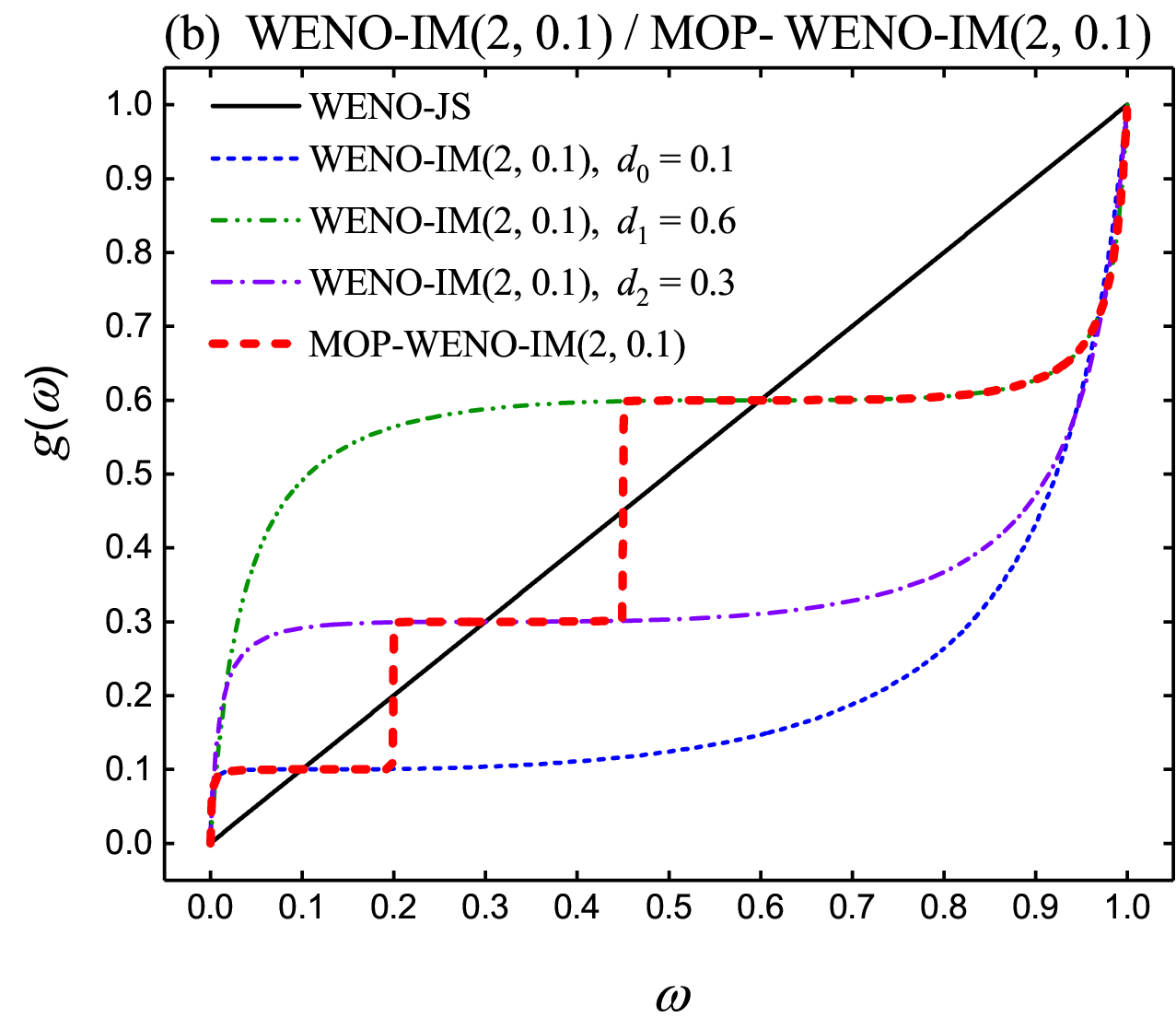}
  \includegraphics[height=0.285\textwidth]
  {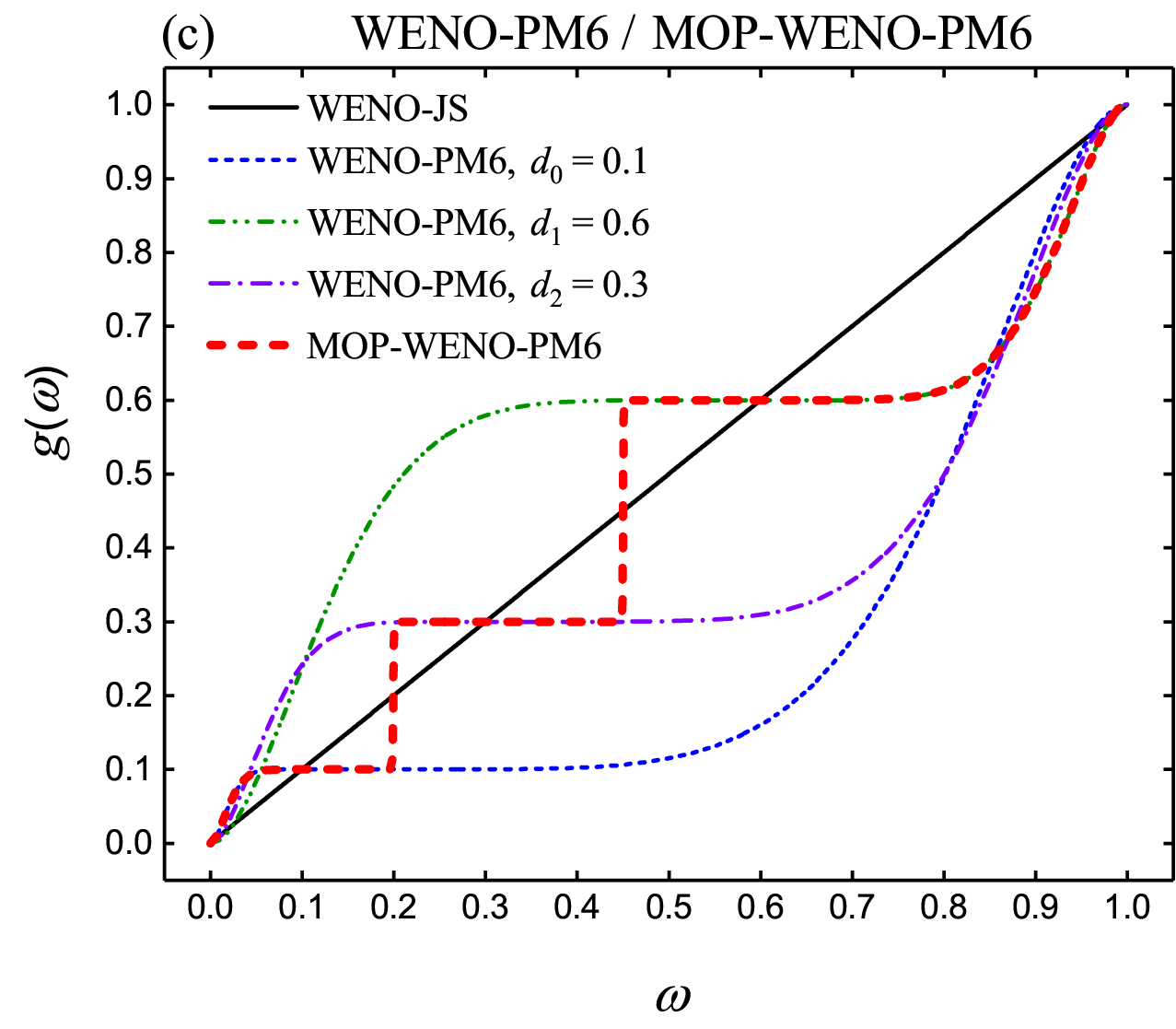}\\  
  \includegraphics[height=0.285\textwidth]
  {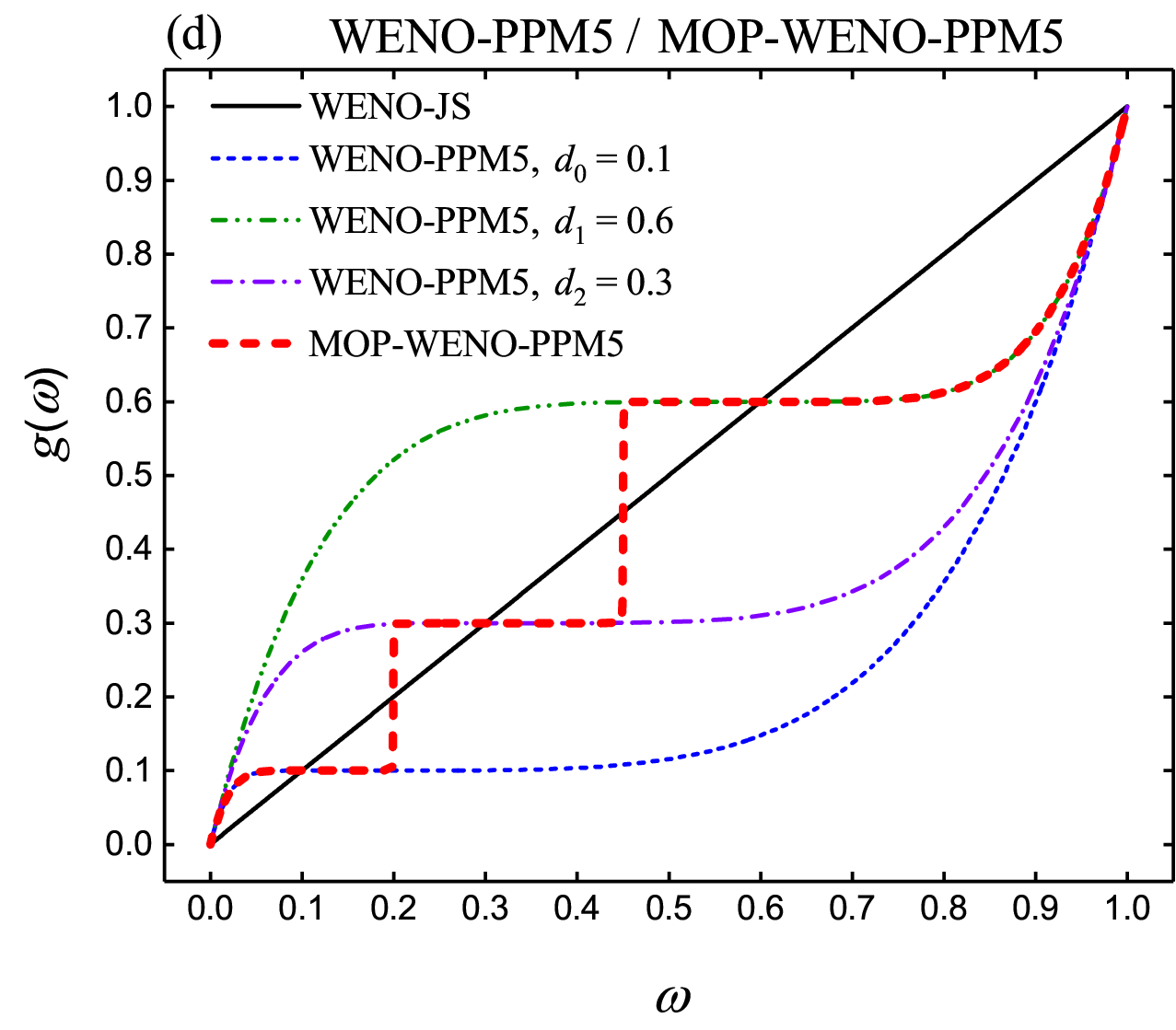}
  \includegraphics[height=0.285\textwidth]
  {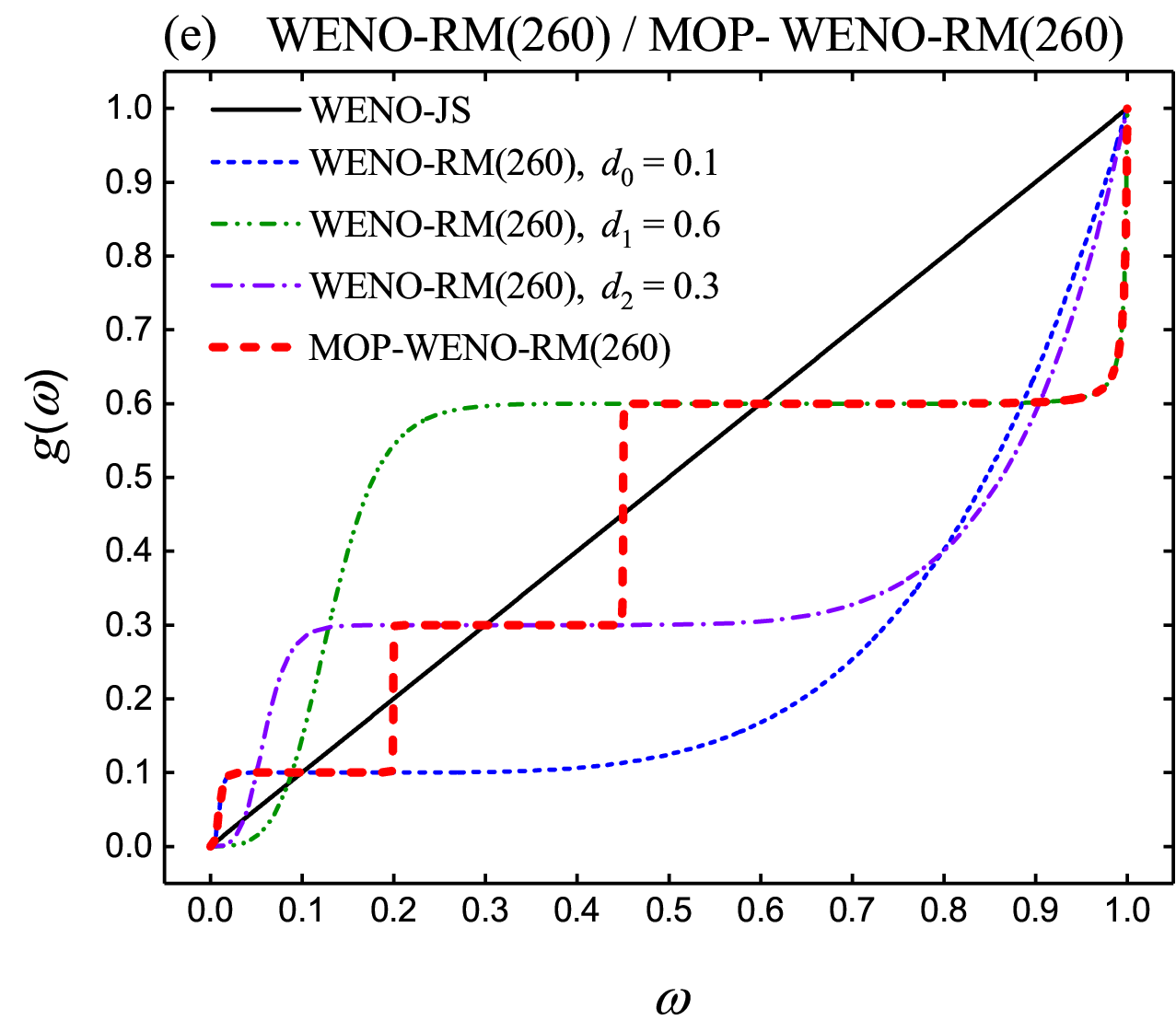}
  \includegraphics[height=0.285\textwidth]
  {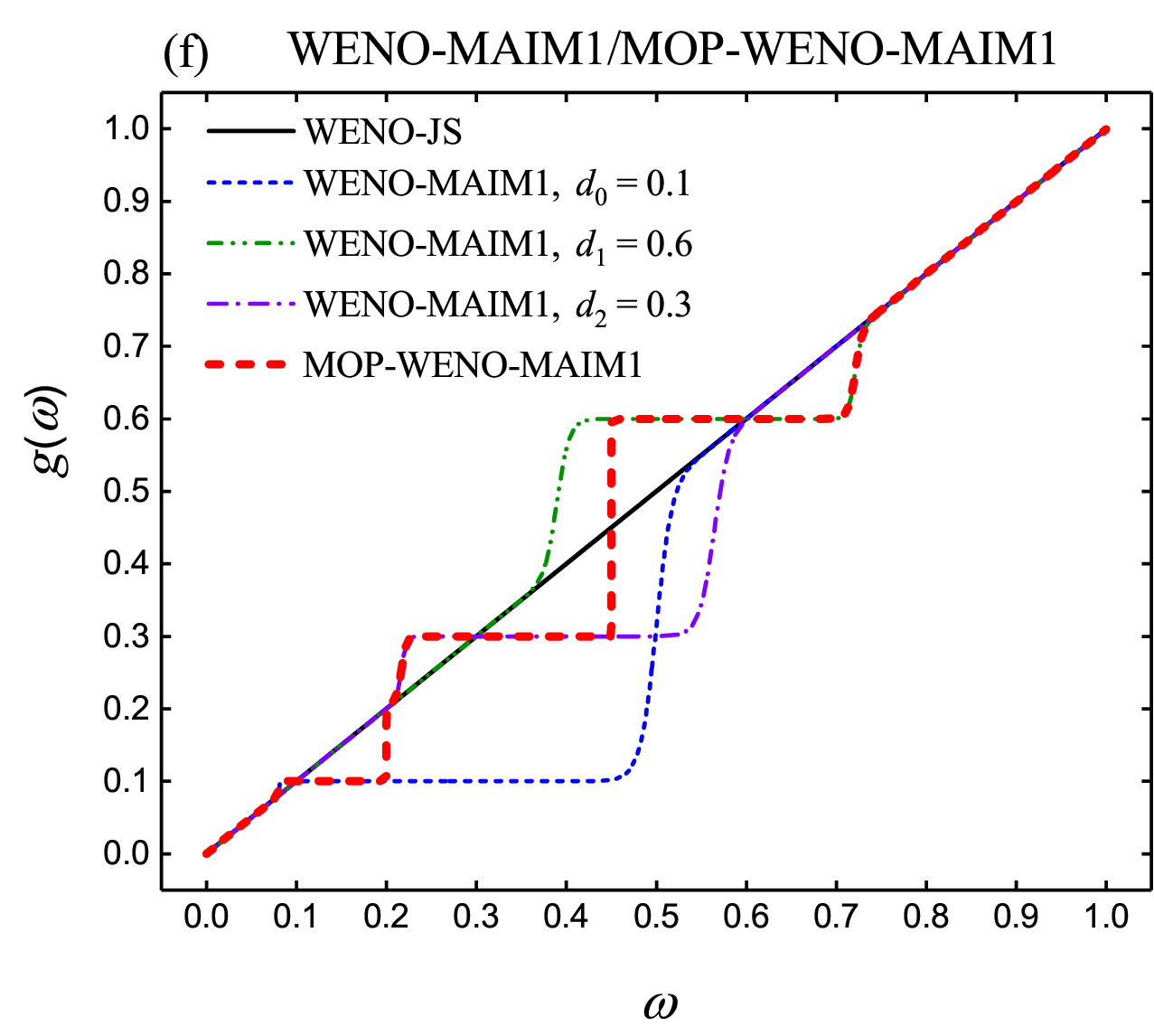}\\
  \includegraphics[height=0.285\textwidth]
  {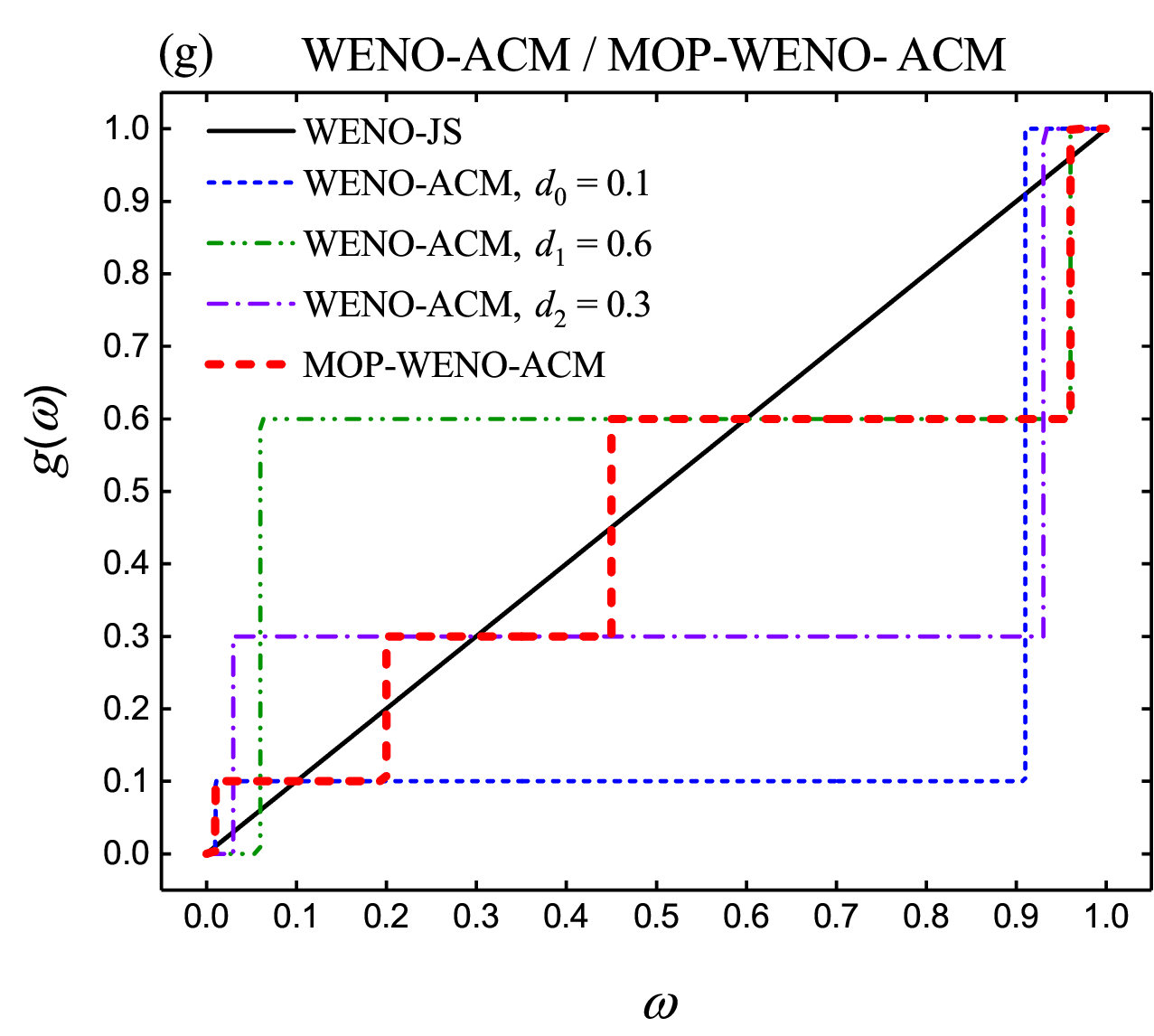}   
  \includegraphics[height=0.285\textwidth]
  {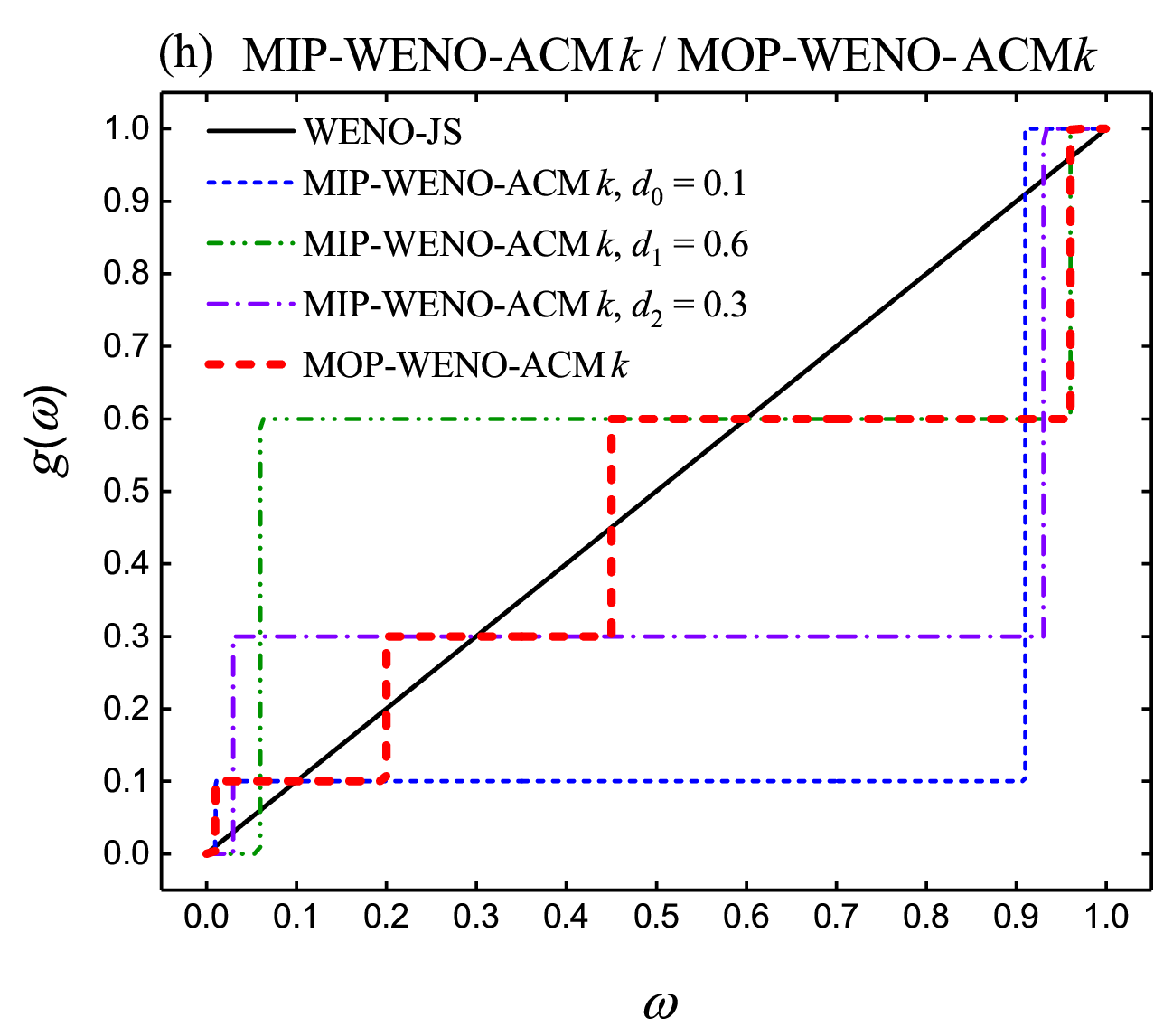} 
  \includegraphics[height=0.285\textwidth]
  {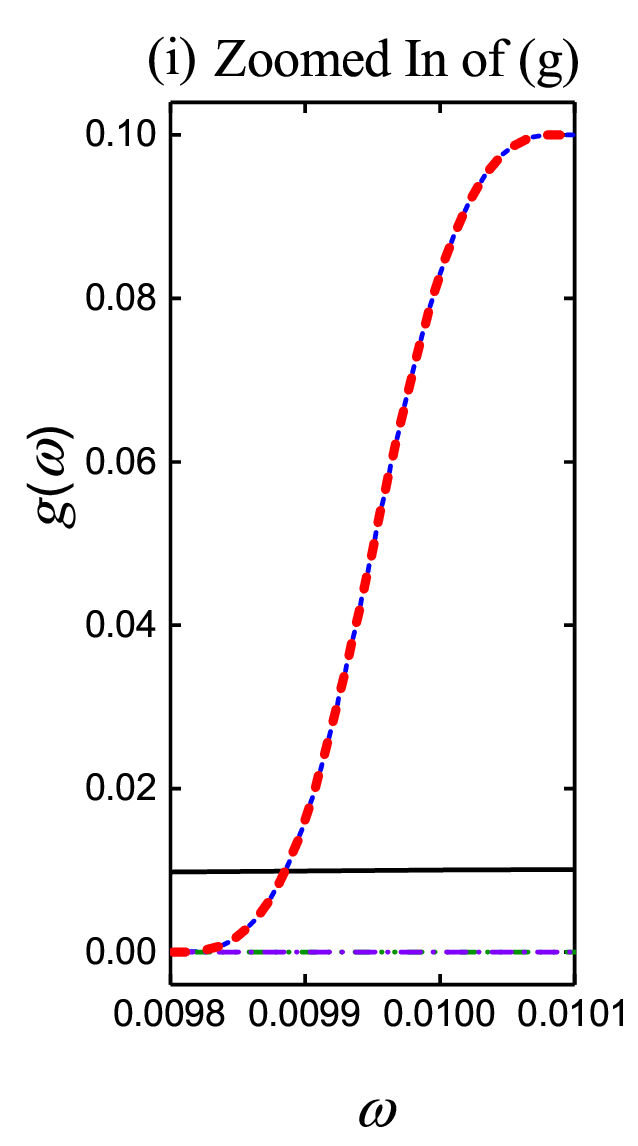}     
  \includegraphics[height=0.285\textwidth]
  {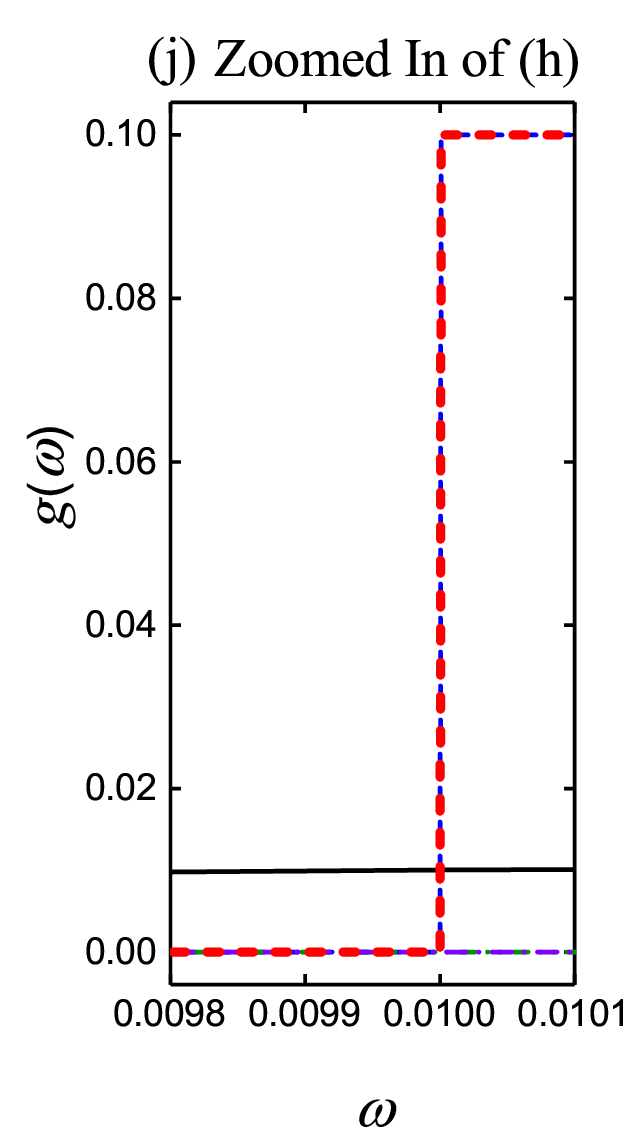}       
  \caption{Comparison for the mapping functions for WENO-X (shown in 
  Table \ref{table:GF_parameters}) and MOP-WENO-X.}
\label{fig:gOmega}
\end{figure}

\begin{lemma}
Let $a, b\in\{0,1,\cdots,r-1\}$ and WENO-X the scheme shown in Table 
\ref{table:GF_parameters}, for $\forall a \geq b$ and 
$\omega_{\alpha} \in \Omega_{a}, \omega_{\beta} \in \Omega_{b}$, we 
have the following properties: C1. if $a = b$ and $\omega_{\alpha} > 
\omega_{\beta}$, then $\widetilde{\big(g^{\mathrm{X}}\big)}_{a}
(\omega_{\alpha}) \geq \widetilde{\big(g^{\mathrm{X}}\big)}_{b}
(\omega_{\beta})$; C2. if $a = b$ and $\omega_{\alpha} = 
\omega_{\beta}$, then $\widetilde{\big(g^{\mathrm{X}}\big)}_{a}
(\omega_{\alpha})=\widetilde{\big(g^{\mathrm{X}}\big)}_{b}
(\omega_{\beta})$; C3. if $a > b$, then $\omega_{\alpha} > 
\omega_{\beta}, \widetilde{\big(g^{\mathrm{X}}\big)}_{a}
(\omega_{\alpha}) > \widetilde{\big(g^{\mathrm{X}}\big)}_{b}
(\omega_{\beta})$.
\label{lemma:Omega_i}
\end{lemma}
\textbf{Proof.}
(1) We can directly get properties \textit{C1} and \textit{C2} from 
Lemma \ref{lemma:general_mapping}. (2) As $a > b$, according to 
Eq.\eqref{eq:tilde_ds}\eqref{def:Omega_i}, we know that the 
interval $\Omega_{a}$ must be on the right side of the interval 
$\Omega_{b}$, while $\omega_{\alpha} \in \Omega_{a}, \omega_{\beta} 
\in \Omega_{b}$ is given, then we get $\omega_{\alpha} > 
\omega_{\beta}$. Trivially, according to Definition 
\ref{def:newMappingSet}, or by intuitively observing the curves of 
the mapping function $\widetilde{\big(g^{\mathrm{X}}\big)}_{s}
(\omega)$ as shown in Fig. \ref{fig:gOmega}, we can obtain 
$\widetilde{\big(g^{\mathrm{X}}\big)}_{a}(\omega_{\alpha}) > 
\widetilde{\big(g^{\mathrm{X}}\big)}_{b}(\omega_{\beta})$. Thus, 
\textit{C3} is proved.
$\hfill\square$ \\

By employing the \textbf{minDist} function, we build a general 
method to introduce the \textit{OP} criterion into the existing 
mappings which are \textit{non-OP}. The general method is stated in 
Algorithm \ref{alg:minDist}.

\begin{algorithm}[htb]
\caption{A general method to construct \textit{OP} mappings.}
\label{alg:minDist}
\SetKwInOut{Input}{input}\SetKwInOut{Output}{output}
\Input{$s$, index indicating the substencil $S_{s}$ and $s = 0,1,
\cdots, r-1$ \\
$d_{s}$, optimal weights \\
$\omega^{\mathrm{JS}}_{s}$, nonlinear weights computed by the 
WENO-JS scheme \\
$m_{\mathrm{P}}$, the number of the parameters related with $s$ \\
$P_{s,j}$, parameters related with $s$ and $j = 1,\cdots, 
m_{\mathrm{P}}$}
\Output{$\Big\{\big(g^{\mathrm{MOP-X}}\big)_{s}(\omega^{
\mathrm{JS}}_{s}), s=0,1,\cdots,r-1\Big\}$, the new set of mapping 
functions that is \textit{OP}}
\BlankLine
\emph{$\big( g^{\mathrm{X}} \big)_{s}(\omega), s = 0, 1,\cdots,r-1$ 
is a monotonically increasing mapping function over $[0, 1]$, and 
the set of mapping functions $\Big\{\big(g^{\mathrm{X}} \big)_{s}(
\omega), s = 0, 1,\cdots,r-1 \Big\}$ is \textit{non-OP}}\;
\tcp{implementation of the ``minDist'' function in Definition 
\ref{definition:minDist}}
\For{$s=0; s \leq r - 1; s ++$}{	
	\tcp{get $k^{*}$ in Eq.(\ref{eq:minDist:02})}
	set $d^{\mathrm{min}}=\lvert\omega_{s}^{\mathrm{JS}}-d_{0}\rvert,
	k_{s}^{*} = 0$\;	
	\For{$i=1; i \leq r - 1; i ++$}{
		\If{$\lvert\omega_{s}^{\mathrm{JS}}-d_{i}\rvert<
		d^{\mathrm{min}}$}{
			$d^{\mathrm{min}}=\lvert \omega_{s}^{\mathrm{JS}} - d_{i}
			\rvert$, \\
			$k_{s}^{*} = i$\;
		}	
	}
	\tcp{remark: the for loop above indicates that 
	$\omega_{s}^{\mathrm{JS}} \in \Omega_{k_{s}^{*}}$}
	\tcp{get $x_{k^{*}}$ in Eq.(\ref{eq:minDist:01})}
	\For{$j = 1; j \leq m_{\mathrm{P}}; j ++$}{
		$\overline{P}_{s,j} = P_{k_{s}^{*},j}$\;
	}
}
\tcp{get $\big(g^{\mathrm{MOP-X}}\big)_{s}(\omega^{\mathrm{JS}}_{s})$
}
\For{$s=0; s \leq r - 1; s ++$}{
	$\big(g^{\mathrm{MOP-X}}\big)_{s}(\omega^{\mathrm{JS}}_{s}) = 
	g^{\mathrm{X}}\Big(\omega^{\mathrm{JS}}_{s}; m_{\mathrm{P}}, 
	\overline{P}_{s,1},\cdots, \overline{P}_{s,m_{\mathrm{P}}}\Big)$.
}
\end{algorithm}

\begin{theorem}
The set of mapping functions $\Big\{\big(g^{\mathrm{MOP-X}}\big)_{s}(
\omega^{\mathrm{JS}}_{s}), s=0,1,\cdots,r-1\Big\}$ obtained through 
Algorithm \ref{alg:minDist} is \textit{OP}.
\label{theorem:g_MOP}
\end{theorem}
\textbf{Proof.} 
Let $\omega_{m}^{\mathrm{JS}}, \omega_{n}^{\mathrm{JS}} \in [0,1]$ 
and $\forall m,n \in \{0,1,\cdots,r-1\}$. According to 
Algorithm \ref{alg:minDist} and without loss of generality, 
we can assume that $\omega_{m}^{\mathrm{JS}} \in \Omega_{k_{m}^{*}}, 
\omega_{n}^{\mathrm{JS}} \in \Omega_{k_{n}^{*}}$, and then we get
\begin{equation*}\left\{
\begin{array}{l}
\big(g^{\mathrm{MOP-X}}\big)_{m}(\omega^{\mathrm{JS}}_{m}) = 
g^{\mathrm{X}}\Big(\omega^{\mathrm{JS}}_{m}; m_{\mathrm{P}}, 
P_{{k_{m}^{*}},1}, \cdots,P_{{k_{m}^{*}},m_{\mathrm{P}}}\Big), \\
\big(g^{\mathrm{MOP-X}}\big)_{n}(\omega^{\mathrm{JS}}_{n}) = 
g^{\mathrm{X}}\Big(\omega^{\mathrm{JS}}_{n}; m_{\mathrm{P}}, 
P_{{k_{n}^{*}},1}, \cdots,P_{{k_{n}^{*}},m_{\mathrm{P}}}\Big).
\end{array}
\right.
\label{eq:gMOP_m-n}
\end{equation*}
It is easy to verify that 
\begin{equation*}\left\{
\begin{array}{l}
g^{\mathrm{X}}\Big(\omega^{\mathrm{JS}}_{m}; m_{\mathrm{P}}, 
P_{{k_{m}^{*}},1}, \cdots,P_{{k_{m}^{*}},m_{\mathrm{P}}}\Big) = 
\widetilde{\big(g^{\mathrm{X}}\big)}_{k_{m}^{*}}
(\omega_{m}^{\mathrm{JS}}), \\
g^{\mathrm{X}}\Big(\omega^{\mathrm{JS}}_{n}; m_{\mathrm{P}}, 
P_{{k_{n}^{*}},1}, \cdots,P_{{k_{n}^{*}},m_{\mathrm{P}}}\Big) = 
\widetilde{\big(g^{\mathrm{X}}\big)}_{k_{n}^{*}}
(\omega_{n}^{\mathrm{JS}}).
\end{array}
\right.
\end{equation*}
Therefore, according to Lemma \ref{lemma:Omega_i}, we can finish the 
proof trivially. 
$\hfill\square$ \\

We now define the modified weights which are \textit{OP} as follows
\begin{equation}
\omega_{s}^{\mathrm{MOP-X}} = \dfrac{\alpha_{s}^{\mathrm{MOP-X}}}{
\sum_{l=0}^{r-1}\alpha_{l}^{\mathrm{MOP-X}}}, \quad \alpha_{s}^{
\mathrm{MOP-X}} = \big(g^{\mathrm{MOP-X}}\big)_{s}(\omega^{
\mathrm{JS}}_{s}), \quad s = 0,\cdots,r-1,
\label{eq:MOP-mapping}
\end{equation}
where $\big(g^{\mathrm{MOP-X}}\big)_{s}
(\omega^{\mathrm{JS}}_{s})$ is obtained from Algorithm 
\ref{alg:minDist}. The associated scheme will be referred to as 
MOP-WENO-X.

The mapping functions of the WENO-X schemes presented in Table 
\ref{table:GF_parameters} and those of the corresponding MOP-WENO-X 
schemes are shown in Fig. \ref{fig:gOmega}. We can find that, for 
the mapping functions of the MOP-WENO-X schemes: (1) the 
monotonicity over the whole domain $(0,1)$ is maintained; (2) the 
differentiability is reduced and limited to the neighborhood of the 
optimal weights $d_{s}$; (3) the \textit{OP} property is obtained. 
We summarize these properties as follows.

\begin{theorem}
Let $\overline{\Omega}_{i} = \Big\{\omega \in \Omega_{i} \cap \omega 
\neq \partial \Omega_{i} \Big\}, i = 0,1,\cdots,r-1$. The mapping 
function $\big(g^{\mathrm{MOP-X}}\big)_{s}(\omega)$ obtained from 
Althorithm \ref{alg:minDist} satisfies the following properties:\\

C1. for $\forall \omega \in \overline{\Omega}_{i}, i = 0, 1, \cdots, 
r-1$, $\big(g^{\mathrm{MOP-X}}\big)'_{s}(\omega) \geq 0$;

C2. for $\forall \omega \in \Omega$, $0 \leq \big(g^{\mathrm{MOP-X}}
\big)_{s}(\omega) \leq 1$;

C3. for $\forall s \in \Big\{ 0, 1, \cdots, r - 1 \Big\}$, 
$\widetilde{d}_{s} \in \Omega_{s}$, and $\big(g^{\mathrm{MOP-X}}
\big)_{s}(\widetilde{d}_{s}) = \widetilde{d}_{s}, 
\big(g^{\mathrm{MOP-X}}\big)'_{s}(\widetilde{d}_{s}) = \cdots = 
\big(g^{\mathrm{MOP-X}}\big)^{(n_{\mathrm{X}})}_{s}
(\widetilde{d}_{s}) = 0$ where $n_{\mathrm{X}}$ is given in Table 
\ref{table:GF_parameters};

C4. $\big(g^{\mathrm{MOP-X}}\big)_{s}(0) = 0, \big(g^{\mathrm{MOP-X}}
\big)_{s}(1) = 1, \big(g^{\mathrm{MOP-X}} \big)'_{s}(0) = 
\big(g^{\mathrm{X}} \big)'_{s}(0), \big(g^{\mathrm{MOP-X}} \big)'_{s}
(1) = \big(g^{\mathrm{X}} \big)'_{s}(1)$;

C5. for $\forall m, n \in \Big\{0, \cdots,r-1 \Big\}$, if $\omega_{m}
> \omega_{n}$, then $\big( g^{\mathrm{MOP-X}}\big)_{m}(\omega_{m}) 
\geq \big( g^{\mathrm{MOP-X}} \big)_{n}(\omega_{n})$, and if 
$\omega_{m} = \omega_{n}$, then $\big( g^{\mathrm{MOP-X}}\big)_{m}
(\omega_{m}) = \big( g^{\mathrm{MOP-X}} \big)_{n}(\omega_{n})$.
\label{theorem:MOP-mapping}
\end{theorem}

\subsection{Convergence properties}
According to Theorem \ref{theorem:MOP-mapping}, we get the 
convergence properties for the $(2r - 1)$th-order MOP-WENO-X schemes 
as given in Theorem \ref{theorem:MOP-convergence}. The proof is 
almost identical to that of the corresponding WENO-X schemes in 
the references presented in Table \ref{table:GF_parameters}.

\begin{theorem}
The requirements for the $(2r-1)$th-order MOP-WENO-X schemes to 
achieve the optimal order of accuracy are identical to that of the 
corresponding $(2r-1)$th-order WENO-X schemes.
\label{theorem:MOP-convergence}
\end{theorem}

For the integrity of this paper and the benefit of the reader, we 
concisely express the following Corollaries of Theorem 
\ref{theorem:MOP-convergence}.

\begin{corollary}
If $n$ mapping is used in the $(2r - 1)$th-order MOP-WENO-M scheme, 
then for different values of $n_{\mathrm{cp}}$, the weights 
$\omega_{s}^{\mathrm{MOP-M}}$ in the $(2r - 1)$th-order MOP-WENO-M 
scheme satisfy
\begin{equation*}
\omega_{s}^{\mathrm{MOP-M}}-d_{s} = O\Big((\Delta x)^{3^{n}\times(r- 
1 - n_{\mathrm{cp}})} \Big), \quad r = 2, 3, \cdots, 9, \quad n_{
\mathrm{cp}} = 0, 1, \cdots, r - 1,
\end{equation*}
and the rate of convergence is
\begin{equation*}
r_{\mathrm{c}} = \left\{
\begin{array}{ll}
2r - 1,  &\mathrm{if} \quad n_{\mathrm{cp}} = 0, \cdots, \Bigg\lfloor
\dfrac{3^{n}-1}{3^{n}}r - 1 \Bigg\rfloor, \\
(3^{n}+1)(r - 1) - 3^{n}\times n_{\mathrm{cp}}, &\mathrm{if} \quad n_
{\mathrm{cp}} = \Bigg\lfloor \dfrac{3^{n}-1}{3^{n}}r - 1 \Bigg\rfloor
+ 1, \cdots, r - 1.
\end{array}\right.
\end{equation*}
where $\lfloor x \rfloor$ is a floor function of $x$.
\label{cor:MOP-M}
\end{corollary}
\textbf{Proof.} The proof is almost identical to that of Lemma 
6 in \cite{WENO-MAIMi}.
$\hfill\square$ \\

\begin{corollary}
When $n_{\mathrm{cp}} = 1$, the $(2r-1)$th-order MOP-WENO-IM($k, A$) 
schemes can achieve the optimal order of accuracy if the mapping 
function $\big(g^{\mathrm{MOP-IM}}\big)_{s}(\omega)$ is applied to 
the original weights in the $(2r-1)$th-order WENO-JS schemes with 
requirement of $k \geq 2$ (except for the case $r = 2$).
\label{cor:MOP-IM}
\end{corollary}
\textbf{Proof.} The proof is almost identical to that of Theorem 2 
in \cite{WENO-IM}.
$\hfill\square$ \\

\begin{corollary}
The $(2r-1)$th-order MOP-WENO-PM$k$ schemes can achieve the optimal 
order of accuracy if the mapping function $\big(g^{\mathrm{MOP-PM}}
\big)_{s}(\omega)$ is applied to the original weights in the $(2r-1)$
th-order WENO-JS schemes with specific requirements for $k$ in 
following different cases: (I) Require $k \geq 1$ for 
$n_{\mathrm{cp}} = 0$; (II) Require $k \geq 1$ for $n_{\mathrm{cp}} =
1$; (III) Require $k \geq 3$ for $n_{\mathrm{cp}} = 2$.
\label{cor:MOP-PMk}
\end{corollary}
\textbf{Proof.} The proof is almost identical to that of Proposition 
1 in \cite{WENO-PM}.
$\hfill\square$ \\

\begin{corollary}
The $(2r - 1)$th-order MOP-WENO-RM($mn0$) schemes can recover the 
optimal order of accuracy if the mapping function 
$\big(g^{\mathrm{MOP-RM}}\big)_{s}(\omega)$ is applied to the 
original weights in the $(2r-1)$th-order WENO-JS schemes with 
requirement of $n \geq \frac{1 + n_{\mathrm{cp}}}{r - 1 - 
n_{\mathrm{cp}}}$ for different values of $n_{\mathrm{cp}}$ with 
$1 \leq n_{\mathrm{cp}} < r-1$.
\label{cor:MOP-RM260}
\end{corollary}
\textbf{Proof.} The proof is almost identical to that of Theorem 3 
in \cite{WENO-RM260}.
$\hfill\square$ \\

\begin{corollary}
Let $\lceil x \rceil$ be a ceiling function of $x$, for 
$n_{\mathrm{cp}} < r - 1$, the $(2r-1)$th-order MOP-WENO-MAIM$1$ 
schemes can achieve the optimal order of accuracy if the mapping 
function $(g^{\mathrm{MOP-MAIM}1})_{s}(\omega)$ is applied to the 
original weights in the $(2r-1)$th-order WENO-JS schemes with 
requirement of $k \geq k^{\mathrm{MAIM}}$, where
\begin{equation*}
k^{\mathrm{MAIM}} = \Bigg\lceil \dfrac{r}{r - 1 - n_{\mathrm{cp}}} - 
2 \Bigg\rceil + \dfrac{1 + \Big( -1 \Big)^{\Bigg\lceil \dfrac{r}{r - 
1 - n_{\mathrm{cp}}} -2 \Bigg\rceil}}{2}.
\end{equation*}
\label{cor:MOP-MAIM1}
\end{corollary}
\textbf{Proof.} The proof is almost identical to that of Theorem 2 
in \cite{WENO-MAIMi}.
$\hfill\square$ \\

\begin{corollary}
For $n_{\mathrm{cp}} < r - 1$, the $(2r-1)$th-order 
MOP-WENO-ACM schemes can achieve the optimal order of accuracy if 
the mapping function $(g^{\mathrm{MOP-ACM}})_{s}(\omega)$ is applied 
to the original weights in the $(2r-1)$th-order WENO-JS schemes.
\label{cor:MOP-ACM}
\end{corollary}
\textbf{Proof.} The proof is almost identical to that of Theorem 2 
in \cite{WENO-ACM}.
$\hfill\square$ \\

\begin{corollary}
When $\mathrm{CFS}_{s} \ll \widetilde{d}_{0}$, for $n_{\mathrm{cp}} <
r - 1$, the $(2r-1)$th-order MOP-WENO-ACM$k$ schemes can achieve the 
optimal order of accuracy if the mapping function 
$(g^{\mathrm{MOP-ACM}k})_{s}(\omega)$ is applied to the 
original weights in the $(2r-1)$th-order WENO-JS schemes.
\label{cor:MOP-ACMk}
\end{corollary}
\textbf{Proof.} The proof is almost identical to that of Theorem 2 
in \cite{MOP-WENO-ACMk}.
$\hfill\square$ \\


%% file: article_numerics.tex
\section{Numerical results}
\label{NumericalTests}
In this section, we compare the numerical performances of the 
MOP-WENO-X schemes with the corresponding existing mapped WENO-X 
schemes shown in Table \ref{table:GF_parameters}, as well as the 
classic WENO-JS scheme. As the performances of the WENO-ACM scheme 
and the MOP-WENO-ACM scheme are almost identical to those of the 
MIP-WENO-ACM$k$ scheme and the MOP-WENO-ACM$k$ scheme respectively, 
we do not present the solutions of the WENO-ACM scheme and the 
MOP-WENO-ACM scheme below for simplicity. Typical one-dimensional 
linear advection equation and two-dimensional Euler equations, with 
different initial conditions, are used to test the considered 
schemes. The presentation of these numerical tests in this section 
starts with the accuracy test of one-dimensional linear advection 
equation with four different initial conditions, followed by the 
long output time simulations of it with two different initial 
conditions including discontinuities, and finishes with 
two-dimensional simulations on the shock-vortex interaction and the 
2D Riemann problem. In all calculations below, $\epsilon$ is taken 
to be $10^{-40}$ for all schemes following the recommendations in 
\cite{WENO-M,WENO-IM}.

In the following numerical tests, the ODEs resulting from the 
semi-discretized PDEs are marched in time using the following 
explicit, third-order, strong stability preserving (SSP) Runge-Kutta 
method \cite{ENO-Shu1988,SSPRK1998,SSPRK2001}
\begin{equation*}
\begin{array}{l}
\begin{aligned}
&\overrightarrow{U}^{*} = \overrightarrow{U}^{n} + \Delta t 
\mathcal{L}(\overrightarrow{U}^{n}), \\
&\overrightarrow{U}^{**} = \dfrac{3}{4} \overrightarrow{U}^{n} + 
\dfrac{1}{4} \overrightarrow{U}^{*} + \dfrac{1}{4}\Delta t 
\mathcal{L}(\overrightarrow{U}^{*}), \\
&\overrightarrow{U}^{n + 1} = \dfrac{1}{3} \overrightarrow{U}^{n} + 
\dfrac{2}{3}\overrightarrow{U}^{**} + \dfrac{2}{3} \Delta t 
\mathcal{L}(\overrightarrow{U}^{**}),
\end{aligned}
\end{array}
\end{equation*}
where $\overrightarrow{U}^{*}$, $\overrightarrow{U}^{**}$ are the 
intermediate stages, $\overrightarrow{U}^{n}$ is the value of 
$\overrightarrow{U}$ at time level $t^{n} = n\Delta t$, and 
$\Delta t$ is the time step satisfying some proper CFL condition. 
The spatial operator $\mathcal{L}$ is defined as in 
Eq.\eqref{eq:discretizedFunction}, and the WENO reconstructions will 
be applied to obtain it.

\subsection{Accuracy test}
In this subsection, we solve the following one-dimensional linear 
advection equation
\begin{equation}
\dfrac{\partial u}{\partial t} + \dfrac{\partial u}{\partial x} = 0,
\quad -1 \leq x \leq 1,
\label{eq:LAE}
\end{equation}
with different initial conditions to test the accuracy of the 
considered WENO schemes. In all accuracy tests, the $L_{1}, L_{2}, 
L_{\infty}$ norms of the error are given as
\begin{equation*}
\displaystyle
\begin{array}{l}
L_{1} = h \cdot \displaystyle\sum\limits_{j} \big\lvert u_{j}^{
\mathrm{exact}} - (u_{h})_{j} \big\rvert, \\
L_{2} = \sqrt{h \cdot \displaystyle\sum\limits_{j} (u_{j}^{
\mathrm{exact}} - (u_{h})_{j})^{2}}, \\
L_{\infty} = \displaystyle\max_{j} \big\lvert u_{j}^{
\mathrm{exact}} - (u_{h})_{j} \big\rvert,
\end{array}
\end{equation*}
where $h = \Delta x$ is the uniform spatial step size, $(u_{h})_{j}$ 
is the numerical solution and $u_{j}^{\mathrm{exact}}$ is the exact 
solution.

\begin{example}
\rm{We calculate Eq.(\ref{eq:LAE}) with the periodic boundary 
condition using the following initial condition \cite{WENO-IM}} 
\label{LAE1}
\end{example} 
\begin{equation}
u(x, 0) = \sin ( \pi x ).
\label{eq:LAE:IC1}
\end{equation}
It is trivial to verify that although the initial condition in 
Eq.(\ref{eq:LAE:IC1}) has two first-order critical points, their 
first and third derivatives vanish simultaneously. The CFL number is 
taken to be $(\Delta x)^{2/3}$ to prevent the convergence rates of 
error from being influenced by time advancement. The calculation is 
run until a time of $t = 2.0$.

In Table \ref{table_LAE1}, we show the $L_{1}, L_{2}, L_{\infty}$ 
errors and corresponding convergence orders of various considered 
WENO schemes. The results of the three rows are $L_{1}$-, $L_{2}$- 
and $L_{\infty}$- norm errors and convergence orders (in brackets) 
in turn (similarly hereinafter). Unsurprisingly, the MOP-WENO-X 
schemes and the corresponding WENO-X schemes provide more accurate 
results than the WENO-JS scheme do in general. Naturally and as 
expected, all the considered schemes have gained the fifth-order 
convergence rate of accuracy. It can be found that the results of 
the MOP-WENO-X schemes are identical to those of the corresponding 
WENO-X schemes for all grid numbers except $N = 10$. As discussed in 
\cite{MOP-WENO-ACMk}, the cause of the accuracy loss for the 
computing cases of all MOP-WENO-X schemes with $N = 10$ is that the 
mapping functions of the MOP-WENO-X schemes have narrower optimal 
weight intervals (standing for the intervals about $\omega = d_{s}$ 
over which the mapping process attempts to use the corresponding 
optimal weights, see \cite{WENO-MAIMi,WENO-ACM}) than the 
corresponding WENO-X schemes.

\begin{table}[ht]
\begin{myFontSize}
\centering
\caption{Convergence properties of considered schemes solving $u_{t} 
+ u_{x} = 0$ with initial condition $u(x, 0) = \sin(\pi x)$.}
\label{table_LAE1}
\begin{tabular*}{\hsize}
{@{}@{\extracolsep{\fill}}lllllll@{}}
\hline
$N$       
& 10                  & 20                  & 40
& 80                  & 160                 & 320    \\
\hline
WENO-JS 
& 6.18628e-02(-)      & 2.96529e-03(4.3821) & 9.27609e-05(4.9985) 
& 2.89265e-06(5.0031) & 9.03392e-08(5.0009) & 2.82330e-09(4.9999)  \\
\space
& 4.72306e-02(-)      & 2.42673e-03(4.2826) & 7.64332e-05(4.9887) 
& 2.33581e-06(5.0322) & 7.19259e-08(5.0213) & 2.23105e-09(5.0107)  \\
\space
& 4.87580e-02(-)      & 2.57899e-03(4.2408) & 9.05453e-05(4.8320) 
& 2.90709e-06(4.9610) & 8.85753e-08(5.0365) & 2.72458e-09(5.0228)  \\
WENO-M 
& 2.01781e-02(-)      & 5.18291e-04(5.2829) & 1.59422e-05(5.0228) 
& 4.98914e-07(4.9979) & 1.56021e-08(4.9990) & 4.99356e-10(4.9977)  \\
\space
& 1.55809e-02(-)      & 4.06148e-04(5.2616) & 1.25236e-05(5.0193) 
& 3.91875e-07(4.9981) & 1.22541e-08(4.9991) & 3.83568e-10(4.9976)  \\
\space
& 1.47767e-02(-)      & 3.94913e-04(5.2256) & 1.24993e-05(4.9816) 
& 3.91808e-07(4.9956) & 1.22538e-08(4.9988) & 3.83541e-10(4.9977)  \\
MOP-WENO-M 
& 3.64427e-02(-)      & 5.18291e-04(6.1357) & 1.59422e-05(5.0228) 
& 4.98914e-07(4.9979) & 1.56021e-08(4.9990) & 4.88356e-10(4.9977)  \\
\space
& 2.95270e-02(-)      & 4.06148e-04(6.1839) & 1.25236e-05(5.0193) 
& 3.91875e-07(4.9981) & 1.22541e-08(4.9991) & 3.83568e-10(4.9976)  \\
\space
& 2.81876e-02(-)      & 3.94913e-04(6.1574) & 1.24993e-05(4.9816) 
& 3.91808e-07(4.9956) & 1.22538e-08(4.9988) & 3.83541e-10(4.9977)  \\
WENO-IM(2,0.1) 
& 1.58051e-02(-)      & 5.04401e-04(4.9697) & 1.59160e-05(4.9860) 
& 4.98863e-07(4.9957) & 1.56020e-08(4.9988) & 4.88355e-10(4.9977)  \\
\space
& 1.23553e-02(-)      & 3.96236e-04(4.9626) & 1.25033e-05(4.9860) 
& 3.91836e-07(4.9959) & 1.22541e-08(4.9989) & 3.83568e-10(4.9976)  \\
\space
& 1.19178e-02(-)      & 3.94458e-04(4.9171) & 1.24963e-05(4.9803) 
& 3.91797e-07(4.9953) & 1.22538e-08(4.9988) & 3.83547e-10(4.9977)  \\
MOP-WENO-IM(2,0.1) 
& 3.35513e-02(-)      & 5.04401e-04(6.0557) & 1.59160e-05(4.9860) 
& 4.98863e-07(4.9957) & 1.56020e-08(4.9988) & 4.88355e-10(4.9977)  \\
\space
& 2.75968e-02(-)      & 3.96236e-04(6.1220) & 1.25033e-05(4.9860) 
& 3.91836e-07(4.9959) & 1.22541e-08(4.9989) & 3.83568e-10(4.9976)  \\
\space
& 2.71898e-02(-)      & 3.94458e-04(6.1071) & 1.24963e-05(4.9803) 
& 3.91797e-07(4.9953) & 1.22538e-08(4.9988) & 3.83547e-10(4.9977)  \\
WENO-PM6 
& 1.74869e-02(-)      & 5.02923e-04(5.1198) & 1.59130e-05(4.9821) 
& 4.98858e-07(4.9954) & 1.56020e-08(4.9988) & 4.88355e-10(4.9977)  \\
\space
& 1.35606e-02(-)      & 3.95215e-04(5.1006) & 1.25010e-05(4.9825) 
& 3.91831e-07(4.9957) & 1.22541e-08(4.9989) & 3.83568e-10(4.9976)  \\
\space
& 1.27577e-02(-)      & 3.94515e-04(5.0151) & 1.24960e-05(4.9805) 
& 3.91795e-07(4.9952) & 1.22538e-08(4.9988) & 3.83543e-10(4.9977)  \\
MOP-WENO-PM6 
& 3.54584e-02(-)      & 5.02923e-04(6.1396) & 1.59130e-05(4.9821) 
& 4.98858e-07(4.9954) & 1.56020e-08(4.9988) & 4.88355e-10(4.9977)  \\
\space
& 2.88246e-02(-)      & 3.95215e-04(6.1885) & 1.25010e-05(4.9825) 
& 3.91831e-07(4.9957) & 1.22541e-08(4.9989) & 3.83568e-10(4.9976)  \\
\space
& 2.76902e-02(-)      & 3.94515e-04(6.1332) & 1.24960e-05(4.9805) 
& 3.91795e-07(4.9952) & 1.22538e-08(4.9988) & 3.83543e-10(4.9977)  \\
WENO-PPM5 
& 1.73978e-02(-)      & 5.03464e-04(5.1109) & 1.59131e-05(4.9836) 
& 4.98858e-07(4.9954) & 1.56020e-08(4.9988) & 4.88356e-10(4.9977)  \\
\space
& 1.34998e-02(-)      & 3.95644e-04(5.0926) & 1.25011e-05(4.9841) 
& 3.91831e-07(4.9957) & 1.22541e-08(4.9989) & 3.83568e-10(4.9976)  \\
\space
& 1.27018e-02(-)      & 3.94865e-04(5.0075) & 1.24961e-05(4.9818) 
& 3.91795e-07(4.9952) & 1.22538e-08(4.9988) & 3.83528e-10(4.9978)  \\
MOP-WENO-PPM5 
& 3.49872e-02(-)      & 5.03464e-04(6.1188) & 1.59131e-05(4.9836) 
& 4.98858e-07(4.9954) & 1.56020e-08(4.9988) & 4.88356e-10(4.9977)  \\
\space
& 2.85173e-02(-)      & 3.95644e-04(6.1715) & 1.25011e-05(4.9841) 
& 3.91831e-07(4.9957) & 1.22541e-08(4.9989) & 3.83568e-10(4.9976)  \\
\space
& 2.75955e-02(-)      & 3.94865e-04(6.1269) & 1.24961e-05(4.9818) 
& 3.91795e-07(4.9952) & 1.22538e-08(4.9988) & 3.83528e-10(4.9978)  \\
WENO-RM(260) 
& 1.52661e-02(-)      & 5.02845e-04(4.9241) & 1.59130e-05(4.9818) 
& 4.98858e-07(4.9954) & 1.56020e-08(4.9988) & 4.88355e-10(4.9977)  \\
\space
& 1.19792e-02(-)      & 3.95138e-04(4.9220) & 1.25010e-05(4.9822) 
& 3.91831e-07(4.9957) & 1.22541e-08(4.9989) & 3.83568e-10(4.9976)  \\
\space
& 1.17698e-02(-)      & 3.94406e-04(4.8993) & 1.24960e-05(4.9801) 
& 3.91795e-07(4.9952) & 1.22538e-08(4.9988) & 3.83543e-10(4.9977)  \\
MOP-WENO-RM(260) 
& 3.29243e-02(-)      & 5.02845e-04(6.0329) & 1.59130e-05(4.9818) 
& 4.98858e-07(4.9954) & 1.56020e-08(4.9988) & 4.88355e-10(4.9977)  \\
\space
& 2.73131e-02(-)      & 3.95138e-04(6.1111) & 1.25010e-05(4.9822) 
& 3.91831e-07(4.9957) & 1.22541e-08(4.9989) & 3.83568e-10(4.9976)  \\
\space
& 2.73015e-02(-)      & 3.94406e-04(6.1132) & 1.24960e-05(4.9801) 
& 3.91795e-07(4.9952) & 1.22538e-08(4.9988) & 3.83543e-10(4.9977)  \\
WENO-MAIM1 
& 6.13264e-02(-)      & 5.08205e-04(6.9150) & 1.59130e-05(4.9971) 
& 4.98858e-07(4.9954) & 1.56020e-08(4.9988) & 4.88355e-10(4.9977)  \\
\space
& 4.81375e-02(-)      & 4.26155e-04(6.8196) & 1.25010e-05(5.0913) 
& 3.91831e-07(4.9957) & 1.22541e-08(4.9989) & 3.83568e-10(4.9976)  \\
\space
& 4.86913e-02(-)      & 5.03701e-04(6.5950) & 1.24960e-05(5.3330) 
& 3.91795e-07(4.9952) & 1.22538e-08(4.9988) & 3.83543e-10(4.9977)  \\
MOP-WENO-MAIM1 
& 6.63923e-02(-)      & 5.08205e-04(7.0295) & 1.59130e-05(4.9971) 
& 4.98858e-07(4.9954) & 1.56020e-08(4.9988) & 4.88355e-10(4.9977)  \\
\space
& 5.17462e-02(-)      & 4.26155e-04(6.9239) & 1.25010e-05(5.0913) 
& 3.91831e-07(4.9957) & 1.22541e-08(4.9989) & 3.83568e-10(4.9976)  \\
\space
& 5.19799e-02(-)      & 5.03701e-04(6.6892) & 1.24960e-05(5.3330) 
& 3.91795e-07(4.9952) & 1.22538e-08(4.9988) & 3.83543e-10(4.9977)  \\
MIP-WENO-ACM$k$ 
& 1.52184e-02(-)      & 5.02844e-04(4.9196) & 1.59130e-05(4.9818) 
& 4.98858e-07(4.9954) & 1.56020e-08(4.9988) & 4.88355e-10(4.9977)  \\
\space
& 1.19442e-02(-)      & 3.95138e-04(4.9178) & 1.25010e-05(4.9822) 
& 3.91831e-07(4.9957) & 1.22541e-08(4.9989) & 3.83568e-10(4.9976)  \\
\space
& 1.17569e-02(-)      & 3.94406e-04(4.8977) & 1.24960e-05(4.9801) 
& 3.91795e-07(4.9952) & 1.22538e-08(4.9988) & 3.83543e-10(4.9977)  \\
MOP-WENO-ACM$k$ 
& 3.29609e-02(-)      & 5.02844e-04(6.0345) & 1.59130e-05(4.9818) 
& 4.98858e-07(4.9954) & 1.56020e-08(4.9988) & 4.88355e-10(4.9977)  \\
\space
& 2.72363e-02(-)      & 3.95138e-04(6.1070) & 1.25010e-05(4.9822) 
& 3.91831e-07(4.9957) & 1.22541e-08(4.9989) & 3.83568e-10(4.9976)  \\
\space
& 2.70295e-02(-)      & 3.94406e-04(6.0987) & 1.24960e-05(4.9801) 
& 3.91795e-07(4.9952) & 1.22538e-08(4.9988) & 3.83543e-10(4.9977)  \\
\hline
\end{tabular*}
\end{myFontSize}
\end{table}

\begin{example}
\rm{We calculate Eq.(\ref{eq:LAE}) with the periodic boundary 
condition using the following initial condition \cite{WENO-M}}
\label{LAE2}
\end{example}
\begin{equation}
u(x, 0) = \sin \bigg( \pi x - \dfrac{\sin(\pi x)}{\pi} \bigg), 
\label{eq:LAE:IC2}
\end{equation}
This particular initial condition has two first-order critical 
points, which both have a non-vanishing third derivative. Again, the 
CFL number is set to be $(\Delta x)^{2/3}$ and the calculation is 
run until a time of $t = 2.0$.

Table \ref{table_LAE2} compares the $L_{1}, L_{2}, L_{\infty}$ 
errors and corresponding convergence orders obtained from the 
considered schemes. It is evident that the WENO-X schemes and the 
corresponding MOP-WENO-X schemes can achieve the optimal convergence 
orders. Unsurprisingly, the WENO-JS scheme gives less accurate 
results than the other schemes, and its $L_{\infty}$ convergence 
order decreases by almost 2 orders leading to the noticeable drops 
of the $L_{1}$ and $L_{2}$ convergence orders. It is noteworthy that 
when the grid number is too small, like $N\leq 40$, in terms of 
accuracy, the MOP-WENO-X schemes provide less accurate results than 
those of the corresponding WENO-X schemes. As mentioned in Example 
\ref{LAE1}, the cause of this kind of accuracy loss is that the 
mapping functions of the MOP-WENO-X schemes have narrower optimal 
weight intervals than the corresponding WENO-X schemes, and this 
issue can surely be addressed by increasing the grid number. 
Therefore, as expected, the MOP-WENO-X schemes show equally accurate 
numerical solutions like those of the corresponding WENO-X schemes 
when the grid number $N \geq 80$.

\begin{table}[ht]
\begin{myFontSize}
\centering
\caption{Convergence properties of considered schemes solving $u_{t} 
+ u_{x} = 0$ with initial condition $u(x, 0) = \sin(\pi x - \sin(\pi 
x)/\pi )$.}
\label{table_LAE2}
\begin{tabular*}{\hsize}
{@{}@{\extracolsep{\fill}}lllllll@{}}
\hline
$N$       
& 10                  & 20                  & 40
& 80                  & 160                 & 320    \\
\hline
WENO-JS 
& 1.24488e-01(-)      & 1.01260e-02(3.6199) & 7.22169e-04(3.8096) 
& 3.42286e-05(4.3991) & 1.58510e-06(4.4326) & 7.95517e-08(4.3165)  \\
\space
& 1.09463e-01(-)      & 8.72198e-03(3.6496) & 6.76133e-04(3.6893) 
& 3.63761e-05(4.2162) & 2.29598e-06(3.9858) & 1.68304e-07(3.7700)  \\
\space
& 1.24471e-01(-)      & 1.43499e-02(3.1167) & 1.09663e-03(3.7099) 
& 9.02485e-05(3.6030) & 8.24022e-06(3.4531) & 8.31702e-07(3.3085)  \\
WENO-M 
& 7.53259e-02(-)      & 3.70838e-03(4.3443) & 1.45082e-04(4.6758) 
& 4.80253e-06(4.9169) & 1.52120e-07(4.9805) & 4.77083e-09(4.9948)  \\
\space
& 6.39017e-02(-)      & 3.36224e-03(4.2484) & 1.39007e-04(4.5962) 
& 4.52646e-06(4.9406) & 1.42463e-07(4.9897) & 4.45822e-09(4.9980)  \\
\space
& 7.49250e-02(-)      & 5.43666e-03(3.7847) & 2.18799e-04(4.6350) 
& 6.81451e-06(5.0049) & 2.14545e-07(4.9893) & 6.71080e-09(4.9987)  \\
MOP-WENO-M 
& 9.41832e-02(-)      & 6.59540e-03(3.8359) & 2.60456e-04(4.6623) 
& 4.80253e-06(5.7611) & 1.52120e-07(4.9805) & 4.77083e-09(4.9948)  \\
\space
& 8.03446e-02(-)      & 6.37937e-03(3.6547) & 2.50868e-04(4.6684) 
& 4.52646e-06(5.7924) & 1.42463e-07(4.9897) & 4.45822e-09(4.9980)  \\
\space
& 9.78919e-02(-)      & 8.97094e-03(3.4479) & 4.10480e-04(4.4499) 
& 6.81451e-06(5.9126) & 2.14545e-07(4.9893) & 6.71080e-09(4.9987)  \\
WENO-IM(2,0.1) 
& 8.38131e-02(-)      & 4.30725e-03(4.2823) & 1.51327e-04(4.8310) 
& 4.85592e-06(4.9618) & 1.52659e-07(4.9914) & 4.77654e-09(4.9982)  \\
\space
& 6.71285e-02(-)      & 3.93700e-03(4.0918) & 1.41737e-04(4.7958) 
& 4.53602e-06(4.9656) & 1.42479e-07(4.9926) & 4.45805e-09(4.9982)  \\
\space
& 7.62798e-02(-)      & 5.84039e-03(3.7072) & 2.10531e-04(4.7940) 
& 6.82606e-06(4.9468) & 2.14534e-07(4.9918) & 6.71079e-09(4.9986)  \\
MOP-WENO-IM(2,0.1) 
& 8.49795e-02(-)      & 7.01287e-03(3.5990) & 2.59767e-04(4.7547) 
& 4.85592e-06(5.7413) & 1.52659e-07(4.9914) & 4.77654e-09(4.9982)  \\
\space
& 7.29388e-02(-)      & 6.80019e-03(3.4230) & 2.51121e-04(4.7591) 
& 4.53602e-06(5.7908) & 1.42479e-07(4.9926) & 4.45805e-09(4.9982)  \\
\space
& 9.47429e-02(-)      & 9.96943e-03(3.2484) & 4.01785e-04(4.6330) 
& 6.82606e-06(5.8792) & 2.14534e-07(4.9918) & 6.71079e-09(4.9986)  \\
WENO-PM6 
& 9.51313e-02(-)      & 4.82173e-03(4.3023) & 1.55428e-04(4.9552) 
& 4.87327e-06(4.9952) & 1.52750e-07(4.9956) & 4.77729e-09(4.9988)  \\
\space
& 7.83600e-02(-)      & 4.29510e-03(4.1894) & 1.43841e-04(4.9001) 
& 4.54036e-06(4.9855) & 1.42488e-07(4.9939) & 4.45807e-09(4.9983)  \\
\space
& 9.32356e-02(-)      & 5.91037e-03(3.9796) & 2.09540e-04(4.8180) 
& 6.83270e-06(4.9386) & 2.14532e-07(4.9932) & 6.71079e-09(4.9986)  \\
MOP-WENO-PM6 
& 1.00298e-01(-)      & 5.84504e-03(4.1009) & 2.51725e-04(4.5373) 
& 4.87327e-06(5.6908) & 1.52750e-07(4.9956) & 4.77729e-09(4.9988)  \\
\space
& 8.49034e-02(-)      & 5.80703e-03(3.8699) & 2.40678e-04(4.5926) 
& 4.54036e-06(5.7282) & 1.42488e-07(4.9939) & 4.45807e-09(4.9983)  \\
\space
& 9.88357e-02(-)      & 9.01779e-03(3.4542) & 3.66822e-04(4.6196) 
& 6.83270e-06(5.7465) & 2.14532e-07(4.9932) & 6.71079e-09(4.9986)  \\
WENO-PPM5 
& 9.22982e-02(-)      & 4.68376e-03(4.3006) & 1.55745e-04(4.9104) 
& 4.88795e-06(4.9938) & 1.52852e-07(4.9990) & 4.77759e-09(4.9997)  \\
\space
& 7.46925e-02(-)      & 4.18882e-03(4.1563) & 1.44018e-04(4.8622) 
& 4.54528e-06(4.9857) & 1.42506e-07(4.9953) & 4.45812e-09(4.9984)  \\
\space
& 8.46229e-02(-)      & 5.92748e-03(3.8356) & 2.09420e-04(4.8229) 
& 6.83617e-06(4.9371) & 2.14527e-07(4.9940) & 6.71080e-09(4.9985)  \\
MOP-WENO-PPM5 
& 9.50369e-02(-)      & 6.27179e-03(3.9215) & 2.52600e-04(4.6340) 
& 4.88795e-06(5.6915) & 1.52852e-07(4.9990) & 4.77759e-09(4.9997)  \\
\space
& 8.08190e-02(-)      & 6.11267e-03(3.7248) & 2.41656e-04(4.6608) 
& 4.54528e-06(5.7324) & 1.42506e-07(4.9953) & 4.45812e-09(4.9984)  \\
\space
& 9.65522e-02(-)      & 8.98120e-03(3.4263) & 3.69338e-04(4.6039) 
& 6.83617e-06(5.7556) & 2.14527e-07(4.9940) & 6.71080e-09(4.9985)  \\
WENO-RM(260) 
& 8.24328e-02(-)      & 4.37642e-03(4.2354) & 1.52200e-04(4.8457) 
& 4.86434e-06(4.9676) & 1.52735e-07(4.9931) & 4.77728e-09(4.9987)  \\
\space
& 6.64590e-02(-)      & 4.00547e-03(4.0524) & 1.42162e-04(4.8164) 
& 4.53769e-06(4.9694) & 1.42486e-07(4.9931) & 4.45807e-09(4.9983)  \\
\space
& 7.64206e-02(-)      & 5.88375e-03(3.6992) & 2.09889e-04(4.8090) 
& 6.83016e-06(4.9416) & 2.14533e-07(4.9926) & 6.71079e-09(4.9986)  \\
MOP-WENO-RM(260) 
& 8.96509e-02(-)      & 6.87612e-03(3.7047) & 2.59418e-04(4.7282) 
& 4.86434e-06(5.7369) & 1.52735e-07(4.9931) & 4.77728e-09(4.9987)  \\
\space
& 7.51169e-02(-)      & 6.65488e-03(3.4967) & 2.51194e-04(4.7275) 
& 4.53769e-06(5.7907) & 1.42486e-07(4.9931) & 4.45807e-09(4.9983)  \\
\space
& 9.20962e-02(-)      & 9.75043e-03(3.2396) & 4.03065e-04(4.5964) 
& 6.83016e-06(5.8829) & 2.14533e-07(4.9926) & 6.71079e-09(4.9986)  \\
WENO-MAIM1 
& 1.24659e-01(-)      & 8.07923e-03(3.9476) & 3.32483e-04(4.6029) 
& 1.01162e-05(5.0385) & 1.52910e-07(6.0478) & 4.77728e-09(5.0003)  \\
\space
& 1.14152e-01(-)      & 7.08117e-03(4.0108) & 3.36264e-04(4.3963) 
& 1.49724e-05(4.4892) & 1.42515e-07(6.7150) & 4.45807e-09(4.9986)  \\
\space
& 1.40438e-01(-)      & 1.03772e-02(3.7584) & 6.62891e-04(3.9685) 
& 4.48554e-05(3.8854) & 2.14522e-07(7.7080) & 6.71079e-09(4.9985)  \\
MOP-WENO-MAIM1 
& 1.27999e-01(-)      & 7.62753e-03(4.0688) & 3.37132e-04(4.4998) 
& 1.01162e-05(5.0586) & 1.52910e-07(6.0478) & 4.77728e-09(5.0003)  \\
\space
& 1.12692e-01(-)      & 6.93240e-03(4.0229) & 3.36497e-04(4.3647) 
& 1.49724e-05(4.4902) & 1.42515e-07(6.7150) & 4.45807e-09(4.9986)  \\
\space
& 1.31113e-01(-)      & 1.27480e-02(3.3625) & 6.40953e-04(4.3139) 
& 4.48554e-05(3.8369) & 2.14522e-07(7.7080) & 6.71079e-09(4.9985)  \\
MIP-WENO-ACM$k$ 
& 8.75629e-02(-)      & 4.39527e-03(4.3163) & 1.52219e-04(4.8517) 
& 4.86436e-06(4.9678) & 1.52735e-07(4.9931) & 4.77728e-09(4.9987)  \\
\space
& 6.98131e-02(-)      & 4.02909e-03(4.1150) & 1.42172e-04(4.8247) 
& 4.53770e-06(4.9695) & 1.42486e-07(4.9931) & 4.45807e-09(4.9983)  \\
\space
& 7.91292e-02(-)      & 5.89045e-03(3.7478) & 2.09893e-04(4.8107) 
& 6.83017e-06(4.9416) & 2.14533e-07(4.9926) & 6.71079e-09(4.9986)  \\
MOP-WENO-ACM$k$ 
& 9.08634e-02(-)      & 7.09246e-03(3.6793) & 2.59429e-04(4.7729) 
& 4.86436e-06(5.7369) & 1.52735e-07(4.9931) & 4.77728e-09(4.9987)  \\
\space
& 7.58160e-02(-)      & 6.88532e-03(3.4609) & 2.51208e-04(4.7766) 
& 4.53770e-06(5.7908) & 1.42486e-07(4.9931) & 4.45807e-09(4.9983)  \\
\space
& 9.29135e-02(-)      & 1.01479e-02(3.1947) & 4.03069e-04(4.6540) 
& 6.83017e-06(5.8830) & 2.14533e-07(4.9926) & 6.71079e-09(4.9986)  \\
\hline
\end{tabular*}
\end{myFontSize}
\end{table}

\begin{example}
\rm{We calculate Eq.(\ref{eq:LAE}) using the following initial 
condition \cite{WENO-PM}} 
\label{LAE4}
\end{example} 
\begin{equation}
u(x, 0) = \sin^{9} ( \pi x ), 
\label{eq:LAE:IC4}
\end{equation}
with the periodic boundary condition. It is trivial to verify that 
this initial condition has high-order critical points. We also set 
the CFL number to be $(\Delta x)^{2/3}$.  

Table \ref{table_LAE4} shows the $L_{1}, L_{2}, L_{\infty}$ errors 
of the considered WENO schemes at several output times with a 
uniform mesh size of $\Delta x = 1/200$. In order to measure the 
dissipations of the schemes, as the authors did in 
\cite{MOP-WENO-ACMk}, we give the corresponding increased errors (in 
percentage) compared to the MIP-WENO-ACM$k$ scheme which gives 
solutions with very low dissipations. Taking the $L_{1}$-norm error 
as an example, let $L_{1}^{*}(t)$ and $L_{1}^{\mathrm{Y}}(t)$ the 
$L_{1}$-norm errors of the MIP-WENO-ACM$k$ scheme and the compared 
scheme, the increased errors at output time $t$ is computed by 
$\frac{L_{1}^{\mathrm{Y}}(t)-L_{1}^{*}(t)}{L_{1}^{*}(t)}\times100\%$
. From Table \ref{table_LAE4}, we can observe that: (1) the 
WENO-JS scheme shows the largest increased errors for no matter 
short or long output times; (2) at short output times, like 
$t \leq 100$, the solutions computed by the WENO-M scheme are 
closer to those of the MIP-WENO-ACM$k$ scheme, leading to smaller 
increased errors, than the corresponding MOP-WENO-M scheme; (3) 
however, when the output time is larger, like $t \geq 200$, the 
solutions computed by the MOP-WENO-M scheme, whose increased errors 
do not get larger but evidently decrease, are closer to those of 
the MIP-WENO-ACM$k$ scheme than the corresponding WENO-M scheme 
whose errors increase dramatically leading to significantly larger 
increased errors; (4) although the errors of the MOP-WENO-X schemes 
except the MOP-WENO-M scheme are not as small as those of the 
corresponding WENO-X schemes, these errors can maintain a 
considerable level leading to acceptable increased errors that are 
much lower than those of the WENO-JS and WENO-M schemes.

Actually, as mentioned in Example \ref{LAE1} and Example \ref{LAE2}, 
the cause of the slight accuracy loss discussed above is that the 
mapping function of the MOP-WENO-X scheme has narrower optimal 
weight intervals than the corresponding WENO-X schemes, and one can 
easily overcome this drawback by increasing the grid number. To 
demonstrate this, we calculate this problem using the same schemes 
at the same output times with a larger grid number of $N = 800$. The 
results are shown in Table \ref{table_LAE4:N800}, and we can see 
that: (1) the errors of the MOP-WENO-X schemes get closer to those 
of the MIP-WENO-ACM$k$ scheme when the grid number increases from 
$N = 200$ to $N = 800$, resulting in the significant decrease of 
the increased errors, and in different words, the errors of the 
MOP-WENO-X schemes and the MIP-WENO-ACM$k$ scheme are so close that 
one can ignore their differences; (2) although the errors of the 
WENO-JS and WENO-M schemes get smaller when the grid number 
increases from $N=200$ to $N=800$, their increased errors become 
very very large; (3) naturally, the increased errors of the 
MOP-WENO-X schemes are extremely smaller than those of the WENO-JS 
and WENO-M schemes. Actually, it is an important advantage of the 
MOP-WENO-X schemes that can maintain comparably high resolution for 
long output times. In the next subsection we will make a further 
discussion for this.

\begin{figure}[ht]
\centering
  \includegraphics[height=0.289\textwidth]
  {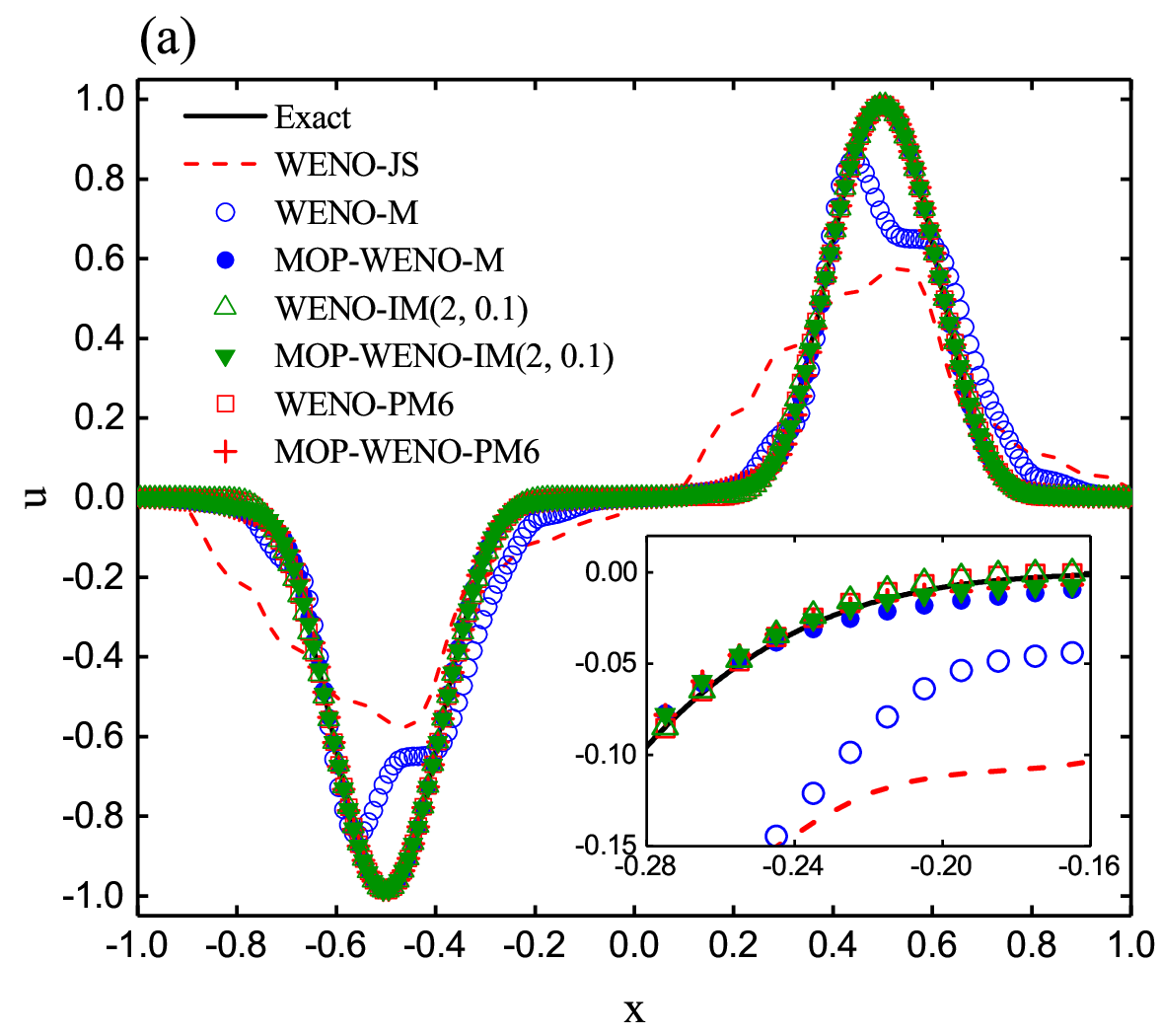}
  \includegraphics[height=0.289\textwidth]
  {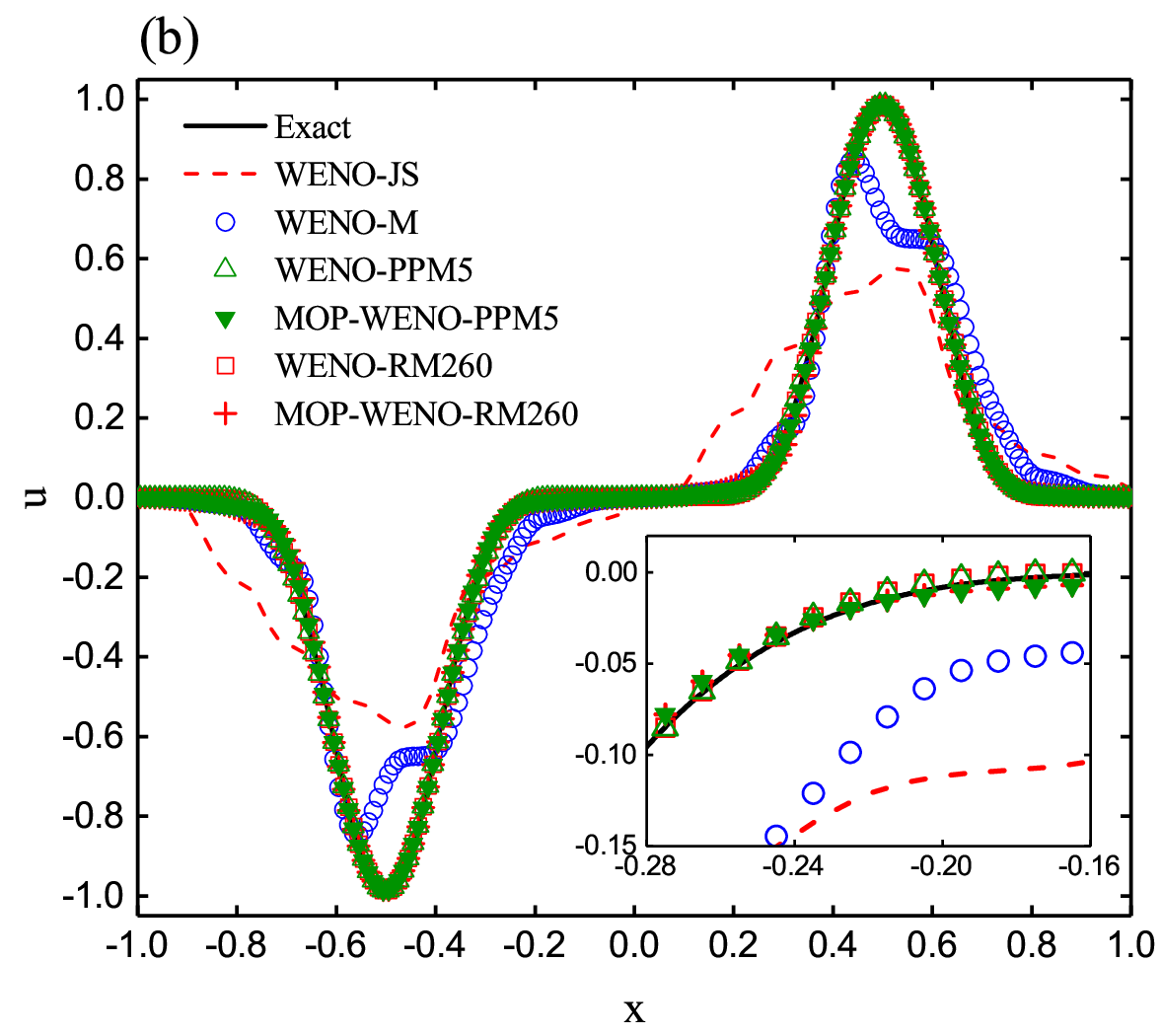}
  \includegraphics[height=0.289\textwidth]
  {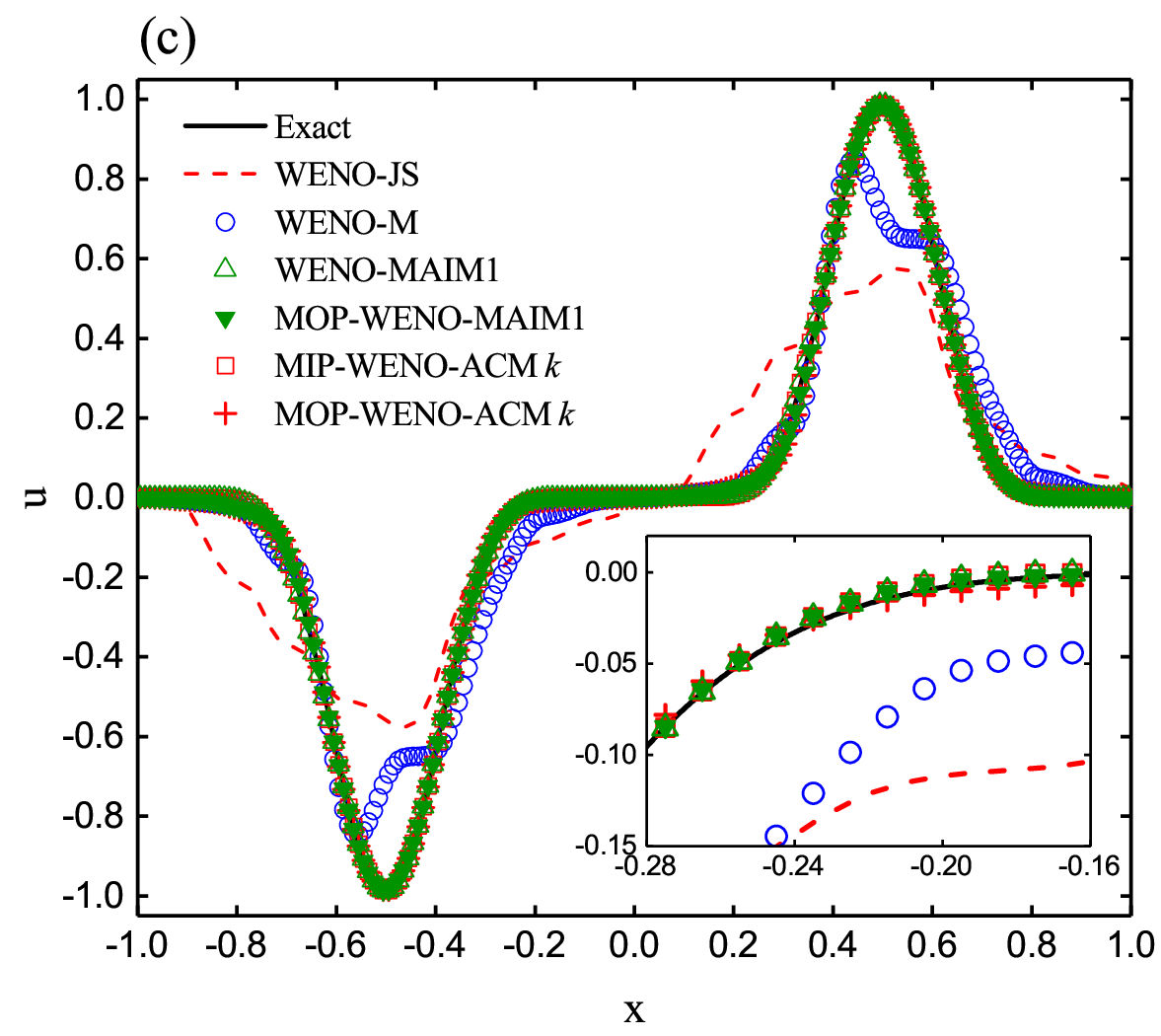}  
  \caption{Performance of various WENO schemes for Example \ref{LAE4}
  at output time $t = 1000$ with a uniform mesh size of $\Delta x= 
  1/200$.}
\label{fig:ex:LAE4:N200}
\end{figure}

\begin{figure}[ht]
\centering
  \includegraphics[height=0.289\textwidth]
  {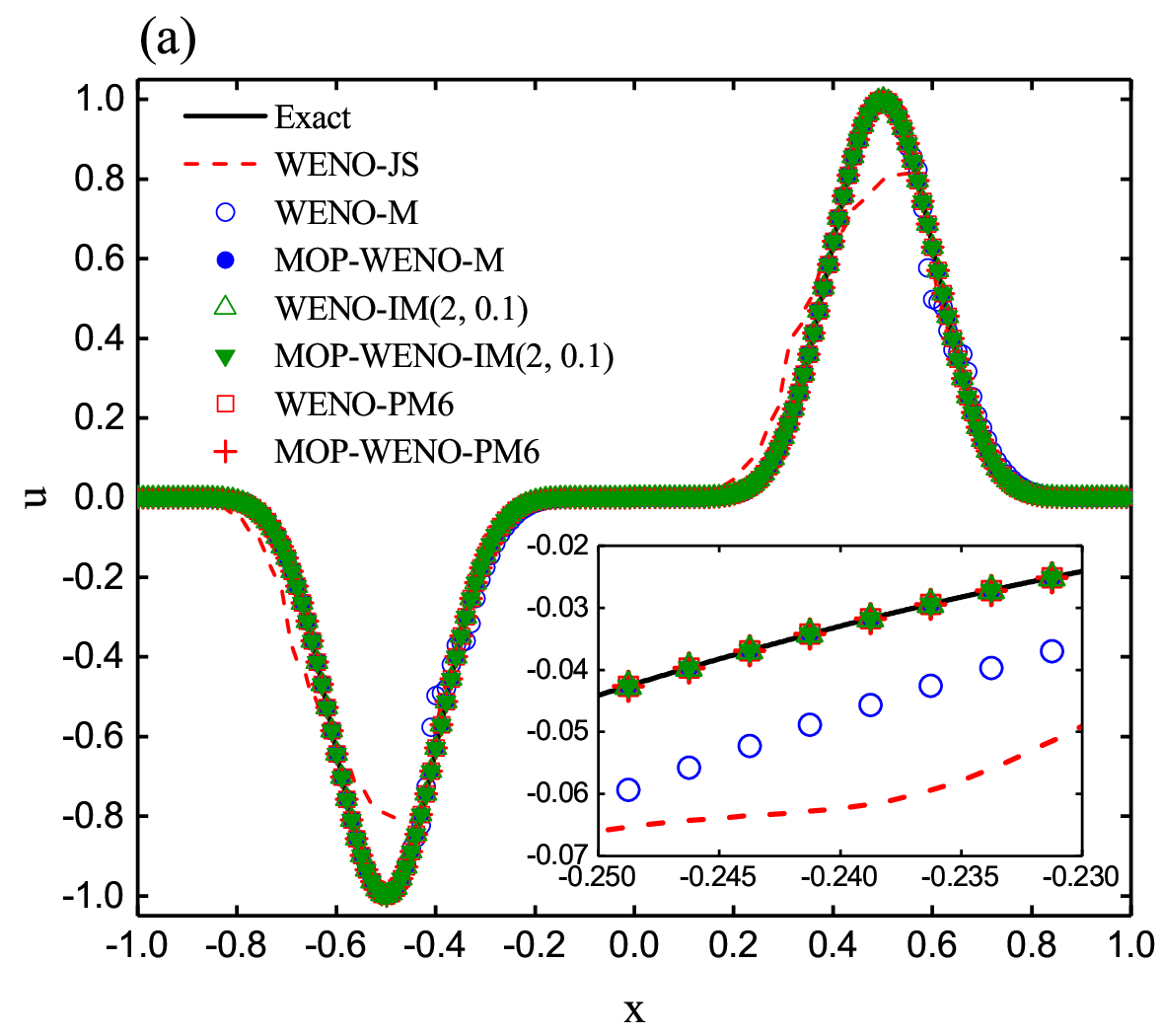}
  \includegraphics[height=0.289\textwidth]
  {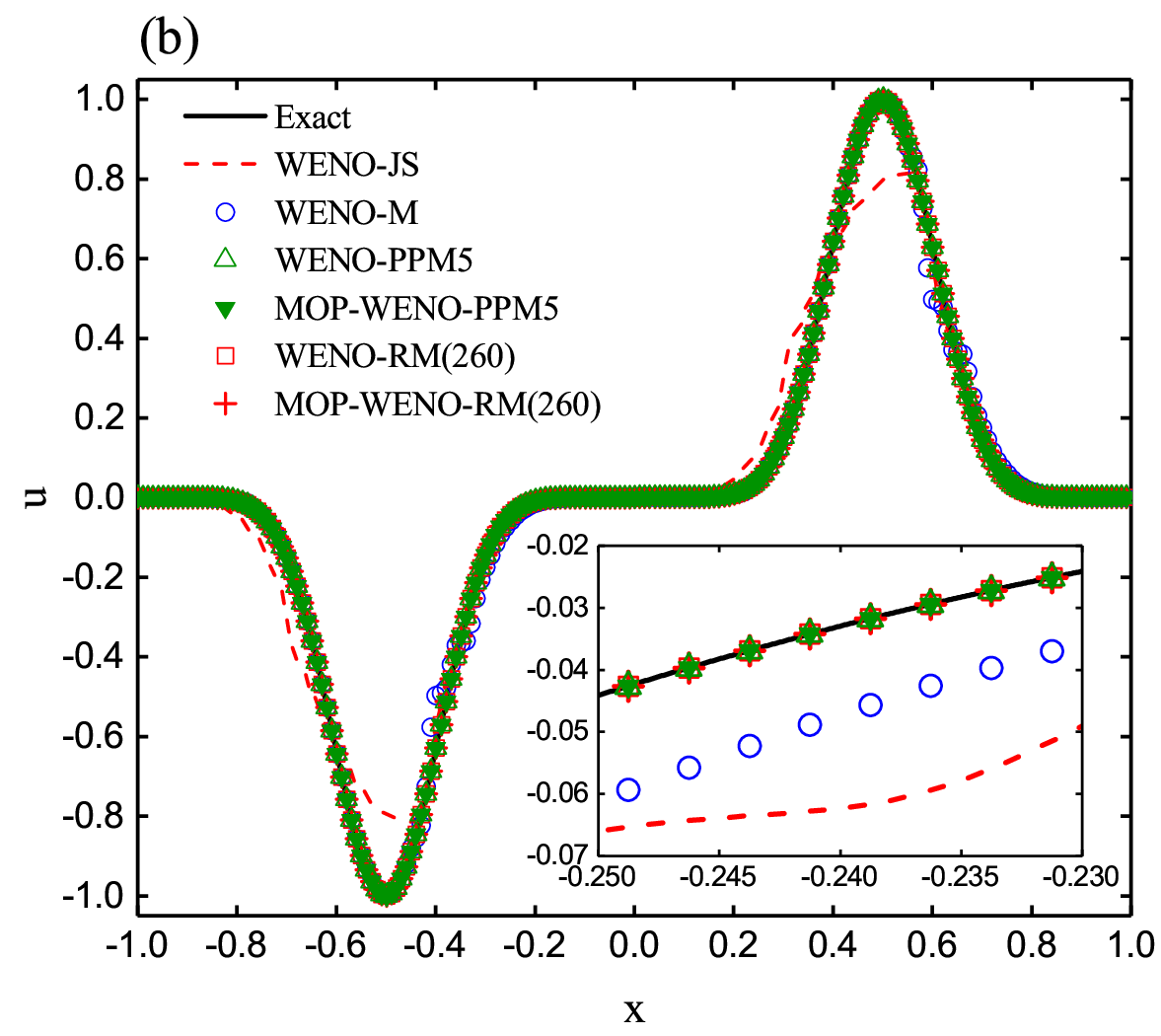}
  \includegraphics[height=0.289\textwidth]
  {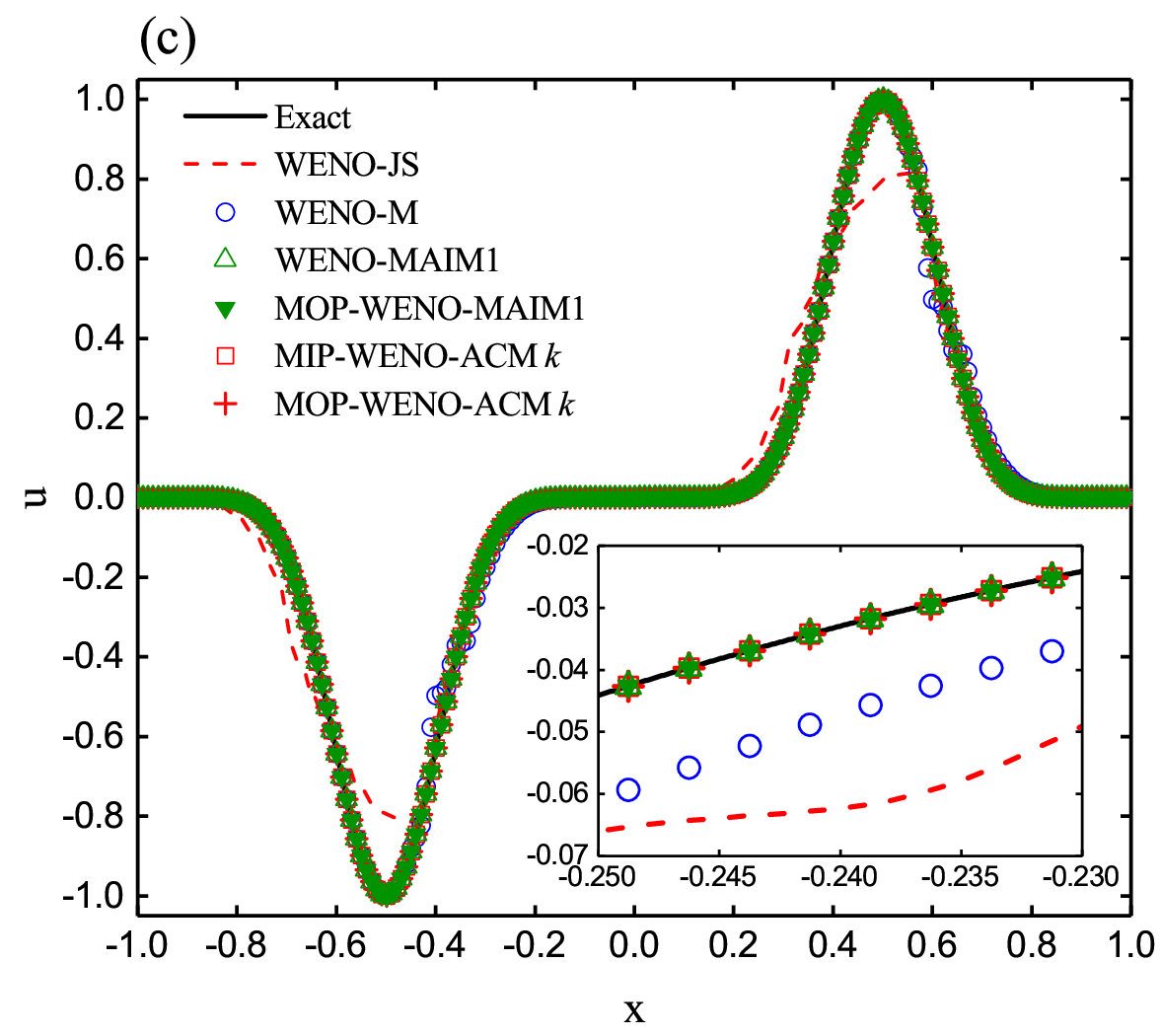}  
  \caption{Performance of various WENO schemes for Example \ref{LAE4}
  at output time $t = 1000$ with a uniform mesh size of $\Delta x= 
  1/800$.}
\label{fig:ex:LAE4:N800}
\end{figure}

In Fig. \ref{fig:ex:LAE4:N200} and Fig. \ref{fig:ex:LAE4:N800}, we 
plot the solutions computed by various schemes at output time 
$t=1000$ with the grid number of $N = 200$ and $N = 800$, 
respectively. For $N=200$, Fig. \ref{fig:ex:LAE4:N200} shows that: 
(1) the MOP-WENO-M scheme provides the result with far higher 
resolution than the corresponding WENO-M scheme which gives result 
with slightly better resolution than the lowest one computed by the 
WENO-JS scheme; (2) the result of the MOP-WENO-MAIM1 scheme is very 
close to that of its corresponding WENO-MAIM1 scheme; (3) the 
results of the other MOP-WENO-X schemes show far better 
resolutions than the WENO-M and WENO-JS schemes, although they give 
results with very slightly lower resolutions than their 
corresponding WENO-X schemes because of the narrower optimal 
weight intervals. Actually, we can amend this minor issue by using a 
larger grid number. Consequently, for $N = 800$, it can be seen from 
Fig. \ref{fig:ex:LAE4:N800} that: (1) all the MOP-WENO-X schemes 
produce results very close to those of their corresponding mapped 
WENO-X schemes with extremely high resolutions except the case of X 
= M; (2) the MOP-WENO-M scheme also produces result with very high 
resolution, while the resolutions of the results from the WENO-M and 
WENO-JS schemes are much lower.

\begin{table}[ht]
\begin{myFontSize}
\centering
\caption{Performance of various considered schemes solving $u_{t} + 
u_{x} = 0$ with $u(x, 0) = \sin^{9} (\pi x), \Delta x= 1/200$.}
\label{table_LAE4}
\begin{tabular*}{\hsize}
{@{}@{\extracolsep{\fill}}llllll@{}}
\hline
Scheme  & $t=10$ & $t=100$ & $t=200$ & $t=500$ & $t=1000$ \\
\hline
WENO-JS         & 3.86931e-04(359.06\%)  
                & 5.42288e-03(548.87\%)
                & 2.35657e-02\textbf{(1323.42\%)}  
                & 1.55650e-01\textbf{(3832.05\%)} 
                & 2.91359e-01\textbf{(3920.28\%)} \\
{}              & 3.52611e-04(330.48\%)  
                & 5.17716e-03(539.41\%)
                & 2.68753e-02\textbf{(1580.45\%)}  
                & 1.46859e-01\textbf{(3716.48\%)} 
                & 2.66692e-01\textbf{(3595.71\%)} \\
{}              & 5.36940e-04(288.51\%)  
                & 1.20056e-02(780.15\%)
                & 6.47820e-02\textbf{(2308.66\%)}  
                & 2.57663e-01\textbf{(3891.29\%)} 
                & 4.44664e-01\textbf{(3556.96\%)} \\
WENO-M          & 8.90890e-05(5.70\%)  
                & 1.29154e-03(54.54\%)
                & 5.74021e-03\textbf{(246.72\%)}
                & 4.89290e-02\textbf{(1136.05\%)}
                & 1.34933e-01\textbf{(1761.86\%)} \\
{}              & 8.32089e-05(1.58\%)  
                & 1.28740e-03(59.00\%)
                & 7.66721e-03\textbf{(379.41\%)}  
                & 6.23842e-02\textbf{(1521.20\%)}
                & 1.46524e-01\textbf{(1930.47\%)} \\
{}              & 1.38348e-04(0.10\%)  
                & 3.32665e-03(143.88\%)
                & 2.37125e-02\textbf{(781.65\%)}
                & 1.78294e-01\textbf{(2661.83\%)}
                & 3.17199e-01\textbf{(2508.69\%)} \\
MOP-WENO-M      & 1.56466e-04(85.63\%)  
                & 2.88442e-03(245.13\%)
                & 5.11795e-03\textbf{(209.14\%)}
                & 9.09352e-03\textbf{(129.72\%)}
                & 1.75990e-02\textbf{(142.84\%)} \\
{}              & 1.63200e-04(99.24\%)  
                & 3.40815e-03(320.93\%)
                & 4.87507e-03\textbf{(204.83\%)}
                & 8.61108e-03\textbf{(123.78\%)}
                & 1.70893e-02\textbf{(136.82\%)} \\
{}              & 5.08956e-04(268.26\%)  
                & 1.01393e-02(643.33\%)
                & 1.02172e-02\textbf{(279.89\%)}
                & 1.98022e-02\textbf{(206.74\%)}
                & 4.01776e-02\textbf{(230.43\%)} \\
WENO-IM(2, 0.1) & 8.46989e-05(0.49\%)  
                & 8.39425e-04(0.44\%)
                & 1.67834e-03\textbf{(1.38\%)}
                & 4.17514e-03\textbf{(5.47\%)}
                & 8.13666e-03\textbf{(12.27\%)} \\
{}              & 8.20061e-05(0.12\%)  
                & 8.10915e-04(0.15\%)
                & 1.60667e-03\textbf{(0.46\%)}
                & 3.93647e-03\textbf{(2.30\%)}
                & 7.71028e-03\textbf{(6.85\%)} \\
{}              & 1.38220e-04(0.01\%)  
                & 1.36420e-03(0.01\%)
                & 2.68977e-03\textbf{(0.01\%)}
                & 6.45231e-03\textbf{(-0.05\%)}
                & 1.21388e-02\textbf{(-0.17\%)} \\
MOP-WENO-IM(2, 0.1)
                & 1.55777e-04(84.82\%) 
                & 2.74109e-03(227.98\%)
                & 4.16210e-03\textbf{(151.40\%)}
                & 8.37898e-03\textbf{(111.67\%)}
                & 1.25166e-02\textbf{(72.71\%)} \\
{}              & 1.62626e-04(98.54\%)  
                & 3.25705e-03(302.26\%)
                & 3.76225e-03\textbf{(135.25\%)}
                & 8.02344e-03\textbf{(108.51\%)}
                & 1.14682e-02\textbf{(58.92\%)} \\
{}              & 5.08361e-04(267.83\%)  
                & 9.88287e-03(624.53\%)
                & 6.81406e-03\textbf{(153.35\%)}
                & 1.84998e-02\textbf{(186.57\%)}
                & 2.02754e-02\textbf{(66.75\%)} \\
WENO-PM6        & 8.40259e-05(-0.31\%)  
                & 8.30374e-04(-0.64\%)
                & 1.63963e-03\textbf{(-0.96\%)}  
                & 3.88864e-03\textbf{(-1.76\%)} 
                & 7.17606e-03\textbf{(-0.98\%)} \\
{}              & 8.19676e-05(0.07\%)  
                & 8.09152e-04(-0.07\%)
                & 1.59697e-03\textbf{(-0.15\%)}  
                & 3.83159e-03\textbf{(-0.43\%)} 
                & 7.19008e-03\textbf{(-0.36\%)} \\
{}              & 1.38205e-04(0.00\%)  
                & 1.36410e-03(0.00\%)
                & 2.68938e-03\textbf{(-0.01\%)}  
                & 6.45650e-03\textbf{(0.01\%)} 
                & 1.21637e-02\textbf{(0.04\%)} \\
MOP-WENO-PM6    & 1.53937e-04(82.63\%)  
                & 2.70283e-03(223.40\%)
                & 4.07454e-03\textbf{(146.11\%)}
                & 8.46326e-03\textbf{(113.80\%)}
                & 1.54196e-02\textbf{(112.77\%)} \\
{}              & 1.59169e-04(94.32\%)  
                & 3.19156e-03(294.18\%)
                & 3.66635e-03\textbf{(129.25\%)}
                & 8.02943e-03\textbf{(108.66\%)}
                & 1.49590e-02\textbf{(107.30\%)} \\
{}              & 4.92116e-04(256.08\%)  
                & 9.52154e-03(598.04\%)
                & 6.49923e-03\textbf{(141.65\%)}
                & 1.83171e-02\textbf{(183.74\%)}
                & 3.15065e-02\textbf{(159.11\%)} \\               
WENO-PPM5       & 8.40198e-05(-0.32\%)  
                & 8.30119e-04(-0.67\%)
                & 1.63931e-03\textbf{(-0.98\%)}
                & 3.89396e-03\textbf{(-1.63\%)}
                & 7.20573e-03\textbf{(-0.57\%)} \\
{}              & 8.19609e-05(0.06\%)  
                & 8.09118e-04(-0.07\%)
                & 1.59692e-03\textbf{(-0.15\%)}
                & 3.83250e-03\textbf{(-0.40\%)}
                & 7.19622e-03\textbf{(-0.28\%)} \\
{}              & 1.38206e-04(0.00\%)  
                & 1.36411e-03(0.01\%)
                & 2.68939e-03\textbf{(-0.01\%)}
                & 6.45658e-03\textbf{(0.01\%)}
                & 1.21629e-02\textbf{(0.03\%)} \\
MOP-WENO-PPM5   & 1.53322e-04(81.90\%)  
                & 2.70476e-03(223.63\%)
                & 4.17894e-03\textbf{(152.42\%)}
                & 8.34997e-03\textbf{(110.94\%)}
                & 1.21149e-02\textbf{(67.17\%)} \\
{}              & 1.59164e-04(94.31\%)  
                & 3.20725e-03(296.11\%)
                & 3.78313e-03\textbf{(136.55\%)}
                & 7.98345e-03\textbf{(107.47\%)}
                & 1.09845e-02\textbf{(52.22\%)} \\
{}              & 4.97691e-04(260.11\%)  
                & 9.71919e-03(612.53\%)
                & 6.89990e-03\textbf{(156.54\%)}
                & 1.83470e-02\textbf{(184.20\%)}
                & 1.87607e-02\textbf{(54.29\%)} \\
WENO-RM(260)    & 8.43348e-05(0.06\%)  
                & 8.35534e-04(-0.03\%)
                & 1.65314e-03\textbf{(-0.15\%)}
                & 3.94006e-03\textbf{(-0.47\%)}
                & 7.25689e-03\textbf{(0.13\%)} \\
{}              & 8.19287e-05(0.02\%)  
                & 8.09704e-04(0.00\%)
                & 1.59924e-03\textbf{(0.00\%)}
                & 3.84494e-03\textbf{(-0.08\%)}
                & 7.22386e-03\textbf{(0.11\%)} \\
{}              & 1.38206e-04(0.00\%)  
                & 1.36404e-03(0.00\%)
                & 2.68956e-03\textbf{(0.00\%)}
                & 6.45544e-03\textbf{(0.00\%)}
                & 1.21576e-02\textbf{(-0.01\%)} \\
MOP-WENO-RM(260)& 1.55787e-04(84.83\%)  
                & 2.72147e-03(225.63\%)
                & 4.13179e-03\textbf{(149.57\%)}
                & 8.32505e-03\textbf{(110.31\%)}
                & 1.57577e-02\textbf{(117.43\%)} \\
{}              & 1.62604e-04(98.51\%)  
                & 3.22596e-03(298.42\%)
                & 3.73345e-03\textbf{(133.44\%)}
                & 7.96768e-03\textbf{(107.06\%)}
                & 1.53795e-02\textbf{(113.12\%)} \\
{}              & 5.05390e-04(265.68\%)  
                & 9.74612e-03(614.50\%)
                & 6.71615e-03\textbf{(149.71\%)}
                & 1.83262e-02\textbf{(183.88\%)}
                & 3.30552e-02\textbf{(171.85\%)} \\
WENO-MAIM1      & 8.24623e-05(-2.17\%)  
                & 8.03920e-04(-3.81\%)
                & 1.58626e-03\textbf{(-4.19\%)}
                & 3.77900e-03\textbf{(-4.53\%)}
                & 7.04287e-03\textbf{(-2.82\%)} \\
{}              & 8.18028e-05(-0.13\%)  
                & 8.04763e-04(-0.61\%)
                & 1.58610e-03\textbf{(-0.82\%)}
                & 3.79788e-03\textbf{(-1.30\%)}
                & 7.14409e-03\textbf{(-1.00\%)} \\
{}              & 1.38215e-04(0.01\%)  
                & 1.36392e-03(-0.01\%)
                & 2.68849e-03\textbf{(-0.04\%)}
                & 6.46356e-03\textbf{(0.12\%)}
                & 1.21473e-02\textbf{(-0.10\%)} \\
MOP-WENO-MAIM1  & 9.97376e-05(18.33\%)  
                & 8.16839e-04(-2.26\%)
                & 1.60912e-03\textbf{(-2.81\%)}
                & 6.83393e-03\textbf{(72.64\%)}
                & 1.24817e-02\textbf{(72.23\%)} \\
{}              & 8.85740e-05(8.13\%)  
                & 8.07516e-04(-0.27\%)
                & 1.59351e-03\textbf{(-0.36\%)}
                & 6.96428e-03\textbf{(80.98\%)}
                & 1.20092e-02\textbf{(66.42\%)} \\
{}              & 1.38172e-04(-0.02\%)  
                & 1.36470e-03(0.05\%)
                & 2.68832e-03\textbf{(-0.05\%)}
                & 1.63188e-02\textbf{(152.78\%)}
                & 2.22178e-02\textbf{(82.72\%)} \\
\rowcolor{gray!40}MIP-WENO-ACM$k$ 
                & 8.42873e-05(-)  
                & 8.35747e-04(-)
                & 1.65557e-03(-)  
                & 3.95849e-03(-) 
                & 7.24723e-03(-) \\
\rowcolor{gray!40}{}
                & 8.19107e-05(-)  
                & 8.09679e-04(-)
                & 1.59929e-03(-)  
                & 3.84802e-03(-) 
                & 7.21626e-03(-) \\
\rowcolor{gray!40}{}              
                & 1.38205e-04(-)  
                & 1.36404e-03(-)
                & 2.68955e-03(-)  
                & 6.45564e-03(-) 
                & 1.21593e-02(-) \\
MOP-WENO-ACM$k$ & 1.55900e-04(84.96\%)  
                & 2.72470e-03(226.02\%)  
                & 4.11740e-03\textbf{(148.70\%)}  
                & 8.34435e-03\textbf{(110.80\%)}  
                & 1.54830e-02\textbf{(113.64\%)}   \\
{}              & 1.63558e-04(99.68\%)  
                & 3.23726e-03(299.82\%)  
                & 3.71649e-03\textbf{(132.38\%)}  
                & 7.96980e-03\textbf{(107.11\%)}  
                & 1.50017e-02\textbf{(107.89\%)}   \\
{}              & 5.22964e-04(278.40\%)  
                & 9.83147e-03(620.76\%)  
                & 6.66166e-03\textbf{(147.69\%)}  
                & 1.83215e-02\textbf{(183.81\%)}  
                & 3.16523e-02\textbf{(160.31\%)}   \\
\hline
\end{tabular*}
\end{myFontSize}
\end{table}

\begin{table}[ht]
\begin{myFontSize}
\centering
\caption{Performance of various considered schemes solving $u_{t} + 
u_{x} = 0$ with $u(x, 0) = \sin^{9} (\pi x), \Delta x= 1/800$.}
\label{table_LAE4:N800}
\begin{tabular*}{\hsize}
{@{}@{\extracolsep{\fill}}lllllll@{}}
\hline
Scheme & $t=10$ & $t=100$ & $t=200$ & $t=500$ & $t=1000$ \\
\hline
WENO-JS         & 4.23531e-07(411.02\%)  
                & 4.74028e-06(471.88\%)
                & 7.29285e-05\textbf{(4299.06\%)}
                & 3.11698e-02\textbf{(751974.43\%)}
                & 1.01278e-01\textbf{(1221783.34\%)} \\
{}              & 3.76804e-07(367.05\%)  
                & 4.32403e-06(435.83\%)
                & 1.60499e-04\textbf{(9844.24\%)}
                & 4.08456e-02\textbf{(1012202.60\%)}
                & 1.13316e-01\textbf{(1404171.46\%)} \\
{}              & 6.95290e-07(410.60\%)  
                & 1.09481e-05(703.79\%)
                & 9.51604e-04\textbf{(34832.14\%)} 
                & 8.63989e-02\textbf{(1268572.78\%)}
                & 2.13485e-01\textbf{(1567406.64\%)} \\              
WENO-M          & 8.28912e-08(0.01\%)  
                & 8.29015e-07(0.01\%)
                & 2.27991e-06\textbf{(37.52\%)}
                & 1.41413e-03\textbf{(34020.56\%)}
                & 1.83325e-02\textbf{(221075.14\%)} \\
{}              & 8.06774e-08(0.00\%)  
                & 8.06977e-07(0.00\%)
                & 2.59031e-06\textbf{(60.49\%)}
                & 3.28891e-03\textbf{(81411.16\%)}
                & 3.30753e-02\textbf{(409786.51\%)} \\
{}              & 1.36173e-07(0.00\%)  
                & 1.36207e-06(0.00\%)
                & 1.22731e-05\textbf{(350.53\%)}
                & 1.90785e-02\textbf{(280046.78\%)}
                & 1.38215e-01\textbf{(1014739.13\%)} \\ 
MOP-WENO-M      & 8.48762e-08(2.41\%)  
                & 9.93577e-07(19.87\%)
                & 1.81123e-06\textbf{(9.25\%)}
                & 4.68314e-06\textbf{(13.00\%)}
                & 8.53126e-06\textbf{(2.93\%)} \\
{}              & 8.11503e-08(0.59\%)  
                & 9.29166e-07(15.14\%)
                & 1.65928e-06\textbf{(2.81\%)}
                & 4.27588e-06\textbf{(5.97\%)}
                & 8.12198e-06\textbf{(0.65\%)} \\
{}              & 1.36173e-07(0.00\%)  
                & 2.03738e-06(49.58\%)
                & 2.72417e-06\textbf{(0.00\%)}
                & 6.81022e-06\textbf{(0.00\%)}
                & 1.36195e-05\textbf{(0.00\%)} \\
WENO-IM(2, 0.1) & 8.28803e-08(0.00\%)  
                & 8.28891e-07(0.00\%)
                & 1.65781e-06\textbf{(0.00\%)}
                & 4.14443e-06\textbf{(0.00\%)}
                & 8.28840e-06\textbf{(0.00\%)} \\
{}              & 8.06769e-08(0.00\%)  
                & 8.06974e-07(0.00\%)
                & 1.61399e-06\textbf{(0.00\%)}
                & 4.03492e-06\textbf{(0.00\%)}
                & 8.06939e-06\textbf{(0.00\%)} \\
{}              & 1.36172e-07(0.00\%)  
                & 1.36206e-06(0.00\%)
                & 2.72415e-06\textbf{(0.00\%)}
                & 6.81019e-06\textbf{(0.00\%)}
                & 1.36194e-05\textbf{(0.00\%)} \\
MOP-WENO-IM(2, 0.1)
                & 8.48292e-08(2.35\%) 
                & 9.80868e-07(18.33\%)
                & 1.79137e-06\textbf{(8.06\%)}
                & 4.88306e-06\textbf{(17.82\%)}
                & 8.63424e-06\textbf{(4.17\%)} \\
{}              & 8.11341e-08(0.57\%)  
                & 9.16723e-07(13.60\%)
                & 1.65133e-06\textbf{(2.31\%)}
                & 4.51320e-06\textbf{(11.85\%)}
                & 8.15362e-06\textbf{(1.04\%)} \\
{}              & 1.36172e-07(0.00\%)  
                & 1.87953e-06(37.99\%)
                & 2.72415e-06\textbf{(0.00\%)}
                & 9.14624e-06\textbf{(34.30\%)}
                & 1.36194e-05\textbf{(0.00\%)} \\
WENO-PM6        & 8.28795e-08(0.00\%)  
                & 8.28892e-07(0.00\%)
                & 1.65782e-06\textbf{(0.00\%)}
                & 4.14452e-06\textbf{(0.00\%)}
                & 8.84565e-06\textbf{(6.72\%)} \\
{}              & 8.06769e-08(0.00\%)  
                & 8.06973e-07(0.00\%)
                & 1.61399e-06\textbf{(0.00\%)}
                & 4.03492e-06\textbf{(0.00\%)}
                & 8.31248e-06\textbf{(3.01\%)} \\
{}              & 1.36172e-07(0.00\%)  
                & 1.36206e-06(0.00\%)
                & 2.72415e-06\textbf{(0.00\%)}
                & 6.81018e-06\textbf{(0.00\%)}
                & 1.38461e-05\textbf{(1.66\%)} \\ 
MOP-WENO-PM6    & 8.47719e-08(2.28\%)  
                & 9.71688e-07(17.23\%)
                & 1.78163e-06\textbf{(7.47\%)}
                & 4.93547e-06\textbf{(19.08\%)}
                & 8.65269e-06\textbf{(4.39\%)} \\
{}              & 8.11105e-08(0.54\%)  
                & 9.05382e-07(12.19\%)
                & 1.64687e-06\textbf{(2.04\%)}
                & 4.61707e-06\textbf{(14.43\%)}
                & 8.15964e-06\textbf{(1.12\%)} \\
{}              & 1.36172e-07(0.00\%)  
                & 1.78452e-06(31.02\%)
                & 2.72415e-06\textbf{(0.00\%)}
                & 1.08735e-05\textbf{(59.67\%)}
                & 1.36194e-05\textbf{(0.00\%)} \\               
WENO-PPM5       & 8.28794e-08(0.00\%)  
                & 8.28890e-07(0.00\%)
                & 1.65781e-06\textbf{(0.00\%)}
                & 4.14448e-06\textbf{(0.00\%)}
                & 8.28862e-06\textbf{(0.00\%)} \\
{}              & 8.06769e-08(0.00\%)  
                & 8.06973e-07(0.00\%)
                & 1.61399e-06\textbf{(0.00\%)}
                & 4.03492e-06\textbf{(0.00\%)}
                & 8.06938e-06\textbf{(0.00\%)} \\
{}              & 1.36172e-07(0.00\%)  
                & 1.36206e-06(0.00\%)
                & 2.72415e-06\textbf{(0.00\%)}
                & 6.81018e-06\textbf{(0.00\%)}
                & 1.36194e-05\textbf{(0.00\%)} \\
MOP-WENO-PPM5   & 8.47367e-08(2.24\%)  
                & 1.04103e-06(25.59\%)
                & 1.83725e-06\textbf{(10.82\%)}
                & 4.30721e-06\textbf{(3.93\%)}
                & 8.27506e-06\textbf{(-0.16\%)} \\
{}              & 8.10958e-08(0.52\%)  
                & 1.00371e-06(24.38\%)
                & 1.67934e-06\textbf{(4.05\%)}
                & 4.06082e-06\textbf{(0.64\%)}
                & 8.07071e-06\textbf{(0.02\%)} \\
{}              & 1.36172e-07(0.00\%)  
                & 1.78285e-06(30.89\%)
                & 2.72415e-06\textbf{(0.00\%)}
                & 6.81018e-06\textbf{(0.00\%)}
                & 1.36194e-05\textbf{(0.00\%)} \\ 
WENO-RM(260)    & 8.28794e-08(0.00\%)  
                & 8.28889e-07(0.00\%)
                & 1.65781e-06\textbf{(0.00\%)}
                & 4.14448e-06\textbf{(0.00\%)}
                & 8.28860e-06\textbf{(0.00\%)} \\
{}              & 8.06769e-08(0.00\%)  
                & 8.06973e-07(0.00\%)
                & 1.61399e-06\textbf{(0.00\%)}
                & 4.03492e-06\textbf{(0.00\%)}
                & 8.06938e-06\textbf{(0.00\%)} \\
{}              & 1.36172e-07(0.00\%)  
                & 1.36206e-06(0.00\%)
                & 2.72415e-06\textbf{(0.00\%)}
                & 6.81018e-06\textbf{(0.00\%)}
                & 1.36194e-05\textbf{(0.00\%)} \\
MOP-WENO-RM(260)& 8.48225e-08(2.34\%)  
                & 9.56819e-07(15.43\%)
                & 1.77008e-06\textbf{(6.77\%)}
                & 4.72311e-06\textbf{(13.96\%)}
                & 8.55573e-06\textbf{(3.22\%)} \\
{}              & 8.11306e-08(0.56\%)  
                & 8.87179e-07(9.94\%)
                & 1.64102e-06\textbf{(1.67\%)}
                & 4.33993e-06\textbf{(7.56\%)}
                & 8.12608e-06\textbf{(0.70\%)} \\
{}              & 1.36172e-07(0.00\%)  
                & 1.58577e-06(16.42\%)
                & 2.72415e-06\textbf{(0.00\%)}
                & 6.81018e-06\textbf{(0.00\%)}
                & 1.36194e-05\textbf{(0.00\%)} \\      
WENO-MAIM1      & 8.28796e-08(0.00\%)  
                & 8.28893e-07(0.00\%)
                & 1.65782e-06\textbf{(0.00\%)}
                & 4.14450e-06\textbf{(0.00\%)}
                & 8.28865e-06\textbf{(0.00\%)} \\
{}              & 8.06776e-08(0.00\%)  
                & 8.06974e-07(0.00\%)
                & 1.61399e-06\textbf{(0.00\%)}
                & 4.03492e-06\textbf{(0.00\%)}
                & 8.06938e-06\textbf{(0.00\%)} \\
{}              & 1.36172e-07(0.00\%)  
                & 1.36206e-06(0.00\%)
                & 2.72415e-06\textbf{(0.00\%)}
                & 6.81018e-06\textbf{(0.00\%)}
                & 1.36194e-05\textbf{(0.00\%)} \\
MOP-WENO-MAIM1  & 8.28791e-08(0.00\%)  
                & 8.28894e-07(0.00\%)
                & 1.65783e-06\textbf{(0.00\%)}
                & 4.14454e-06\textbf{(0.00\%)}
                & 8.28830e-06\textbf{(0.00\%)} \\
{}              & 8.06770e-08(0.00\%)  
                & 8.06973e-07(0.00\%)
                & 1.61399e-06\textbf{(0.00\%)}
                & 4.03491e-06\textbf{(0.00\%)}
                & 8.06937e-06\textbf{(0.00\%)} \\
{}              & 1.36172e-07(0.00\%)  
                & 1.36206e-06(0.00\%)
                & 2.72415e-06\textbf{(0.00\%)}
                & 6.81018e-06\textbf{(0.00\%)}
                & 1.36194e-05\textbf{(0.00\%)} \\
\rowcolor{gray!40}MIP-WENO-ACM$k$ 
                & 8.28794e-08(-)  
                & 8.28891e-07(-)
                & 1.65782e-06(-)  
                & 4.14451e-06(-) 
                & 8.28868e-06(-) \\
\rowcolor{gray!40}{}
                & 8.06769e-08(-)  
                & 8.06973e-07(-)
                & 1.61399e-06(-)  
                & 4.03492e-06(-) 
                & 8.06938e-06(-) \\
\rowcolor{gray!40}{}
                & 1.36172e-07(-)  
                & 1.36206e-06(-)
                & 2.72415e-06(-)  
                & 6.81018e-06(-) 
                & 1.36194e-05(-) \\
MOP-WENO-ACM$k$ & 8.47930e-08(2.31\%)  
                & 9.73202e-07(17.41\%)  
                & 1.78369e-06\textbf{(7.59\%)}
                & 4.84739e-06\textbf{(16.96\%)}
                & 8.61232e-06\textbf{(3.90\%)} \\
{}              & 8.11193e-08(0.55\%)  
                & 9.07161e-07(12.42\%)  
                & 1.64768e-06\textbf{(2.09\%)}
                & 4.47345e-06\textbf{(10.87\%)}
                & 8.14436e-06\textbf{(0.93\%)} \\
{}              & 1.36172e-07(0.00\%)  
                & 1.79160e-06(31.54\%)  
                & 2.72415e-06\textbf{(0.00\%)}
                & 8.79296e-06\textbf{(29.11\%)}
                & 1.36194e-05\textbf{(0.00\%)} \\                                
\hline
\end{tabular*}
\end{myFontSize}
\end{table}

\begin{example}
\rm{We calculate Eq.(\ref{eq:LAE}) using the following initial 
condition \cite{WENO-JS}}
\label{LAE3}
\end{example} 
\begin{equation}
\begin{array}{l}
u(x, 0) = \left\{
\begin{array}{ll}
\dfrac{1}{6}\big[ G(x, \beta, z - \hat{\delta}) + 4G(x, \beta, z) + G
(x, \beta, z + \hat{\delta}) \big], & x \in [-0.8, -0.6], \\
1, & x \in [-0.4, -0.2], \\
1 - \big\lvert 10(x - 0.1) \big\rvert, & x \in [0.0, 0.2], \\
\dfrac{1}{6}\big[ F(x, \alpha, a - \hat{\delta}) + 4F(x, \alpha, a) +
 F(x, \alpha, a + \hat{\delta}) \big], & x \in [0.4, 0.6], \\
0, & \mathrm{otherwise},
\end{array}\right. 
\end{array}
\label{eq:LAE:SLP}
\end{equation}
where $G(x, \beta, z) = \mathrm{e}^{-\beta (x - z)^{2}}, F(x, \alpha
, a) = \sqrt{\max \big(1 - \alpha ^{2}(x - a)^{2}, 0 \big)}$, and 
the constants are $z = -0.7, \hat{\delta} = 0.005, \beta = \dfrac{
\log 2}{36\hat{\delta} ^{2}}, a = 0.5$ and $\alpha = 10$. The 
periodic boundary condition is used and the CFL number is taken to 
be $0.1$. For brevity in the presentation, we call this 
\textit{\textbf{L}inear \textbf{P}roblem} SLP as it is presented by 
\textit{\textbf{S}hu} et al. in \cite{WENO-JS}. It is known that 
this problem consists of a Gaussian, a square wave, a sharp triangle 
and a semi-ellipse. 

In Table \ref{table_LAE3}, we present the $L_{1}, L_{2}, L_{\infty}$ 
errors and the corresponding convergence rates of accuracy with 
$t = 2, 2000$. For the case of $t = 2$, it can be seen that: (1) the 
$L_{1}$ and $L_{2}$ orders of all considered schemes are 
approximately $1.0$ and about $0.35$ to $0.5$, respectively; (2) 
negative values of the $L_{\infty}$ orders of all considered schemes 
are generated; (3) in terms of accuracy, the MOP-WENO-X schemes 
produce less accurate results than the corresponding WENO-X schemes. 
For the case of $t=2000$, it can be seen that: (1) the $L_{1}$,
$L_{2}$ orders of the WENO-JS and WENO-M schemes decrease to very 
small values and even become negative; (2) however, the $L_{1}$ and
$L_{2}$ orders of all the MOP-WENO-X schemes, as well as the 
associated mapped WENO-X schemes without WENO-M, are clearly 
larger than $1.0$ and around $0.5$ to $0.9$ respectively; (3) the
$L_{\infty}$ orders of all WENO-X schemes are very small and some of 
then even become negative (e.g., the WENO-JS, WENO-PPM5 and 
MIP-WENO-ACM$k$ schemes), while those of the MOP-WENO-X schemes are 
all positive although they are also very small; (4) in terms of 
accuracy, on the whole, the MOP-WENO-X schemes produce accurate and 
comparable results as the corresponding WENO-X schemes except the 
WENO-M scheme. However, if we take a closer look, we can find that 
the resolution of the result computed by the WENO-M scheme is 
significantly lower than that of the MOP-WENO-M scheme, and the 
other mapped WENO-X schemes generate spurious oscillations but the 
corresponding MOP-WENO-X schemes do not. Detailed tests will be 
conducted and the solutions will be presented carefully to 
demonstrate this in the following subsection.

\begin{table}[ht]
\begin{myFontSize}
\centering
\caption{Convergence properties of various considered schemes 
solving $u_{t} + u_{x} = 0$ with initial condition 
Eq.(\ref{eq:LAE:SLP}).}
\label{table_LAE3}
\begin{tabular*}{\hsize}
{@{}@{\extracolsep{\fill}}lllllll@{}}
\hline
\space  &\multicolumn{3}{l}{t = 2}  &\multicolumn{3}{l}{t = 2000}  \\
\cline{2-4}  \cline{5-7}
$N$
& 200            & 400                  & 800 
& 200            & 400                  & 800   \\
\hline
WENO-JS     
& 6.30497e-02(-) & 2.81654e-02(1.2103)  & 1.41364e-02(0.9945)  
& 6.12899e-01(-) & 5.99215e-01(0.0326)  & 5.50158e-01(0.1232)   \\
\space
& 1.08621e-01(-) & 7.71111e-02(0.4943)  & 5.69922e-02(0.4362)  
& 5.08726e-01(-) & 5.01160e-01(0.0216)  & 4.67585e-01(0.1000)   \\
\space
& 4.09733e-01(-) & 4.19594e-01(-0.0343) & 4.28463e-01(-0.0302)  
& 7.99265e-01(-) & 8.20493e-01(-0.0378) & 8.14650e-01(0.0103)   \\
WENO-M     
& 4.77201e-02(-) & 2.23407e-02(1.0949)  & 1.11758e-02(0.9993)  
& 3.81597e-01(-) & 3.25323e-01(0.2302)  & 3.48528e-01(-0.0994)  \\
\space
& 9.53073e-02(-) & 6.91333e-02(0.4632)  & 5.09232e-02(0.4411)  
& 3.59205e-01(-) & 3.12970e-01(0.1988)  & 3.24373e-01(-0.0516)  \\
\space
& 3.94243e-01(-) & 4.05856e-01(-0.0419) & 4.16937e-01(-0.0389)  
& 6.89414e-01(-) & 6.75473e-01(0.0295)  & 6.25645e-01(0.1106)   \\
MOP-WENO-M     
& 5.72690e-02(-) & 2.72999e-02(1.0689)  & 1.42908e-02(0.9338)  
& 3.85134e-01(-) & 1.74987e-01(1.1381)  & 6.40251e-02(1.4505)  \\
\space
& 1.00827e-01(-) & 7.33765e-02(0.4585)  & 5.57886e-02(0.3953)  
& 3.48164e-01(-) & 1.86418e-01(0.9012)  & 1.07629e-01(0.7925)  \\
\space
& 4.14785e-01(-) & 4.45144e-01(-0.1019) & 4.64024e-01(-0.0599)  
& 7.41230e-01(-) & 5.04987e-01(0.5537)  & 4.81305e-01(0.0693)  \\
WENO-IM(2,0.1)     
& 4.40293e-02(-) & 2.02331e-02(1.1217)  & 1.01805e-02(0.9909)  
& 2.17411e-01(-) & 1.12590e-01(0.9493)  & 5.18367e-02(1.1190)   \\
\space
& 9.19118e-02(-) & 6.68479e-02(0.4594)  & 4.95333e-02(0.4325)  
& 2.30000e-01(-) & 1.64458e-01(0.4839)  & 9.98968e-02(0.7192)   \\
\space
& 3.86789e-01(-) & 3.98769e-01(-0.0441) & 4.09515e-01(-0.0383)  
& 5.69864e-01(-) & 4.82180e-01(0.2410)  & 4.73102e-01(0.02784)   \\
MOP-WENO-IM(2,0.1)     
& 6.09985e-02(-) & 2.86731e-02(1.0891)  & 1.45601e-02(0.9777)  
& 3.83289e-01(-) & 1.67452e-01(1.1947)  & 6.44253e-02(1.3780)  \\
\space
& 1.03438e-01(-) & 7.56598e-02(0.4512)  & 5.61842e-02(0.4294)  
& 3.47817e-01(-) & 1.76550e-01(0.9783)  & 1.05858e-01(0.7379)  \\
\space
& 4.35238e-01(-) & 4.62098e-01(-0.0864) & 4.64674e-01(-0.0080)  
& 7.25185e-01(-) & 5.24538e-01(0.4673)  & 5.19333e-01(0.0144)  \\
WENO-PM6     
& 4.66681e-02(-) & 2.13883e-02(1.1256)  & 1.06477e-02(1.0063)  
& 2.17323e-01(-) & 1.05197e-01(1.0467)  & 4.47030e-02(1.2347)   \\
\space
& 9.45566e-02(-) & 6.82948e-02(0.4694)  & 5.03724e-02(0.4391)  
& 2.28655e-01(-) & 1.47518e-01(0.6323)  & 9.34250e-02(0.6590)   \\
\space
& 3.96866e-01(-) & 4.06118e-01(-0.0332) & 4.15277e-01(-0.0322)  
& 5.63042e-01(-) & 5.04977e-01(0.1570)  & 4.71368e-01(0.0994)   \\
MOP-WENO-PM6     
& 5.45129e-02(-) & 2.61755e-02(1.0584)  & 1.38981e-02(0.9133)  
& 4.51487e-01(-) & 1.75875e-01(1.3601)  & 6.32990e-02(1.4743)  \\
\space
& 9.95654e-02(-) & 7.16656e-02(0.4744)  & 5.44733e-02(0.3957)  
& 4.01683e-01(-) & 1.83478e-01(1.1305)  & 1.04688e-01(0.8095)  \\
\space
& 4.02785e-01(-) & 4.26334e-01(-0.0820) & 4.63134e-01(-0.1194)  
& 7.71539e-01(-) & 5.06314e-01(0.6077)  & 4.76091e-01(0.0888)  \\
WENO-PPM5     
& 4.54081e-02(-) & 2.07948e-02(1.1267)  & 1.04018e-02(0.9994)  
& 2.17174e-01(-) & 1.03201e-01(1.0734)  & 4.81637e-02(1.0994)  \\
\space
& 9.33165e-02(-) & 6.76172e-02(0.4647)  & 4.99580e-02(0.4367)  
& 2.29008e-01(-) & 1.46610e-01(0.6434)  & 9.47748e-02(0.6294)  \\
\space
& 3.91076e-01(-) & 4.02214e-01(-0.0405) & 4.12113e-01(-0.0351)  
& 5.65575e-01(-) & 5.06463e-01(0.1593)  & 5.14402e-01(-0.0224)  \\
MOP-WENO-PPM5     
& 5.51553e-02(-) & 2.65464e-02(1.0550)  & 1.41381e-02(0.9089)  
& 3.86292e-01(-) & 1.75232e-01(1.1404)  & 6.36336e-02(1.4614)  \\
\space
& 9.94592e-02(-) & 7.19973e-02(0.4662)  & 5.52704e-02(0.3814)  
& 3.49072e-01(-) & 1.88491e-01(0.8890)  & 1.06801e-01(0.8196)  \\
\space
& 4.04763e-01(-) & 4.32887e-01(-0.0969) & 4.68577e-01(-0.1143)  
& 7.36405e-01(-) & 5.14732e-01(0.5167)  & 4.98424e-01(0.0464)  \\
WENO-RM(260)     
& 4.63072e-02(-) & 2.13545e-02(1.1167)  & 1.06392e-02(1.0052)  
& 2.17363e-01(-) & 1.04347e-01(1.0587)  & 4.45176e-02(1.2289)  \\
\space
& 9.40674e-02(-) & 6.81954e-02(0.4640)  & 5.03289e-02(0.4383)  
& 2.28662e-01(-) & 1.47093e-01(0.6365)  & 9.33066e-02(0.6567)  \\
\space
& 3.96762e-01(-) & 4.08044e-01(-0.0405) & 4.16722e-01(-0.0304)  
& 5.62933e-01(-) & 4.98644e-01(0.1750)  & 4.71450e-01(0.0809)  \\
MOP-WENO-RM(260)     
& 5.54343e-02(-) & 2.71415e-02(1.0303)  & 1.45563e-02(0.8989)  
& 4.56942e-01(-) & 2.25420e-01(1.0194)  & 8.02414e-02(1.4902)  \\
\space
& 9.93009e-02(-) & 7.22823e-02(0.4582)  & 5.66845e-02(0.3507)  
& 4.06524e-01(-) & 2.25814e-01(0.8482)  & 1.18512e-01(0.9301)  \\
\space
& 4.04041e-01(-) & 4.38358e-01(-0.1176) & 4.70380e-01(-0.1017)  
& 7.71747e-01(-) & 5.12018e-01(0.5919)  & 4.90610e-01(0.0616)  \\
WENO-MAIM1     
& 5.71142e-02(-) & 2.48065e-02(1.2031)  & 1.21078e-02(1.0348)  
& 2.18238e-01(-) & 1.09902e-01(0.9897)  & 4.41601e-02(1.3154)  \\
\space
& 1.03257e-01(-) & 7.29236e-02(0.5018)  & 5.32803e-02(0.4528)  
& 2.29151e-01(-) & 1.51024e-01(0.6015)  & 9.35506e-02(0.6910)  \\
\space
& 4.15051e-01(-) & 4.23185e-01(-0.0280) & 4.28710e-01(-0.0187)  
& 5.63682e-01(-) & 4.94657e-01(0.1885)  & 4.72393e-01(0.0664)  \\
MOP-WENO-MAIM1     
& 5.98640e-02(-) & 2.64819e-02(1.1767)  & 1.33647e-02(0.9866)  
& 2.39900e-01(-) & 1.41890e-01(0.7577)  & 5.43475e-02(1.3845)  \\
\space
& 1.05066e-01(-) & 7.38102e-02(0.5094)  & 5.44089e-02(0.4400)  
& 2.47191e-01(-) & 1.71855e-01(0.5244)  & 1.02170e-01(0.7502)  \\
\space
& 4.12365e-01(-) & 4.26841e-01(-0.0498) & 4.38310e-01(-0.0383)  
& 6.06985e-01(-) & 5.61908e-01(0.1113)  & 5.10242e-01(0.1392)  \\
MIP-WENO-ACM$k$     
& 4.45059e-02(-) & 2.03667e-02(1.1278)  & 1.02183e-02(0.9951)  
& 2.21312e-01(-) & 1.10365e-01(1.0038)  & 4.76589e-02(1.2115)   \\
\space
& 9.24356e-02(-) & 6.70230e-02(0.4638)  & 4.96081e-02(0.4341)  
& 2.28433e-01(-) & 1.48498e-01(0.6213)  & 9.40843e-02(0.6584)   \\
\space
& 3.92505e-01(-) & 4.04024e-01(-0.0417) & 4.13511e-01(-0.0335)  
& 5.36242e-01(-) & 5.13503e-01(0.0625)  & 5.15898e-01(-0.0067)   \\
MOP-WENO-ACM$k$     
& 5.56533e-02(-) & 2.79028e-02(0.9961)  & 1.43891e-02(0.9554)  
& 3.83033e-01(-) & 1.77114e-01(1.1128)  & 6.70535e-02(1.4013)   \\
\space
& 9.94223e-02(-) & 7.33101e-02(0.4396)  & 5.51602e-02(0.4104)  
& 3.46814e-01(-) & 1.87369e-01(0.8883)  & 1.09368e-01(0.7767)   \\
\space
& 4.03765e-01(-) & 4.48412e-01(-0.1513) & 4.67036e-01(-0.0587)  
& 7.18464e-01(-) & 5.05980e-01(0.5058)  & 4.80890e-01(0.0734)   \\
\hline
\end{tabular*}
\end{myFontSize}
\end{table}

\subsection{1D linear advection problems with long output times}
The objective of this subsection is to demonstrate the advantage of 
the MOP-WENO-X schemes on long output time simulations that can 
obtain high resolution and meanwhile do not generate spurious 
oscillations. 

The one-dimensional linear advection problem Eq.(\ref{eq:LAE}) is 
solved with the periodic boundary condition by taking the following 
two initial conditions. 

Case 1. (SLP) The initial condition is given by 
Eq.(\ref{eq:LAE:SLP}).

Case 2. (BiCWP) The initial condition is given by
\begin{equation}
\begin{array}{l}
u(x, 0) = \left\{
\begin{array}{ll}
0,   & x \in [-1.0, -0.8] \bigcup (-0.2, 0.2] \bigcup (0.8, 1.0], \\
0.5, & x \in (-0.6, -0.4] \bigcup (0.2, 0.4]  \bigcup (0.6, 0.8], \\
1,   & x \in (-0.8, -0.6] \bigcup (-0.4, -0.2] \bigcup (0.4, 0.6].
\end{array}\right. 
\end{array}
\label{eq:LAE:BiCWP}
\end{equation}

Case 1 and Case 2 were carefully simulated in \cite{MOP-WENO-ACMk}. 
Case 1 is called SLP as mentioned earlier in this paper. Case 2 
consists of several constant states separated by sharp 
discontinuities at $x = \pm 0.8, \pm 0.6, \pm 0.4, \pm 0.2$ and it 
was called BiCWP for brevity in the presentation as the profile of 
the exact solution for this \textit{\textbf{P}roblem} looks like the 
\textit{\textbf{B}reach \textbf{i}n \textbf{C}ity \textbf{W}all}.

\begin{figure}[ht]
\centering
\includegraphics[height=0.32\textwidth]
{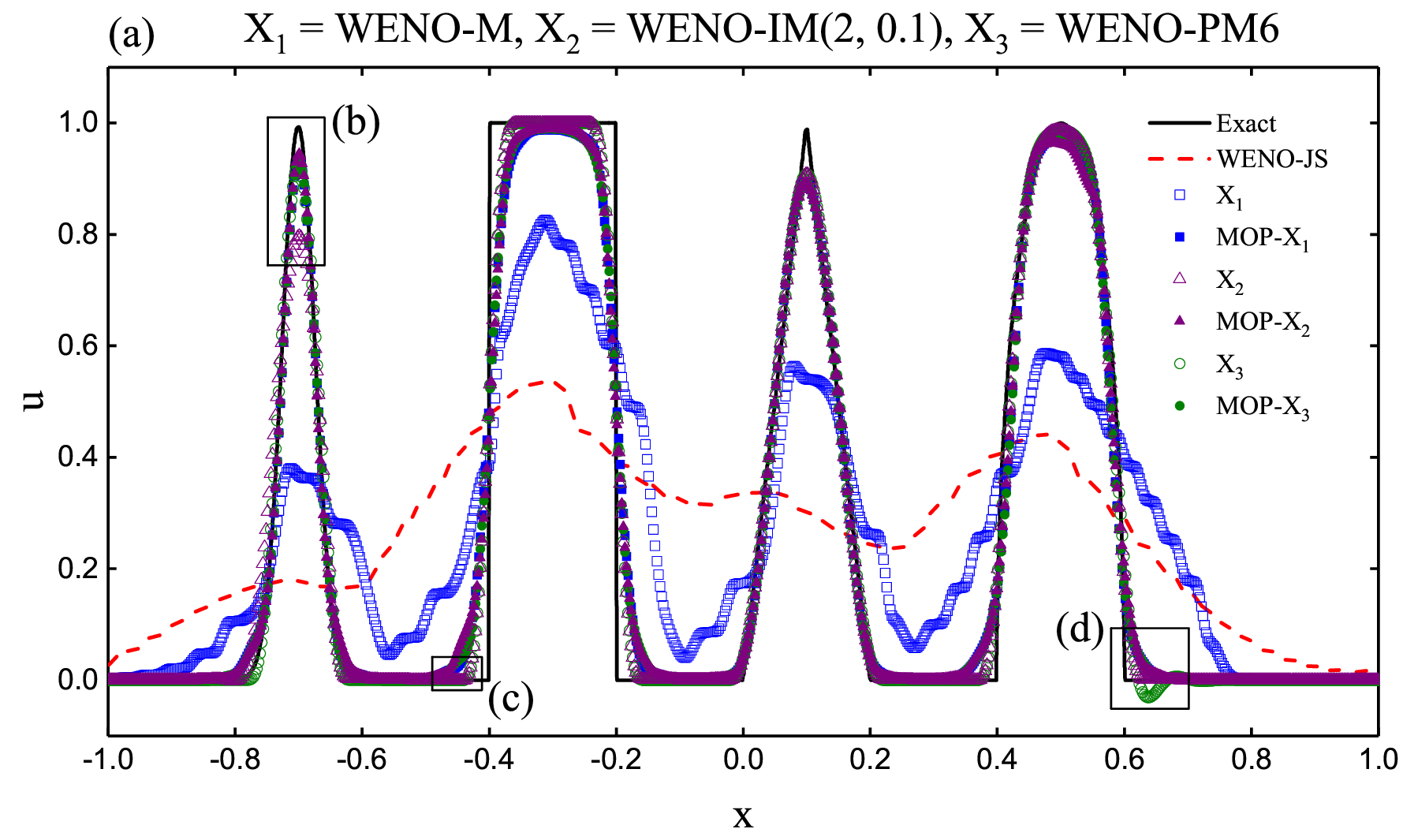}
\includegraphics[height=0.32\textwidth]
{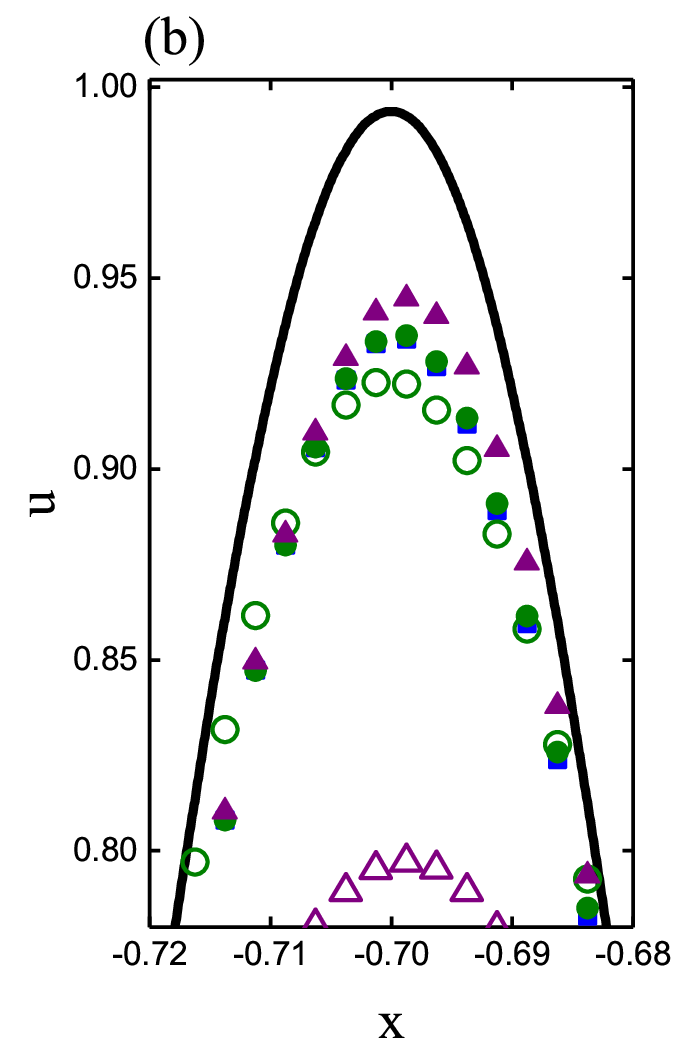}
\includegraphics[height=0.32\textwidth]
{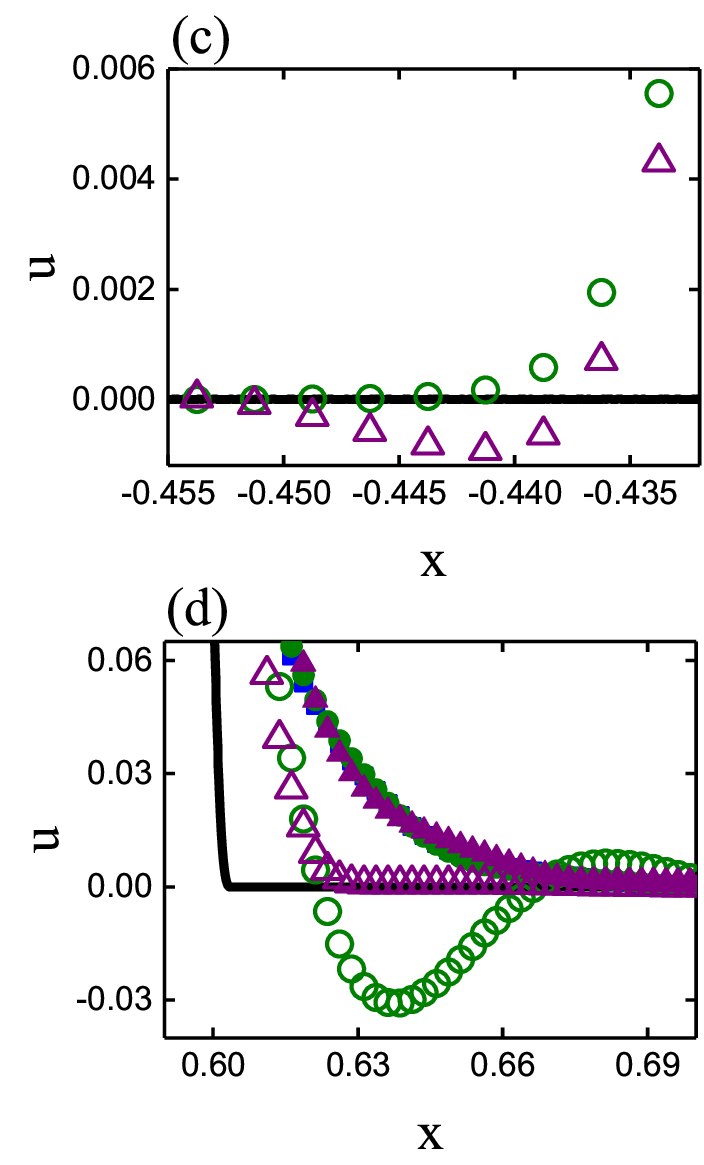}
\caption{Performance of the WENO-JS, WENO-M, MOP-WENO-M, 
WENO-IM($2,0.1$), MOP-WENO-IM($2,0.1$), WENO-PM6 and MOP-WENO-PM6 
schemes for the SLP at output time $t = 2000$ with a uniform mesh 
size of $N = 800$.}
\label{fig:SLP:N800:PM6}
\end{figure}

\begin{figure}[ht]
\centering
\includegraphics[height=0.32\textwidth]
{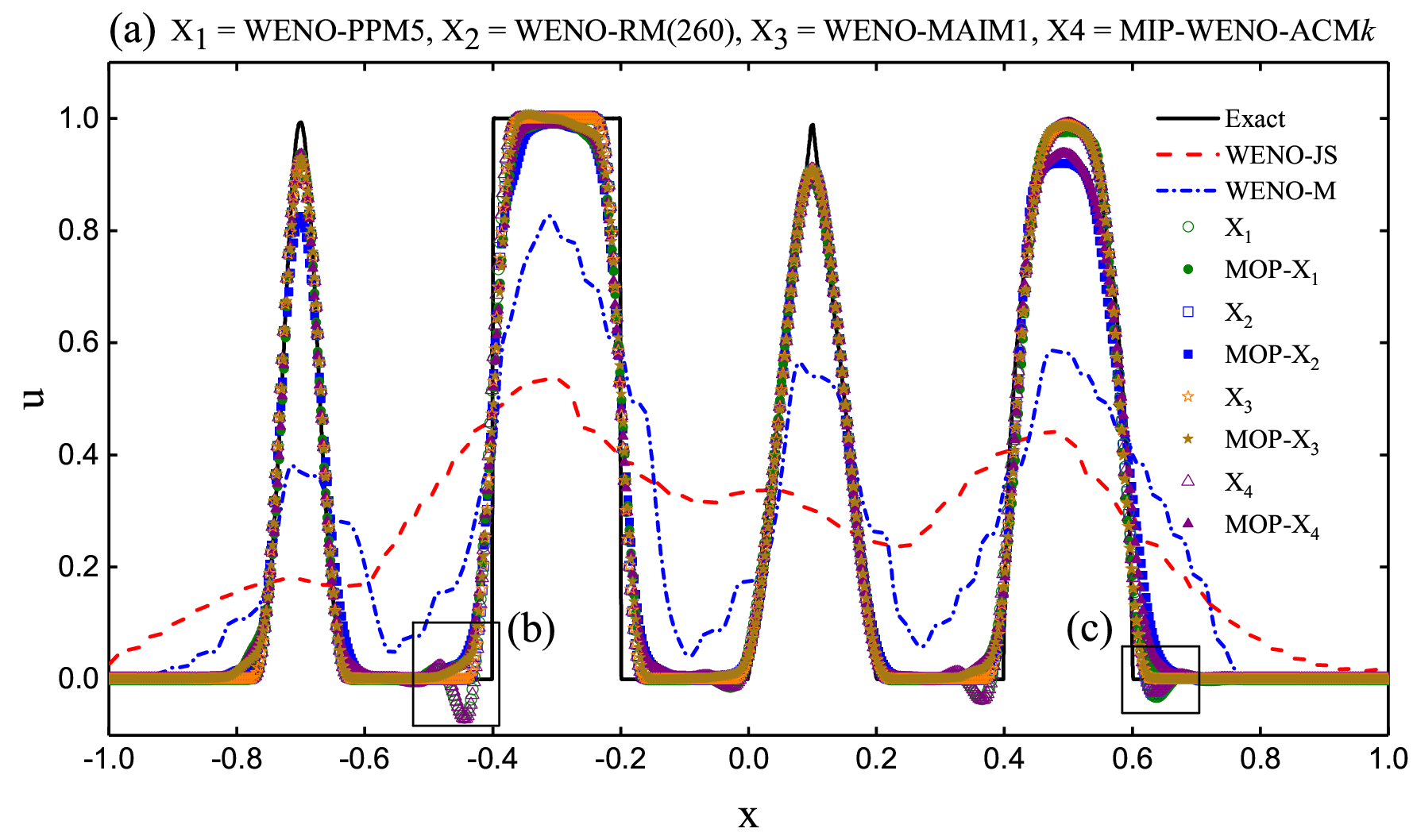}
\includegraphics[height=0.32\textwidth]
{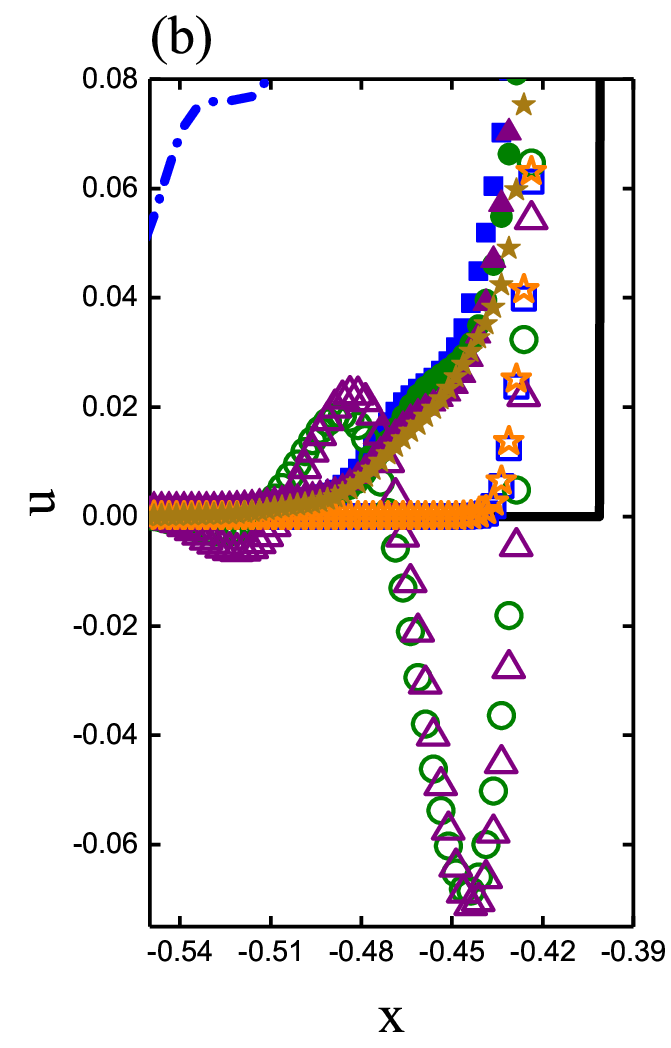}
\includegraphics[height=0.32\textwidth]
{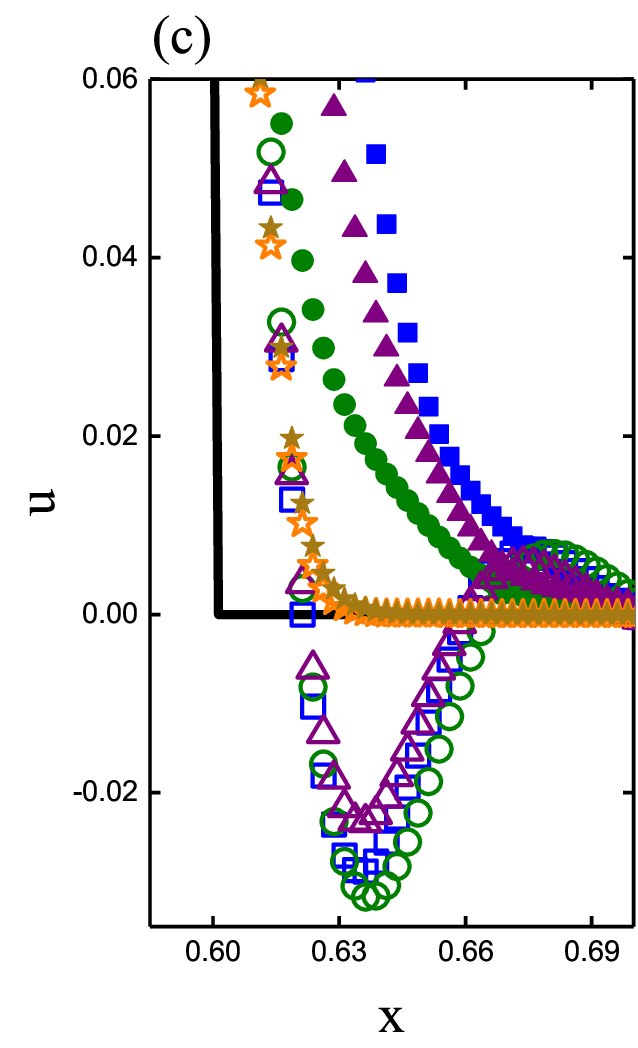}
\caption{Performance of the WENO-JS, WENO-M, WENO-PPM5, 
MOP-WENO-PPM5, WENO-RM260, MOP-WENO-RM260, WENO-MAIM1,
MOP-WNEO-MAIM1, MIP-WENO-ACM$k$ and MOP-WENO-ACM$k$ schemes for the 
SLP at output time $t = 2000$ with a uniform mesh size of $N = 800$.}
\label{fig:SLP:N800:RM260}
\end{figure}

\begin{figure}[ht]
\centering
\includegraphics[height=0.32\textwidth]
{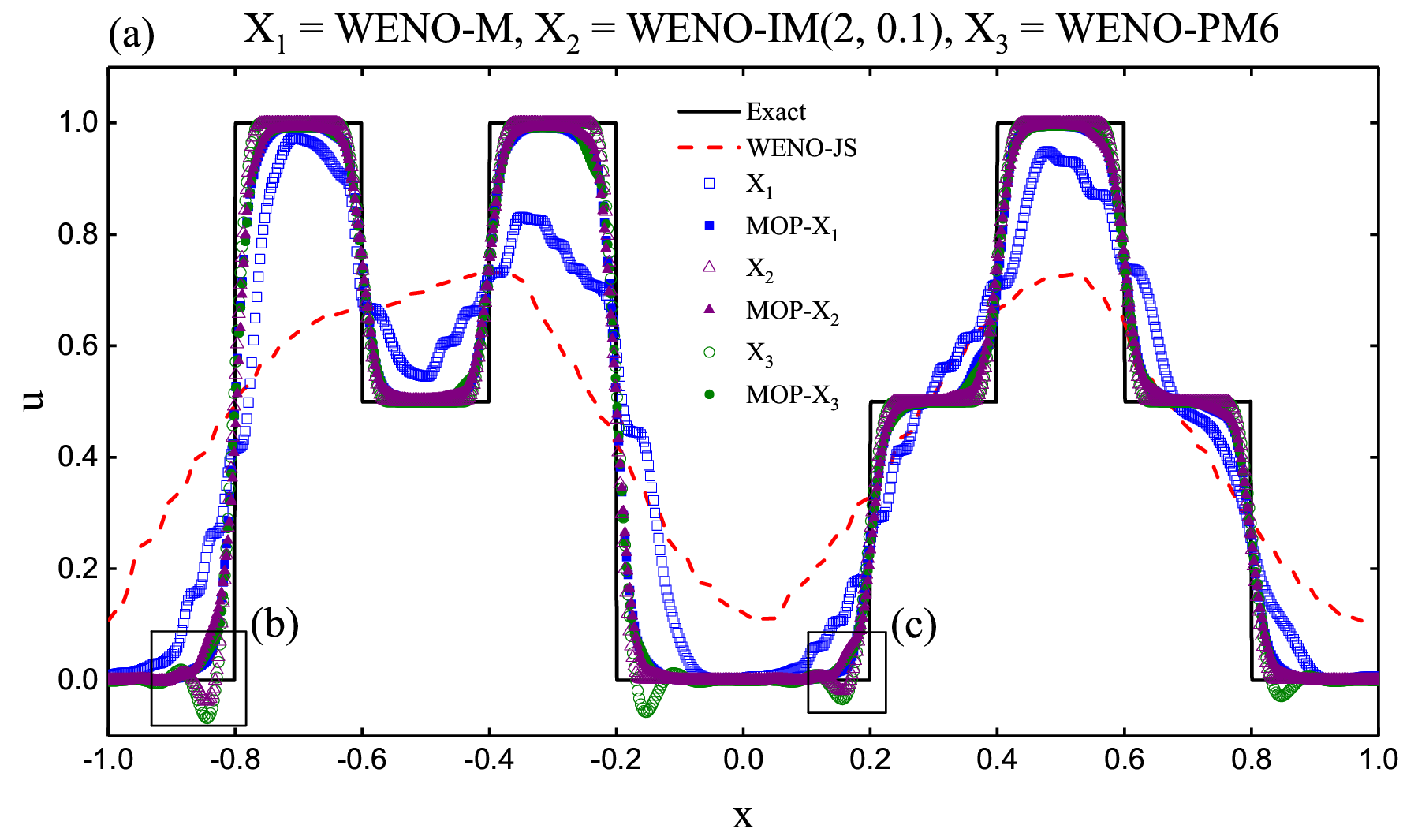}
\includegraphics[height=0.32\textwidth]
{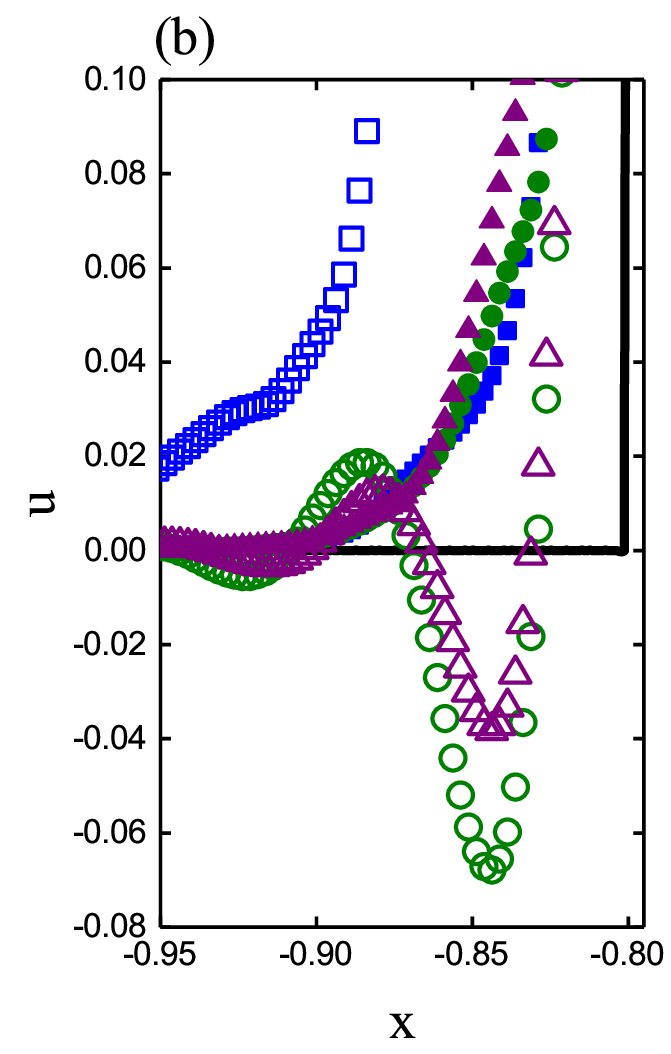}
\includegraphics[height=0.32\textwidth]
{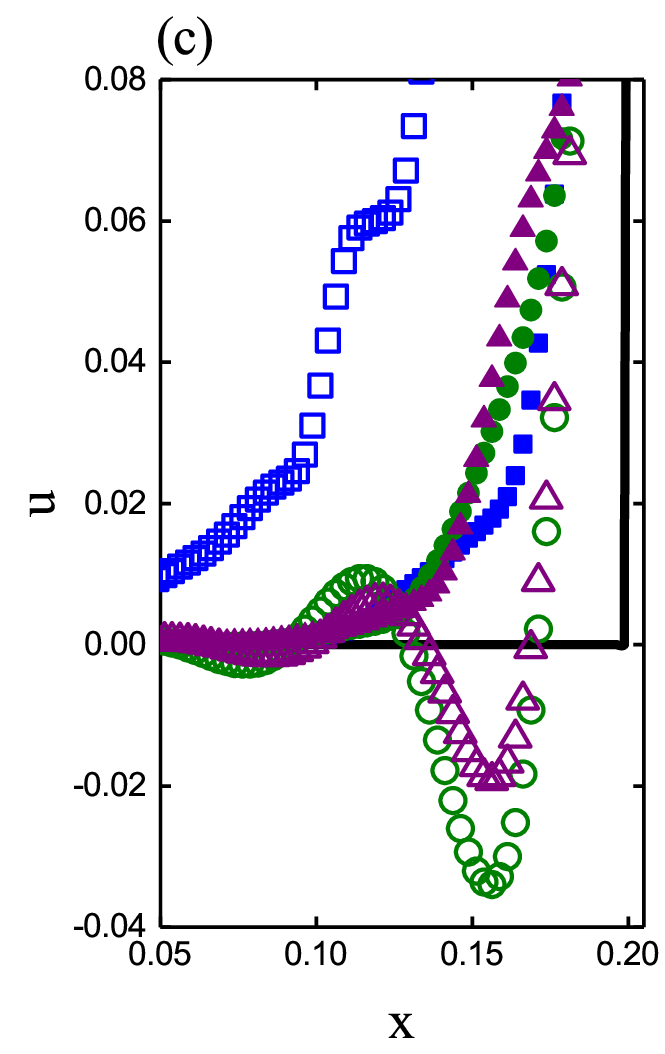}
\caption{Performance of the WENO-JS, WENO-M, MOP-WENO-M, 
WENO-IM($2,0.1$), MOP-WENO-IM($2,0.1$), WENO-PM6 and MOP-WENO-PM6 
schemes for the BiCWP at output time $t = 2000$ with a uniform mesh 
size of $N = 800$.}
\label{fig:BiCWP:N800:PM6}
\end{figure}

\begin{figure}[ht]
\centering
\includegraphics[height=0.32\textwidth]
{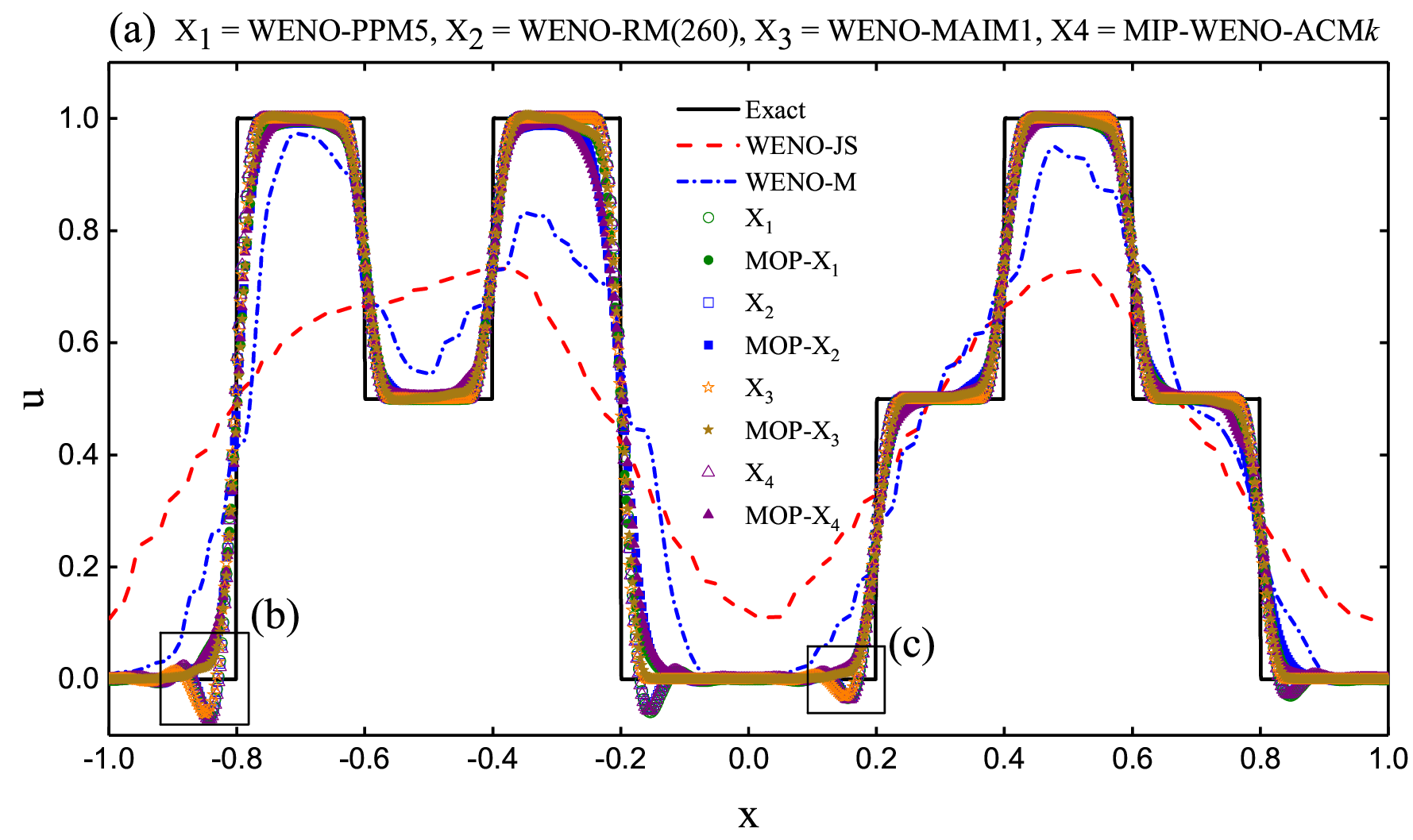}
\includegraphics[height=0.32\textwidth]
{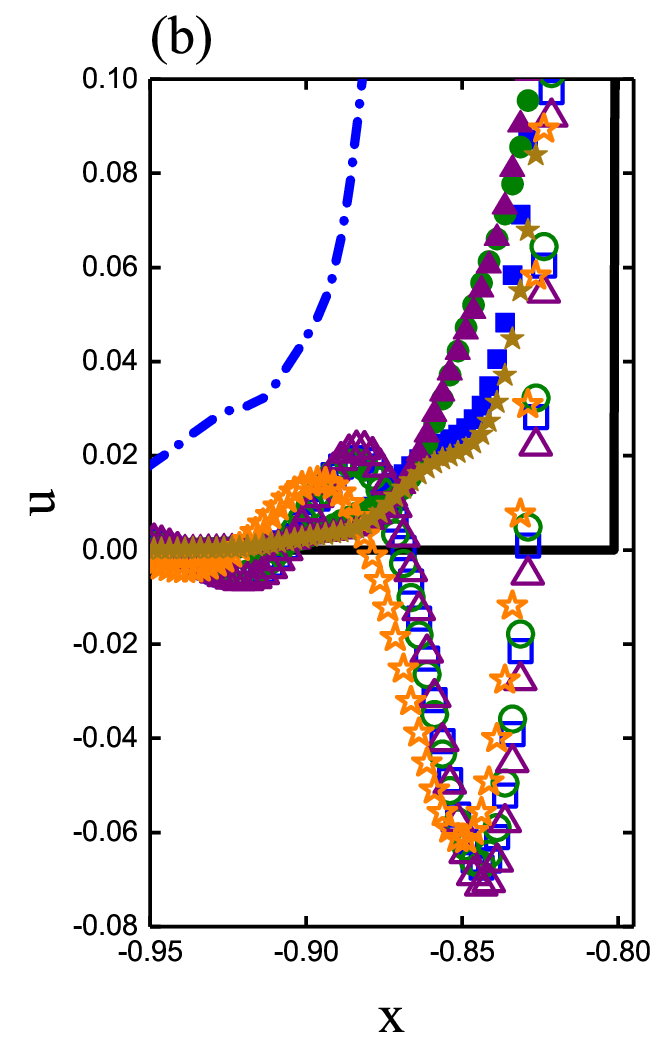}
\includegraphics[height=0.32\textwidth]
{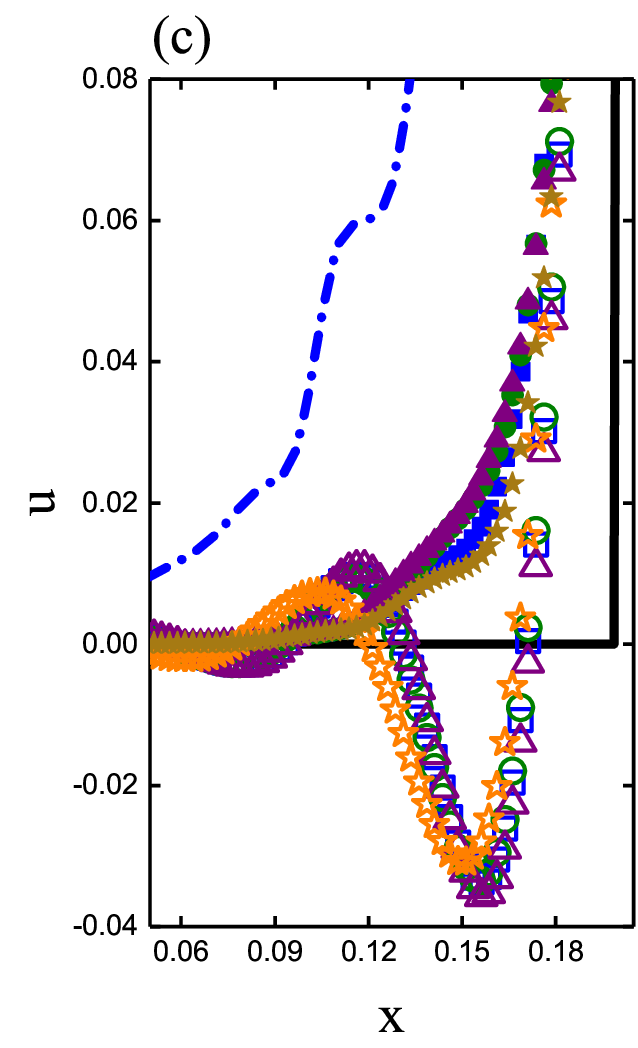}
\caption{Performance of the WENO-JS, WENO-M, WENO-PPM5, 
MOP-WENO-PPM5, WENO-RM260, MOP-WENO-RM260, WENO-MAIM1,
MOP-WNEO-MAIM1, MIP-WENO-ACM$k$ and MOP-WENO-ACM$k$ schemes for the 
BiCWP at output time $t = 2000$ with a uniform mesh size of $N=800$.}
\label{fig:BiCWP:N800:RM260}
\end{figure}

In Figs. \ref{fig:SLP:N800:PM6}, \ref{fig:SLP:N800:RM260} and Figs. 
\ref{fig:BiCWP:N800:PM6}, \ref{fig:BiCWP:N800:RM260}, we show the 
comparison of considered schemes for SLP and BiCWP respectively, by 
taking $t = 2000$ and $N = 800$. It can be seen that: (1) all the 
MOP-WENO-X schemes produce results with considerable resolutions 
which are significantly higher than those of the WENO-JS and WENO-M 
schemes, and what's more, they all do not generate spurious 
oscillations, while most of their corresponding WENO-X schemes do, 
when solving both SLP and BiCWP; (2) it should be reminded that the 
WENO-IM(2, 0.1) scheme appears not to generate spurious oscillations 
and it gives better resolution than the MOP-WENO-IM(2, 0.1) scheme 
in most of the region when solving SLP on present computing 
condition, however, from Fig. \ref{fig:SLP:N800:PM6}(b), one can 
observe that the MOP-WENO-IM(2, 0.1) scheme gives a better 
resolution of the Gaussian than the WENO-IM(2, 0.1) scheme, and if 
taking a closer look, one can see that the WENO-IM(2, 0.1) scheme 
generates a very slight spurious oscillation near $x = - 0.435$ as 
shown in Fig. \ref{fig:SLP:N800:PM6}(c); (3) it is very evident as 
shown in Fig. \ref{fig:BiCWP:N800:PM6} that, when solving BiCWP, the 
WENO-IM(2, 0.1) scheme generates the spurious oscillations. 

\begin{figure}[ht]
\centering
\includegraphics[height=0.32\textwidth]
{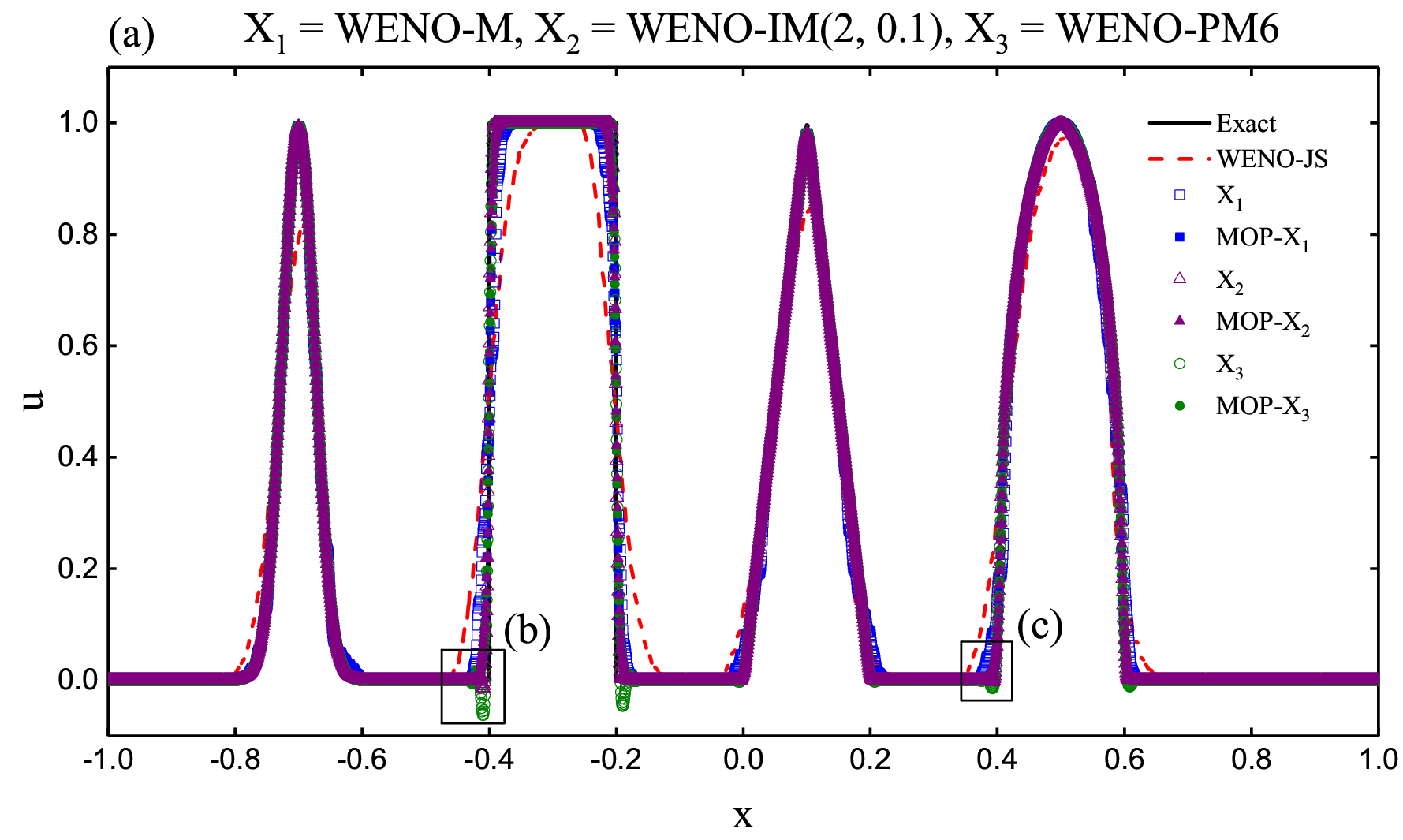}
\includegraphics[height=0.32\textwidth]
{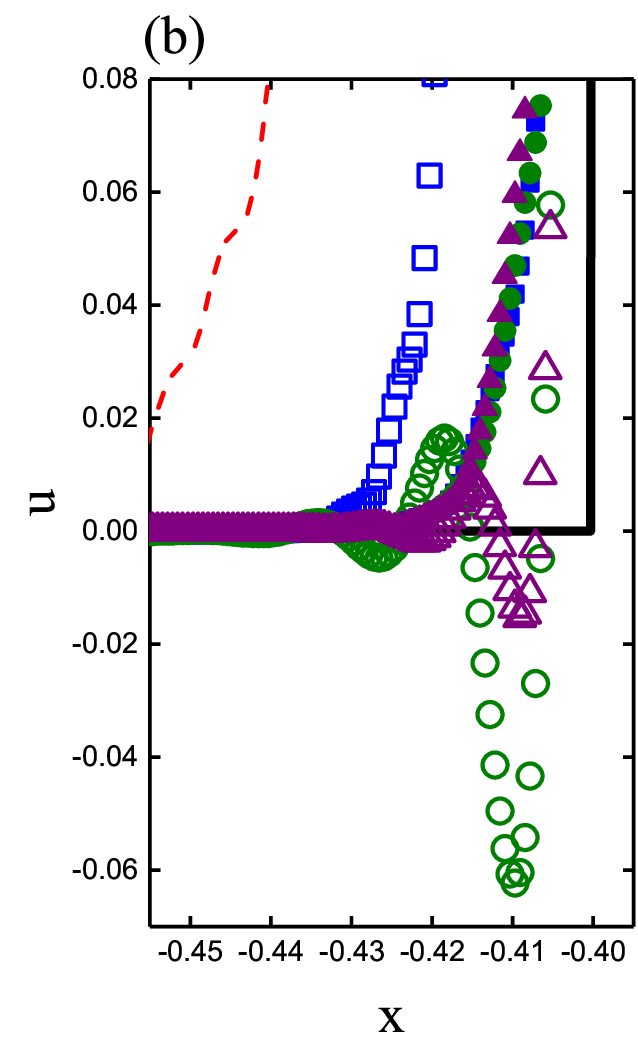}
\includegraphics[height=0.32\textwidth]
{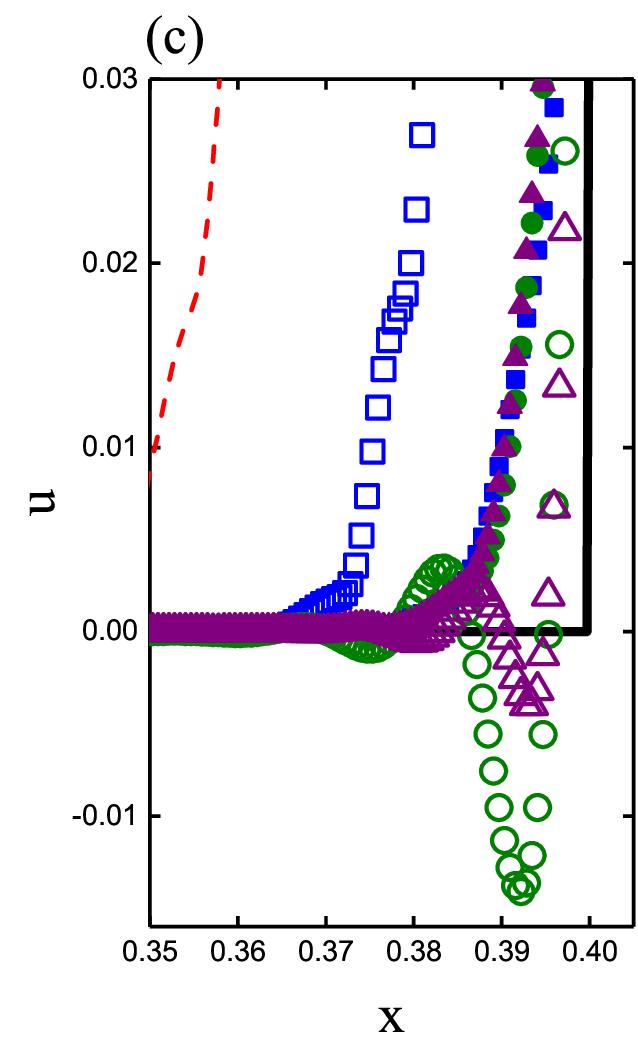}
\caption{Performance of the WENO-JS, WENO-M, MOP-WENO-M, 
WENO-IM($2,0.1$), MOP-WENO-IM($2,0.1$), WENO-PM6 and MOP-WENO-PM6 
schemes for the SLP at output time $t = 200$ with a uniform mesh 
size of $N = 3200$.}
\label{fig:SLP:N3200:PM6}
\end{figure}

\begin{figure}[ht]
\centering
\includegraphics[height=0.32\textwidth]
{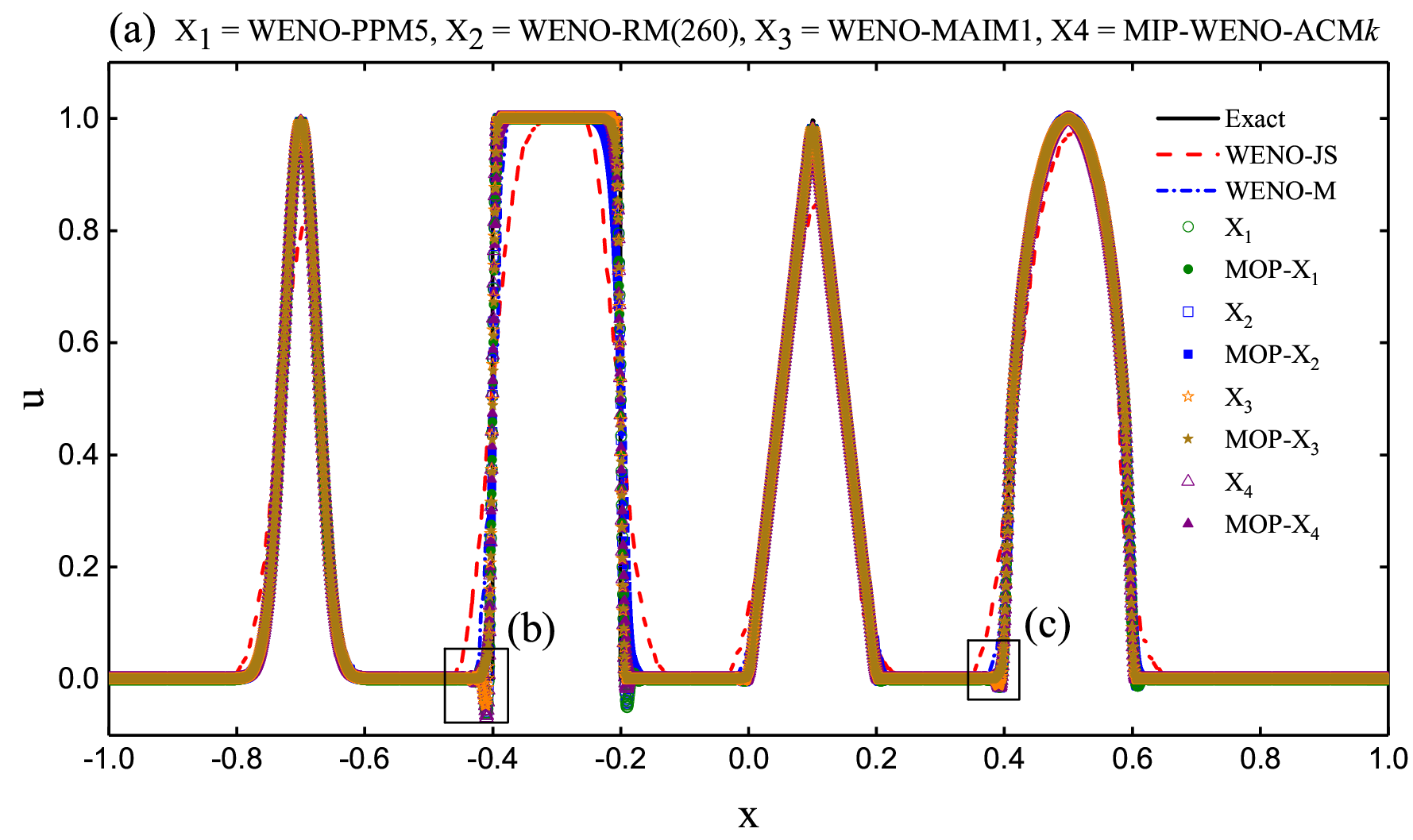}
\includegraphics[height=0.32\textwidth]
{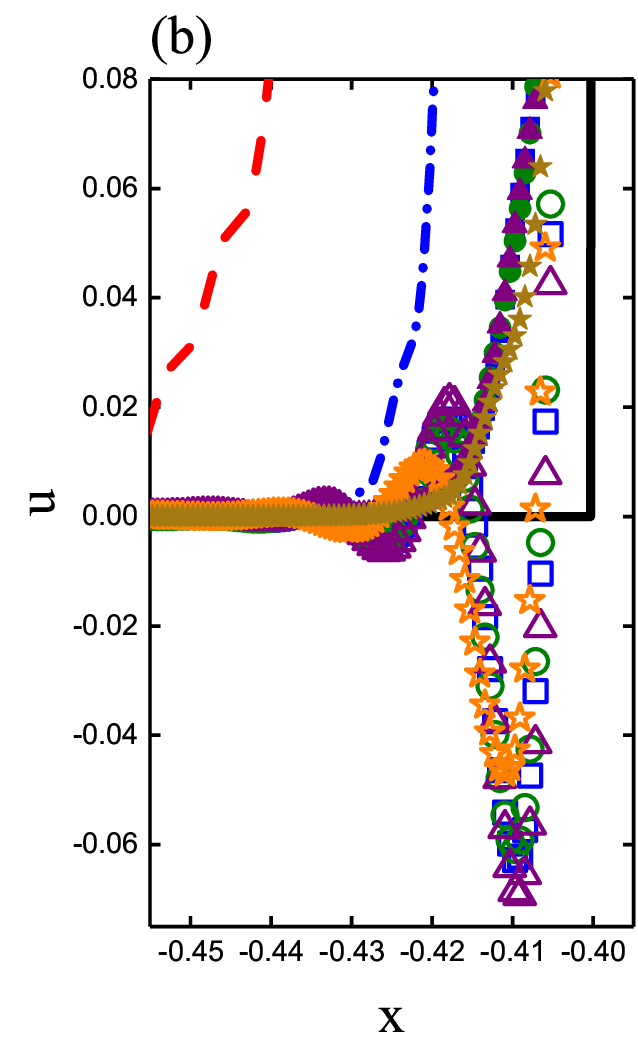}
\includegraphics[height=0.32\textwidth]
{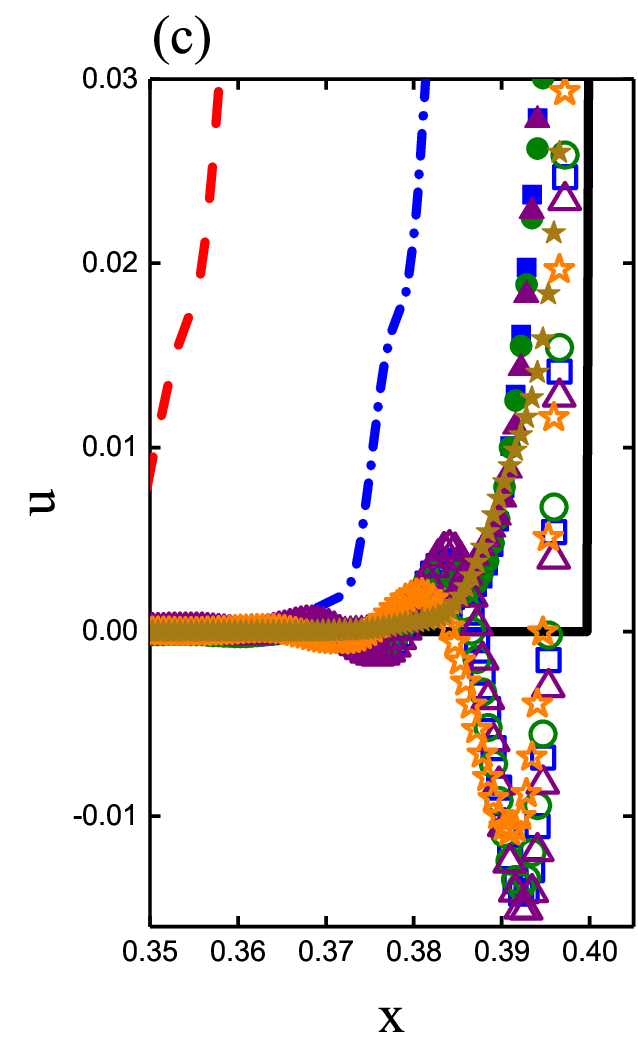}
\caption{Performance of the WENO-JS, WENO-M, WENO-PPM5, 
MOP-WENO-PPM5, WENO-RM260, MOP-WENO-RM260, WENO-MAIM1,
MOP-WNEO-MAIM1, MIP-WENO-ACM$k$ and MOP-WENO-ACM$k$ schemes for the 
SLP at output time $t = 200$ with a uniform mesh size of $N = 3200$.}
\label{fig:SLP:N3200:RM260}
\end{figure}

\begin{figure}[ht]
\centering
\includegraphics[height=0.32\textwidth]
{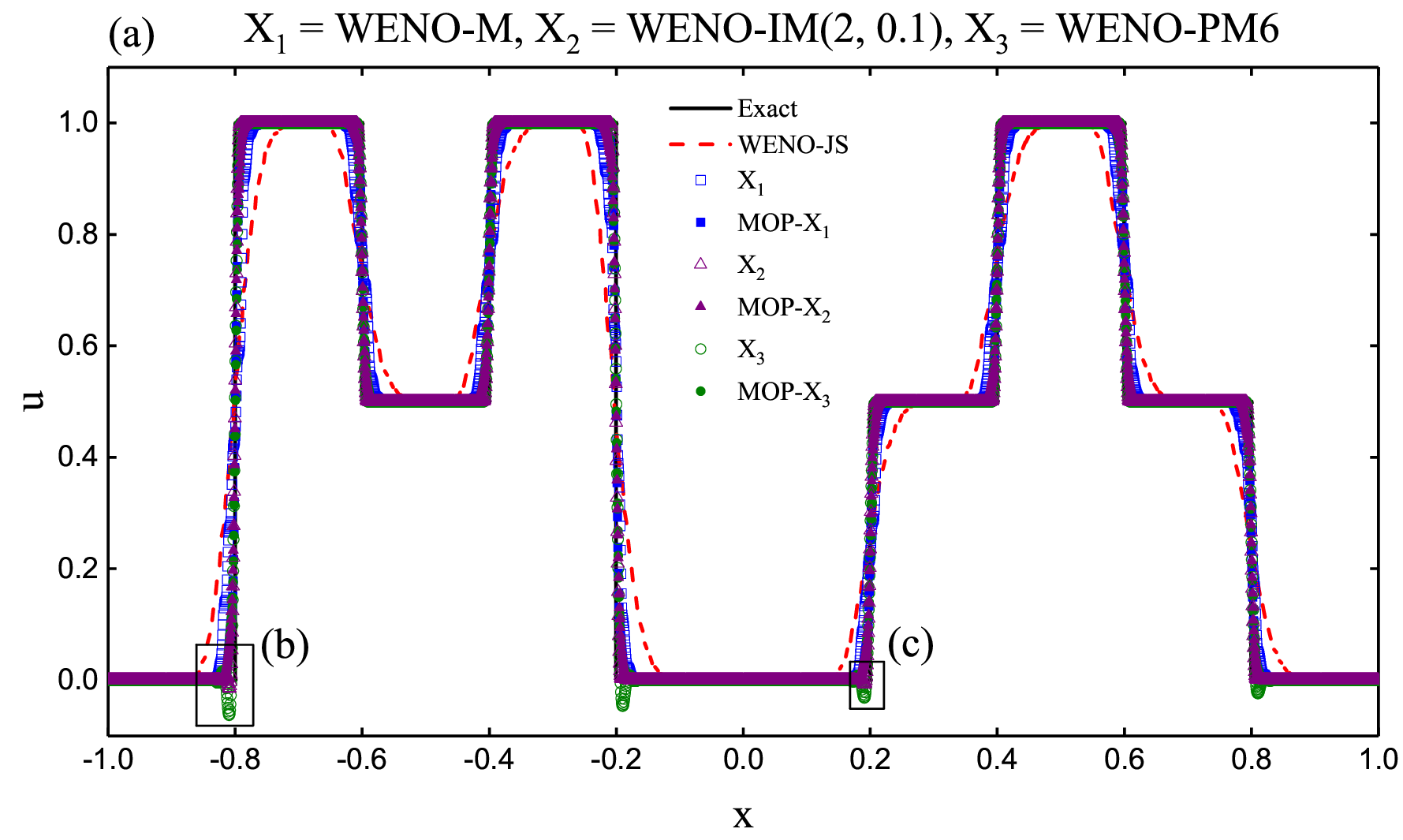}
\includegraphics[height=0.32\textwidth]
{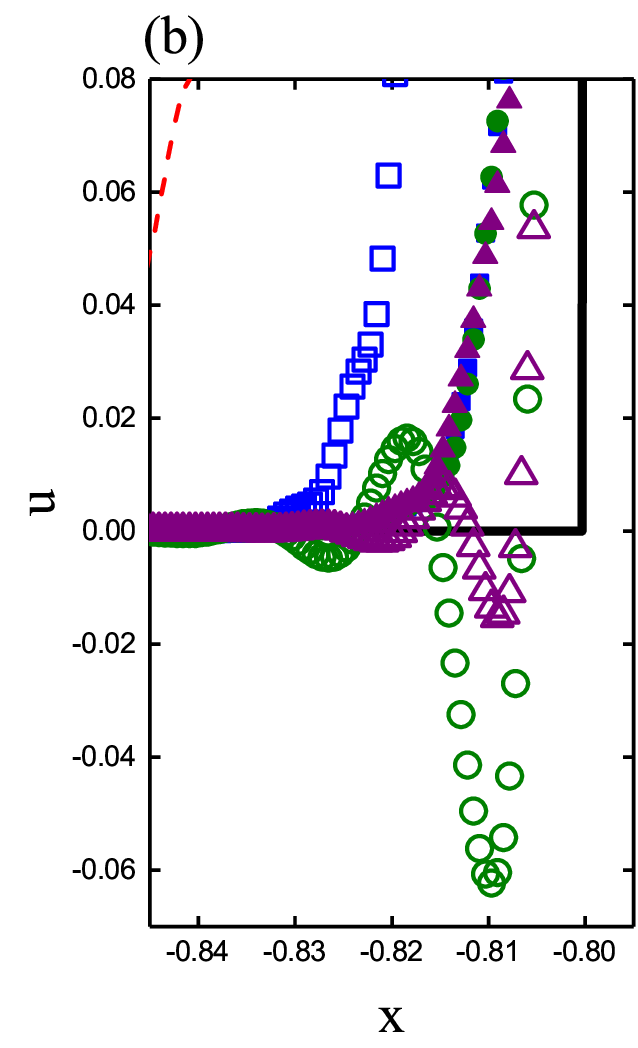}
\includegraphics[height=0.32\textwidth]
{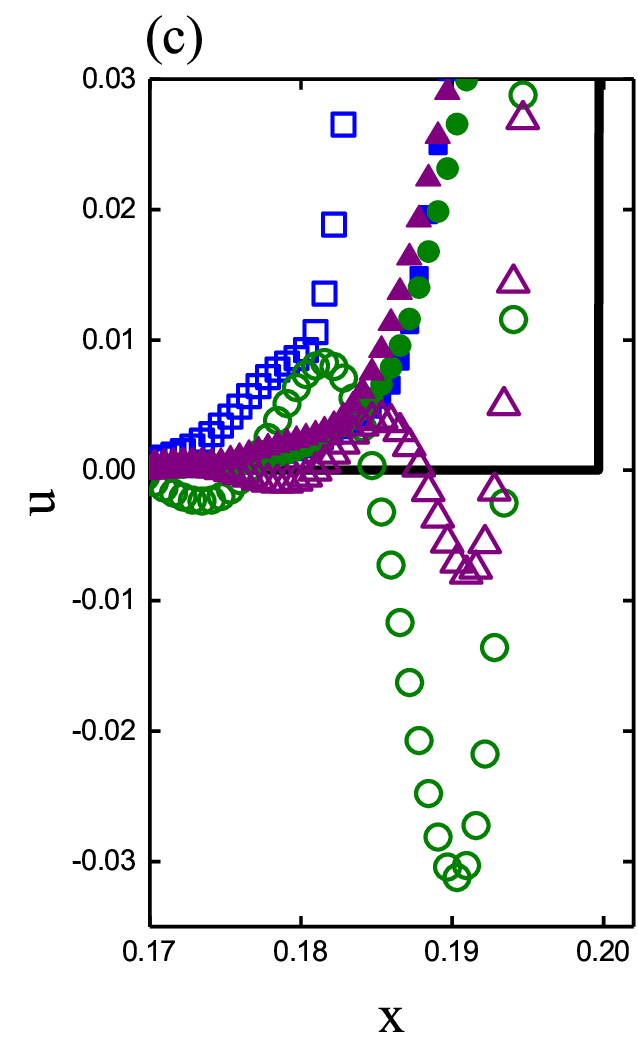}
\caption{Performance of the WENO-JS, WENO-M, MOP-WENO-M, 
WENO-IM($2,0.1$), MOP-WENO-IM($2,0.1$), WENO-PM6 and MOP-WENO-PM6 
schemes for the BiCWP at output time $t = 200$ with a uniform mesh 
size of $N = 3200$.}
\label{fig:BiCWP:N3200:PM6}
\end{figure}

\begin{figure}[ht]
\centering
\includegraphics[height=0.32\textwidth]
{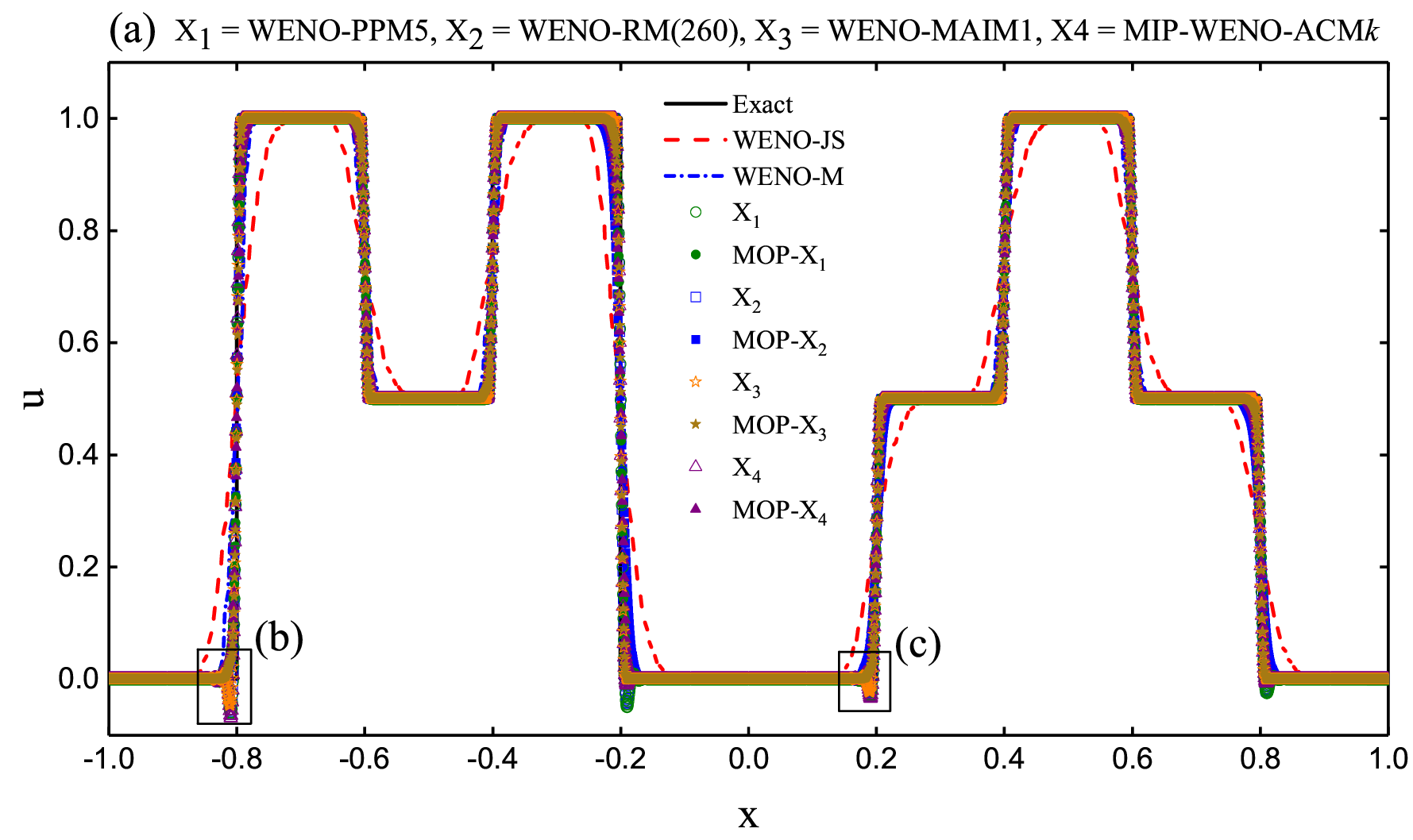}
\includegraphics[height=0.32\textwidth]
{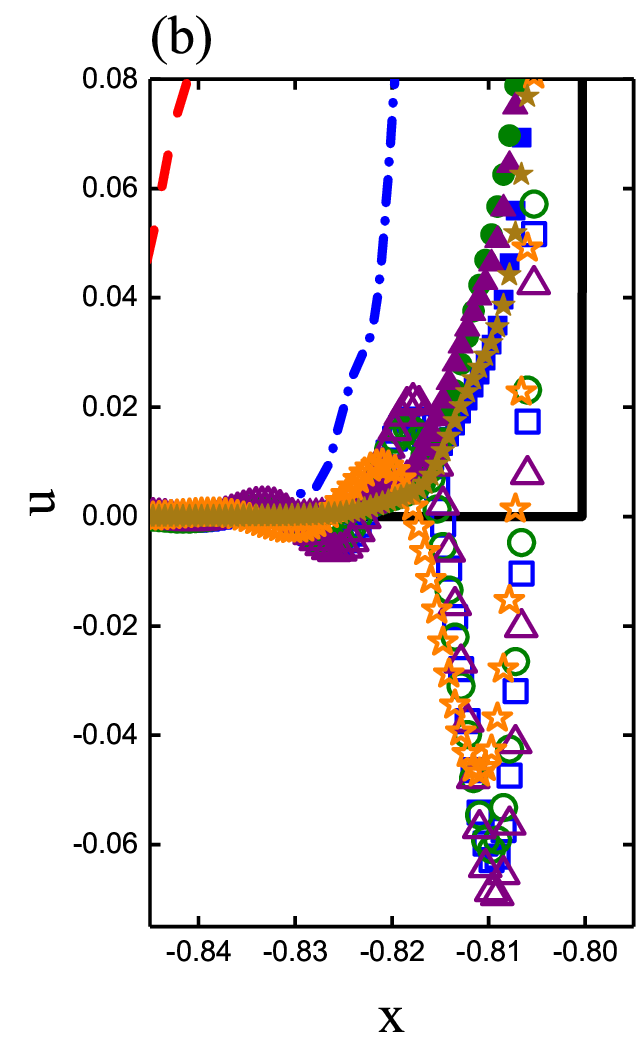}
\includegraphics[height=0.32\textwidth]
{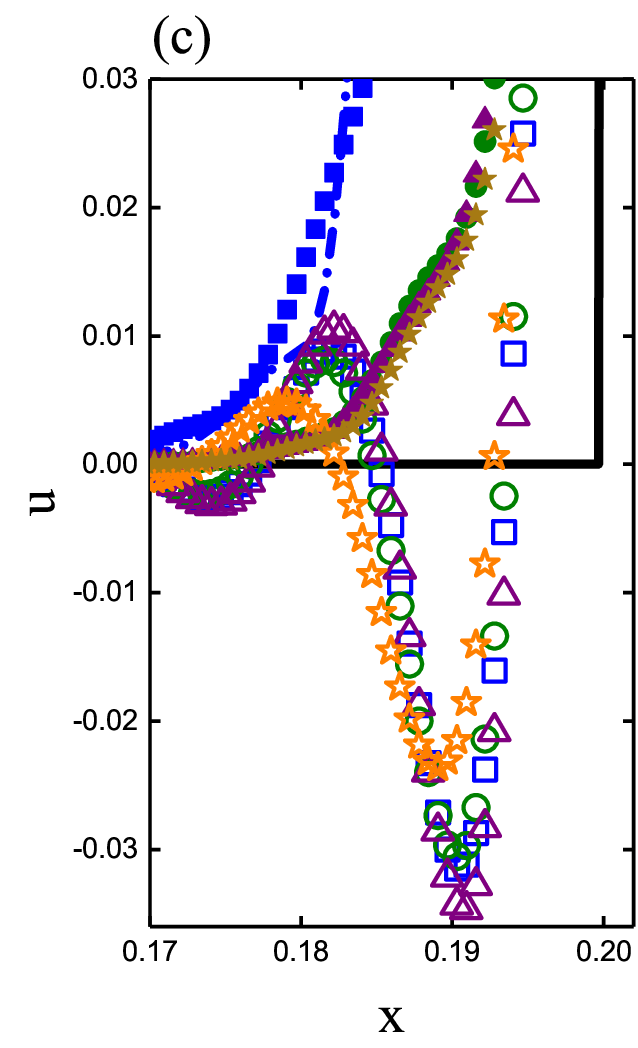}
\caption{Performance of the WENO-JS, WENO-M, WENO-PPM5, 
MOP-WENO-PPM5, WENO-RM260, MOP-WENO-RM260, WENO-MAIM1,
MOP-WNEO-MAIM1, MIP-WENO-ACM$k$ and MOP-WENO-ACM$k$ schemes for the 
BiCWP at output time $t = 200$ with a uniform mesh size of $N=3200$.}
\label{fig:BiCWP:N3200:RM260}
\end{figure}

In Figs. \ref{fig:SLP:N3200:PM6}, \ref{fig:SLP:N3200:RM260} and 
Figs. \ref{fig:BiCWP:N3200:PM6}, \ref{fig:BiCWP:N3200:RM260}, we 
show the comparison of considered schemes for SLP and BiCWP 
respectively, by taking $t = 200$ and $N = 3200$. From these 
solutions computed with larger grid numbers and a reduced but still 
long output time, it can be seen that: (1) firstly, the WENO-IM(2,
0.1) scheme generates spurious oscillations but the MOP-WENO-IM(2, 
0.1) scheme does not while provides an improved resolution when 
solving SLP; (2) although the resolutions of the results computed by 
the WENO-JS and WENO-M schemes are significantly improved for both 
SLP and BiCWP, the MOP-WENO-X schemes still evidently provide much 
better resolutions; (3) the spurious oscillations generated by the 
WENO-X schemes appear to be more evident and more intense as the 
grid number increases, while the corresponding MOP-WENO-X schemes 
can still avoid spurious oscillations but obtain higher resolutions, 
when solving both SLP and BiCWP.

For the further interpretation, without loss of generality, in Fig. 
\ref{fig:x-Omega:SLP_BiCWP}, we present the \textit{non-OP points} 
of the numerical solutions of SLP computed by the WENO-M and 
MOP-WENO-M schemes with $N=800, t=2000$, and the 
\textit{non-OP points} of the numerical solutions of BiCWP computed 
by the WENO-PM6 and MOP-WENO-PM6 schemes with $N=3200, t=200$. 
We can find that there are a great many \textit{non-OP points} in 
the solutions of the WENO-M and WENO-PM6 schemes while the numbers 
of the \textit{non-OP points} in the solutions of the MOP-WENO-M and 
MOP-WENO-PM6 schemes are zero. Actually, there are many 
\textit{non-OP points} for all considered mapped WENO-X schemes. And 
as expected, there are no \textit{non-OP points} for the 
corresponding MOP-WENO-X schemes and the WENO-JS scheme for all 
computing cases here. We do not show the results of the 
\textit{non-OP points} for all computing cases here just for the 
simplicity of illustration. 

In summary, it could be indicated that the general method to 
introduce the \textit{OP} mapping can help to gain the advantage of 
achieving high resolutions and in the meantime preventing spurious 
oscillations when solving problems with discontinuities for long 
output times. And this is the most important point we want to 
report in this paper.

\subsection{Euler system in two dimension}
\label{subsec:examples_2D_Euler}
In this subsection, we focus on the numerical simulations of the 
shock-vortex interaction problem \cite{Shock-vortex_interaction-1,
Shock-vortex_interaction-2,Shock-vortex_interaction-3} and the 2D 
Riemann problem \cite{Riemann-2D-01,Riemann2D-02,Riemann2D-03}. They 
are governed by the two-dimensional Euler system of gas dynamics, 
taking the following strong conservation form of mass, momentum and 
energy
\begin{equation}
\begin{array}{ll}
\begin{aligned}
&\dfrac{\partial \rho}{\partial t} + \dfrac{\partial (\rho u)}
{\partial x} + \dfrac{\partial (\rho v)}{\partial y} = 0, \\
&\dfrac{\partial (\rho u)}{\partial t} + \dfrac{\partial (\rho u^{2}
+ p)}{\partial x} + \dfrac{\partial (\rho uv)}{\partial y}= 0, \\
&\dfrac{\partial (\rho v)}{\partial t} + \dfrac{\partial (\rho vu)}
{\partial x} + \dfrac{\partial (\rho v^{2} + p)}{\partial y}= 0, \\
&\dfrac{\partial E}{\partial t} + \dfrac{\partial (uE + up)}
{\partial x} + \dfrac{\partial (vE + vp)}{\partial y} = 0, \\
\end{aligned}
\end{array}
\label{2DEulerEquations}
\end{equation}
where $\rho, u, v, p$ and $E$ are the density, components of 
velocity in the $x$ and $y$ coordinate directions, pressure and 
total energy, respectively. The following equation of state for an 
ideal polytropic gas is used to close the two-dimensional Euler 
system Eq.(\ref{2DEulerEquations})
\begin{equation*}
p = (\gamma - 1)\Big( E - \dfrac{1}{2}\rho(u^{2} + v^{2})\Big),
\end{equation*}
where $\gamma$ is the ratio of specific heat, and we set $\gamma=1.4$
in this paper. In the computations below, the CFL number is taken to 
be $0.5$. All the considered WENO schemes are applied 
dimension-by-dimension to solve the two-dimensional Euler system 
and the local characteristic decomposition \cite{WENO-JS} is used. 
In \cite{FVMaccuracyProofs03}, Zhang et al. investigated two 
commonly used classes of finite volume WENO schemes in 
two-dimensional Cartesian meshes, and we employ the one denoted as 
class A in this subsection.

\begin{example}
\bf{(Shock-vortex interaction)}
\rm{We consider the shock-vortex interaction problem used in 
\cite{Shock-vortex_interaction-1,Shock-vortex_interaction-2,
Shock-vortex_interaction-3}. It consists of the interaction of a 
left moving shock wave with a right moving vortex. The computational 
domain is initialized by}
\label{ex:shock-vortex}
\end{example}
\begin{equation*}
\big( \rho, u, v, p \big)(x, y, 0) = \left\{
\begin{aligned}
\begin{array}{ll}
\textbf{U}_{\mathrm{L}}, & x < 0.5, \\
\textbf{U}_{\mathrm{R}}, & x \geq 0.5, \\
\end{array}
\end{aligned}
\right.
\label{eq:initial_Euer2D:shock-vortex-interaction}
\end{equation*}
where $\textbf{U}_{\mathrm{L}} = (\rho_{\mathrm{L}}, u_{\mathrm{L}}, 
v_{\mathrm{L}}, p_{\mathrm{L}}) = (1, \sqrt{\gamma}, 0, 1)$, and 
$\textbf{U}_{\mathrm{R}} = (\rho_{\mathrm{R}}, u_{\mathrm{R}}, 
v_{\mathrm{R}}, p_{\mathrm{R}})$ taking the form
\begin{equation*}
\begin{array}{l}
p_{\mathrm{R}} = 1.3, \rho_{\mathrm{R}} = \rho_{\mathrm{L}}\bigg( 
\dfrac{\gamma - 1 + (\gamma + 1)p_{\mathrm{R}}}{\gamma + 1 + (\gamma
- 1)p_{\mathrm{R}}} \bigg)\\
u_{\mathrm{R}} = u_{\mathrm{L}}\bigg( \dfrac{1 - p_{\mathrm{R}}}{
\sqrt{\gamma - 1 + p_{\mathrm{R}}(\gamma + 1)}}\bigg), v_{\mathrm{R}}
= 0.
\end{array}
\label{eq:Euler2D:shock-vortex-interaction:rightState}
\end{equation*}
The vortex $\delta \textbf{U} = (\delta \rho, \delta u, \delta v, 
\delta p)$, defined by the following perturbations, is superimposed 
onto the left state $\textbf{U}_{\mathrm{L}}$,
\begin{equation*}
\delta \rho = \dfrac{\rho_{\mathrm{L}}^{2}}{(\gamma - 1)
p_{\mathrm{L}}}\delta T, 
\delta u = \epsilon \dfrac{y - y_{\mathrm{c}}}{r_\mathrm{c}}
\mathrm{e}^{\alpha(1-r^{2})}, 
\delta v = - \epsilon \dfrac{x - x_{\mathrm{c}}}{r_\mathrm{c}}
\mathrm{e}^{\alpha(1-r^{2})}, 
\delta p = \dfrac{\gamma \rho_{\mathrm{L}}^{2}}{(\gamma - 1)
\rho_{\mathrm{L}}}\delta T,
\label{eq:Euler2D:shock-vortex-interactions:Perturbations}
\end{equation*}
where $\epsilon = 0.3, r_{\mathrm{c}} = 0.05, \alpha = 0.204, 
x_{\mathrm{c}} = 0.25, y_{\mathrm{c}} = 0.5,
r = \sqrt{((x - x_{\mathrm{c}})^{2} + (y - y_{\mathrm{c}})^{2})/r_{
\mathrm{c}}^{2}}, \delta T = - (\gamma - 1)\epsilon^{2}\mathrm{e}^{2
\alpha (1 - r^{2})}/(4\alpha \gamma)$.
The transmissive boundary condition is used on all boundaries. A 
uniform mesh size of $800 \times 800$ is used and the output time is 
set to be $t = 0.35$.

We calculate this problem using all the considered mapped WENO-X 
schemes in Table \ref{table:GF_parameters} and their corresponding 
MOP-WENO-X schemes. For the sake of brevity though, we only present
the solutions of the WENO-M, WENO-IM(2, 0.1), WENO-PPM5, WENO-MAIM1 
schemes and their corresponding MOP-WENO-X schemes in 
Figs. \ref{fig:ex:SVI:1} and \ref{fig:ex:SVI:2}, where the first 
rows give the final structures of the shock and vortex in density 
profile of the existing mapped WENO-X schemes, the second rows give 
those of the corresponding MOP-WENO-X schemes, and the third rows 
give the cross-sectional slices of density plot along the plane 
$y = 0.65$ where $x \in [0.70, 0.76]$. We find that all the 
considered schemes perform well in capturing the main structure of 
the shock and vortex after the interaction. It can be seen that 
there are clear post-shock oscillations in the solutions of the 
WENO-M, WENO-IM(2, 0.1) and WENO-PPM5 schemes. However,
in the solutions of the MOP-WENO-M, MOP-WENO-IM(2, 0.1) and 
MOP-WENO-PPM5 schemes, the post-shock oscillations are either 
removed or significantly reduced. The post-shock oscillations of the 
WENO-MAIM1 scheme are very slight and even hard to be noticed. 
Actually, it seems difficult to distinguish the solutions of the 
WENO-MAIM1 scheme from that of the MOP-WENO-MAIM1 scheme only 
according to the structure of the shock and vortex in the density 
profile. Nevertheless, when taking a closer look from the 
cross-sectional slices of the density profile along the plane 
$y=0.65$ at the bottom right picture of Fig. \ref{fig:ex:SVI:2} 
where the reference solution is obtained using the WENO-JS scheme 
with a uniform mesh size of $1600 \times 1600$, we can see that the 
post-shock oscillation of the WENO-MAIM1 scheme is very remarkable 
while it is imperceptible for the MOP-WENO-MAIM1 scheme. Also, from 
the third rows of Figs. \ref{fig:ex:SVI:1} and \ref{fig:ex:SVI:2}, 
we find that the WENO-IM(2, 0.1) and WENO-PPM5 schemes generate the 
post-shock oscillations with much bigger amplitudes than that of the 
WENO-MAIM1 scheme. The WENO-M scheme also generates clear post-shock 
oscillations with the amplitudes slightly smaller than that of the 
WENO-IM(2, 0.1) and WENO-PPM5 schemes. Evidently, the solutions of 
the MOP-WENO-M, MOP-WENO-IM(2, 0.1) and MOP-WENO-PPM5 schemes almost 
generate no post-shock oscillations or only generate some 
imperceptible numerical oscillations and their solutions are very 
close to the reference solution, and this should be an advantage of 
the mapped WENO schemes whose mapping functions are \textit{OP}.

\begin{example}
\bf{(2D Riemann problem)} 
\rm{It is very favorable to test the high-resolution numerical 
methods \cite{Riemann2D-03,WENO-PPM5,Riemann2D-04} using the series 
of 2D Riemann problems \cite{Riemann-2D-01,Riemann2D-02}. In 
\cite{Riemann2D-03}, Lax et al. classified a total of 19 genuinely 
different Configurations for 2D Riemann problem and calculated all 
the numerical solutions. Configuration 4 is chosen here for the 
test, and the computational domain is initialized by}
\label{ex:Riemann2D}
\end{example}
\begin{equation*}
\big( \rho, u, v, p \big)(x, y, 0) = \left\{
\begin{aligned}
\begin{array}{ll}
(1.1, 0.0, 0.0, 1.1),         & 0.5 \leq x \leq 1.0, 
                                0.5 \leq y \leq 1.0, \\
(0.5065, 0.8939, 0.0, 0.35),  & 0.0 \leq x \leq 0.5, 
                                0.5 \leq y \leq 1.0, \\
(1.1, 0.8939, 0.8939, 1.1),   & 0.0 \leq x \leq 0.5, 
                                0.0 \leq y \leq 0.5, \\
(0.5065, 0.0, 0.8939, 0.35),  & 0.5 \leq x \leq 1.0, 
                                0.0 \leq y \leq 0.5. \\
\end{array}
\end{aligned}
\right.
\label{eq:initial_Euer2D:Riemann2D}
\end{equation*}
The transmission boundary condition is used on all boundaries, and 
the numerical solutions are calculated on a uniform mesh size of 
$800 \times 800$. The computations proceed to $t = 0.25$.

Similarly, although we calculate this problem using all the 
considered mapped WENO-X schemes in Table \ref{table:GF_parameters} 
and their corresponding MOP-WENO-X schemes, we only present the 
solutions of the WENO-M, WENO-PM6, WENO-RM260 and MIP-WENO-ACM$k$ 
schemes and their corresponding MOP-WENO-X schemes here for the sake 
of brevity. We have shown the numerical results of density obtained 
by using these schemes in Figs. \ref{fig:ex:Riemann2D:1} and 
\ref{fig:ex:Riemann2D:2}, where the first rows give the structures 
of the 2D Riemann problem in density profile of the existing mapped 
WENO-X schemes, the second rows give those of the corresponding 
MOP-WENO-X schemes, and the third rows give the cross-sectional 
slices of density plot along the plane $y = 0.5$ where $x \in 
[0.65, 0.692]$. We can see that all schemes can capture the main 
structure of the solution. However, we can also observe that there 
are obvious post-shock oscillations (as marked by the pink boxes), 
which are unfavorable for the fidelity of the results, in the 
solutions of the WENO-M, WENO-PM6, WENO-RM(260) and MIP-WENO-ACM$k$ 
schemes. These post-shock oscillations can be seen more clearly from 
the cross-sectional slices of density profile as presented in the 
third rows of Figs. \ref{fig:ex:Riemann2D:1} and 
\ref{fig:ex:Riemann2D:2}, where the reference solution is obtained 
by using the WENO-JS scheme with a uniform mesh size of 
$3000 \times 3000$. Noticeably, there are either almost no 
or imperceptible post-shock oscillations in the solutions of the 
MOP-WENO-M, MOP-WENO-PM6, MOP-WENO-RM(RM260) and MOP-WENO-ACM$k$ 
schemes. Again, we believe that this should be an advantage of the 
mapped WENO schemes whose mapping functions are \textit{OP}.


%% file: article_conclusion.tex
\section{Conclusions}
\label{secConclusions} 
The concept of \textit{OP-Mapped WENO} schemes standing for the 
family of the mapped WENO schemes with \textit{order-preserving (OP)}
mappings, as well as a general way to build one group of this kind 
of schemes, has been proposed in this paper. Specifically, we extend 
the \textit{OP} mapping introduced in \cite{MOP-WENO-ACMk} to 
various existing mapped WENO schemes in references by providing a 
general formula of their mapping functions. A systematic analysis 
has been made to prove that the improved mapped WENO scheme based on 
the existing mapped WENO-X scheme, denoted as MOP-WENO-X, generates 
numerical solutions with the same convergence rates of accuracy in 
smooth regions as the corresponding WENO-X scheme. Furthermore, 
numerical experiments are run to show that the MOP-WENO-X schemes 
have the same advantage as the mapped WENO scheme proposed in 
\cite{MOP-WENO-ACMk} in calculating the one-dimensional linear 
advection problems including discontinuities with long output times. 
The mapping functions of the MOP-WENO-X schemes are \textit{OP} and 
hence able to attain high resolutions and avoid spurious 
oscillations meanwhile. Moreover, numerical results with the 2D 
Euler system problems were presented to show that the MOP-WENO-X 
schemes perform well in simulating the two-dimensional steady 
problems with strong shock waves to capture the main flow structures 
and remove or significantly reduce the post-shock oscillations.

\begin{figure}[ht]
\centering
\includegraphics[height=0.305\textwidth]
{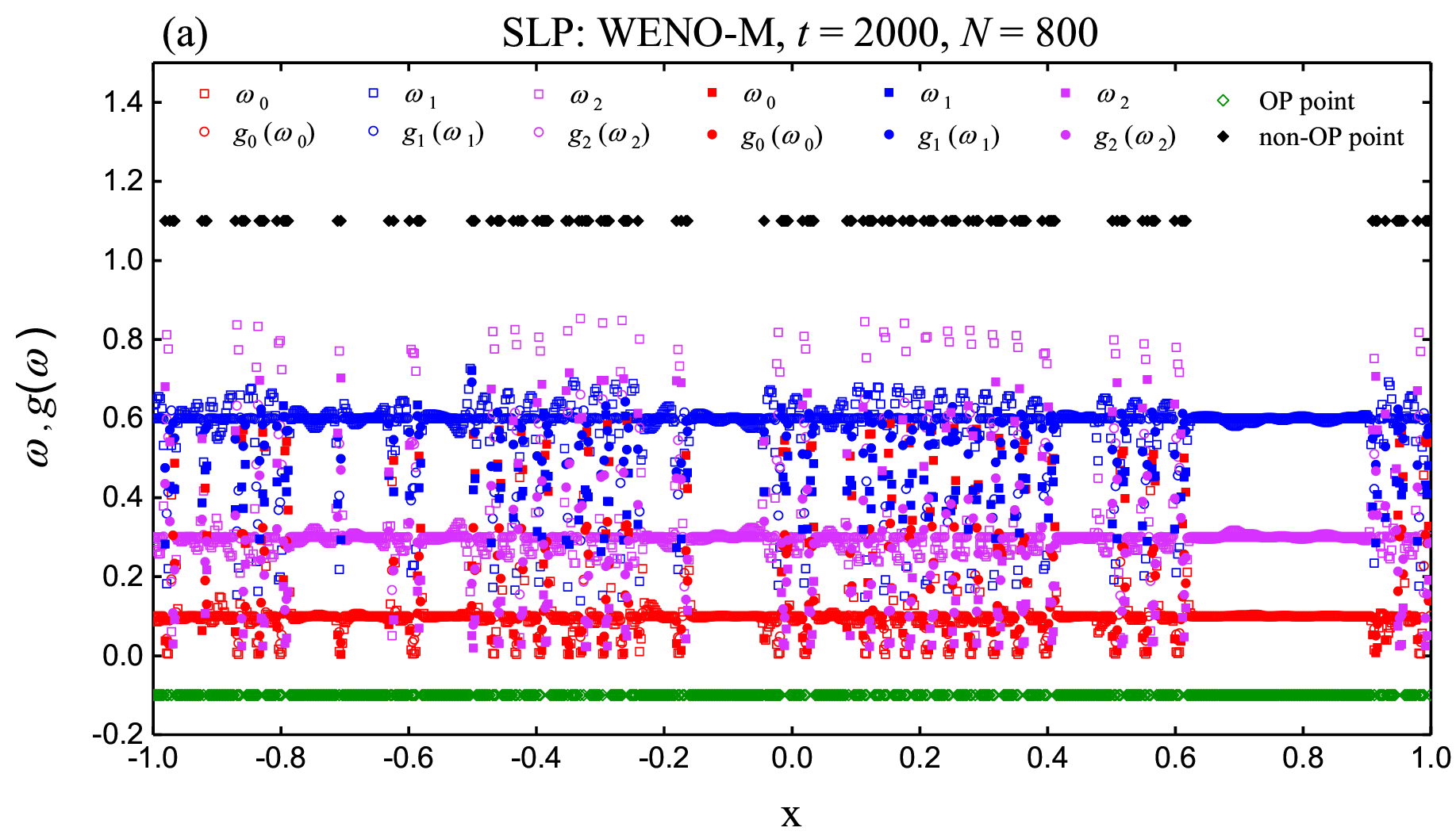}
\includegraphics[height=0.305\textwidth]
{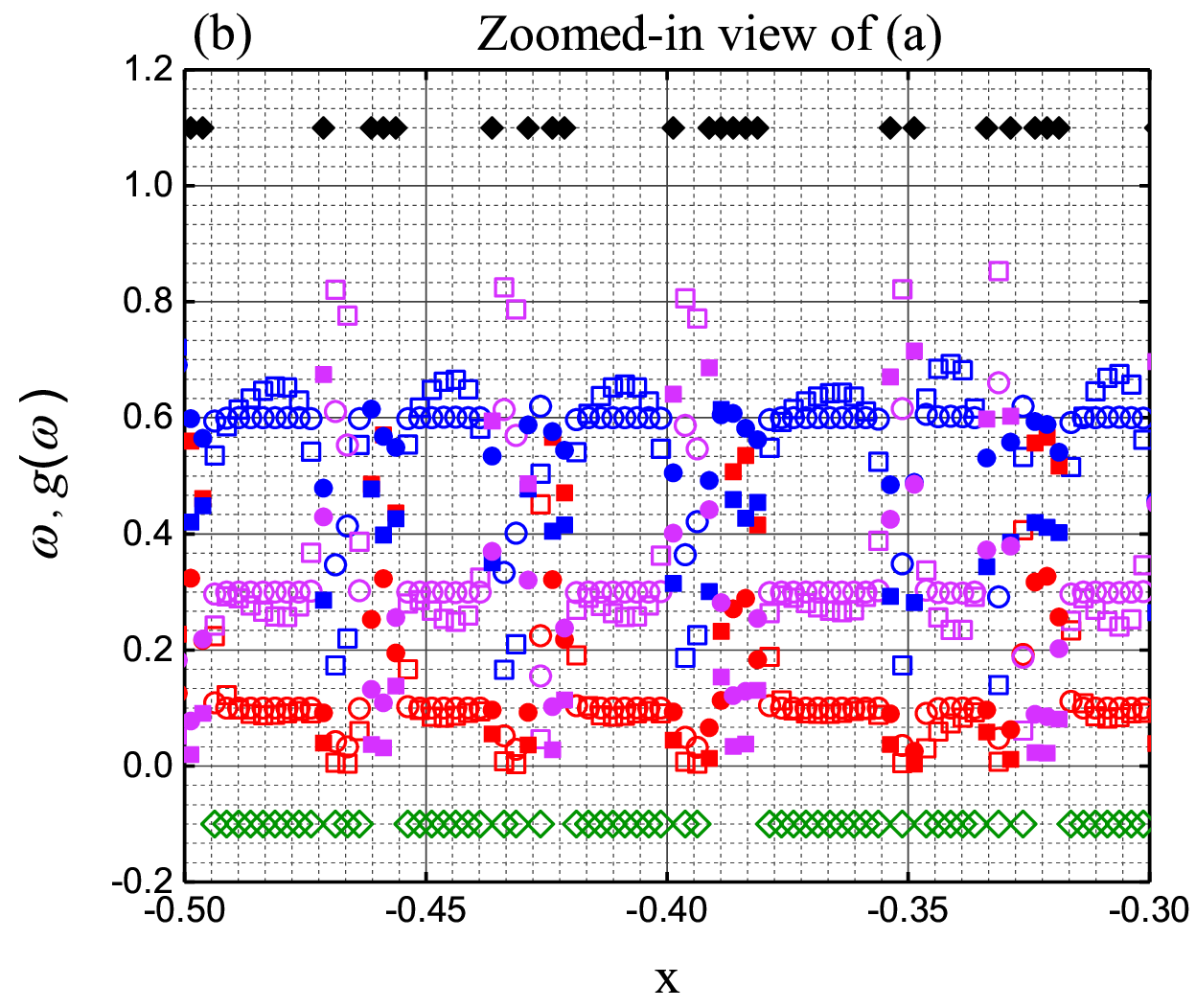}\\
\includegraphics[height=0.305\textwidth]
{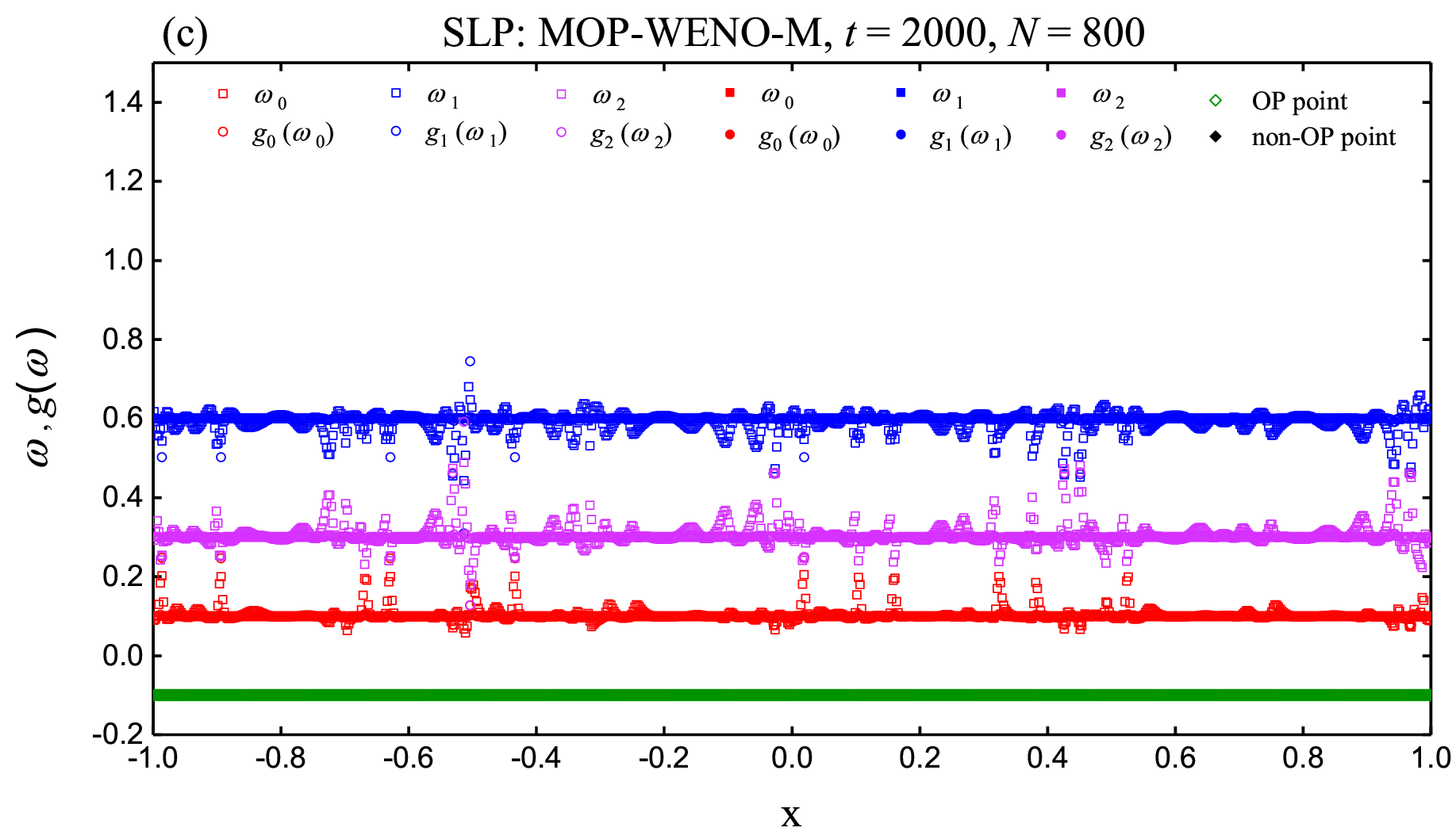}
\includegraphics[height=0.305\textwidth]
{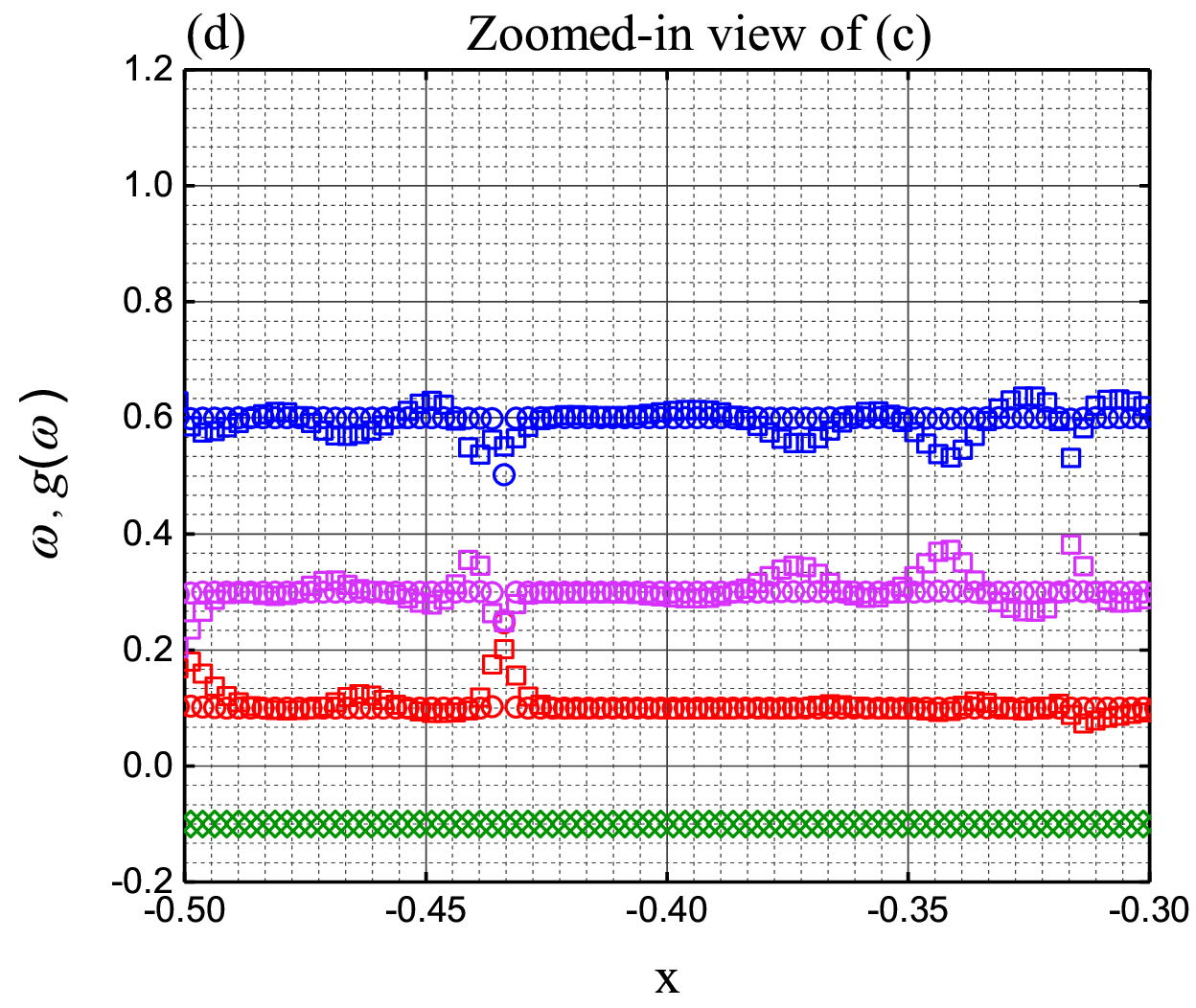}\\
\includegraphics[height=0.305\textwidth]
{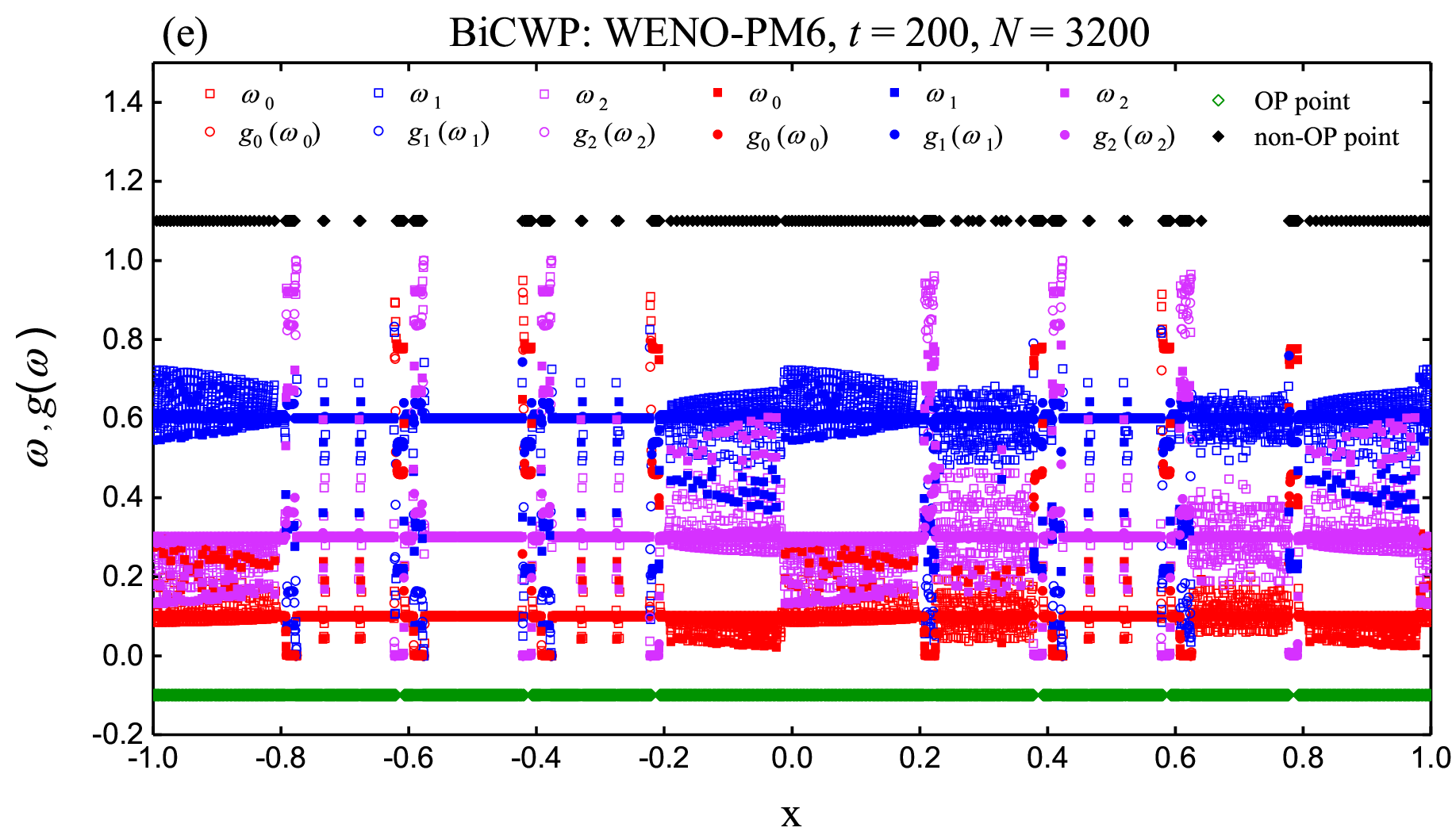}
\includegraphics[height=0.305\textwidth]
{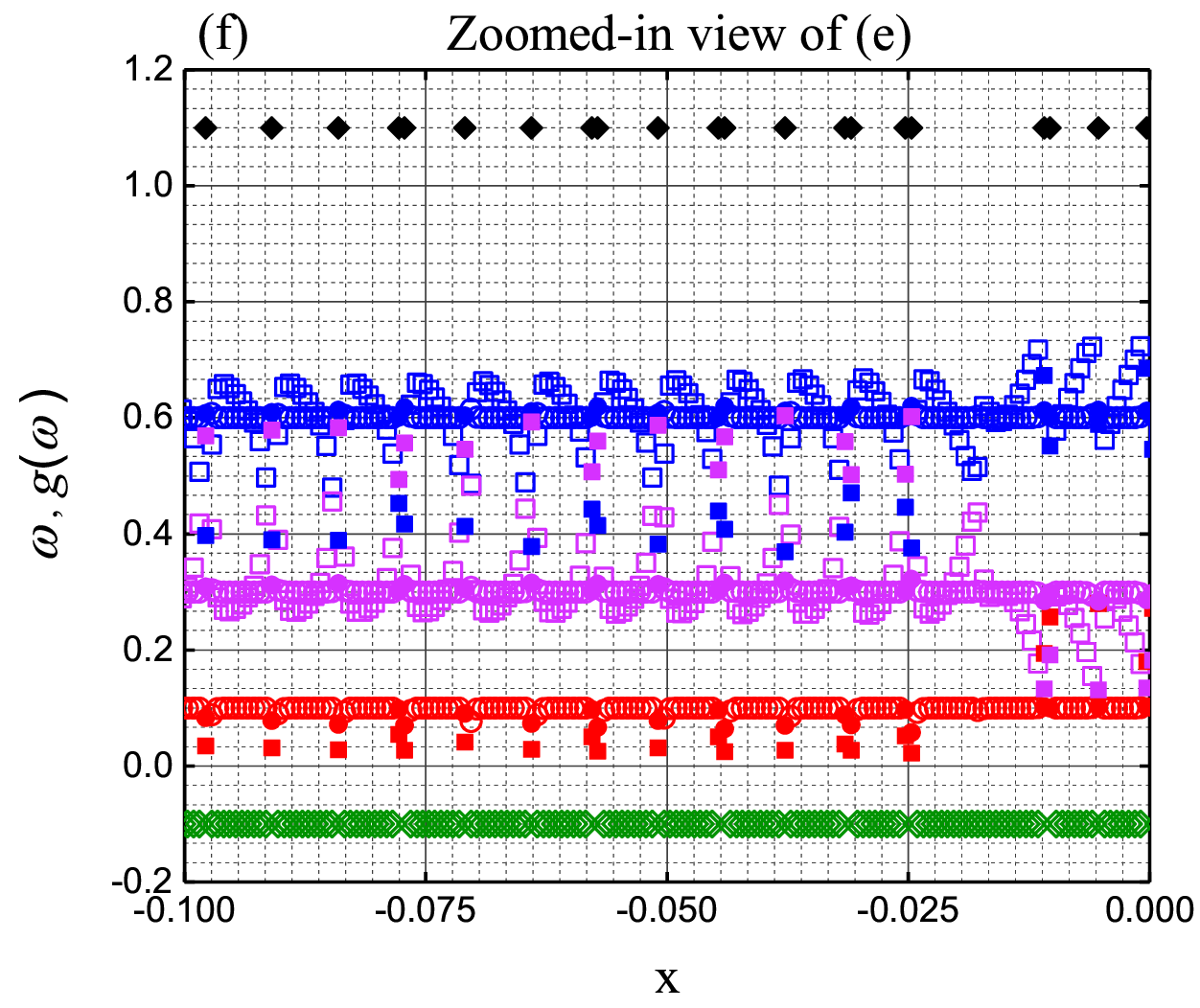}\\
\includegraphics[height=0.305\textwidth]
{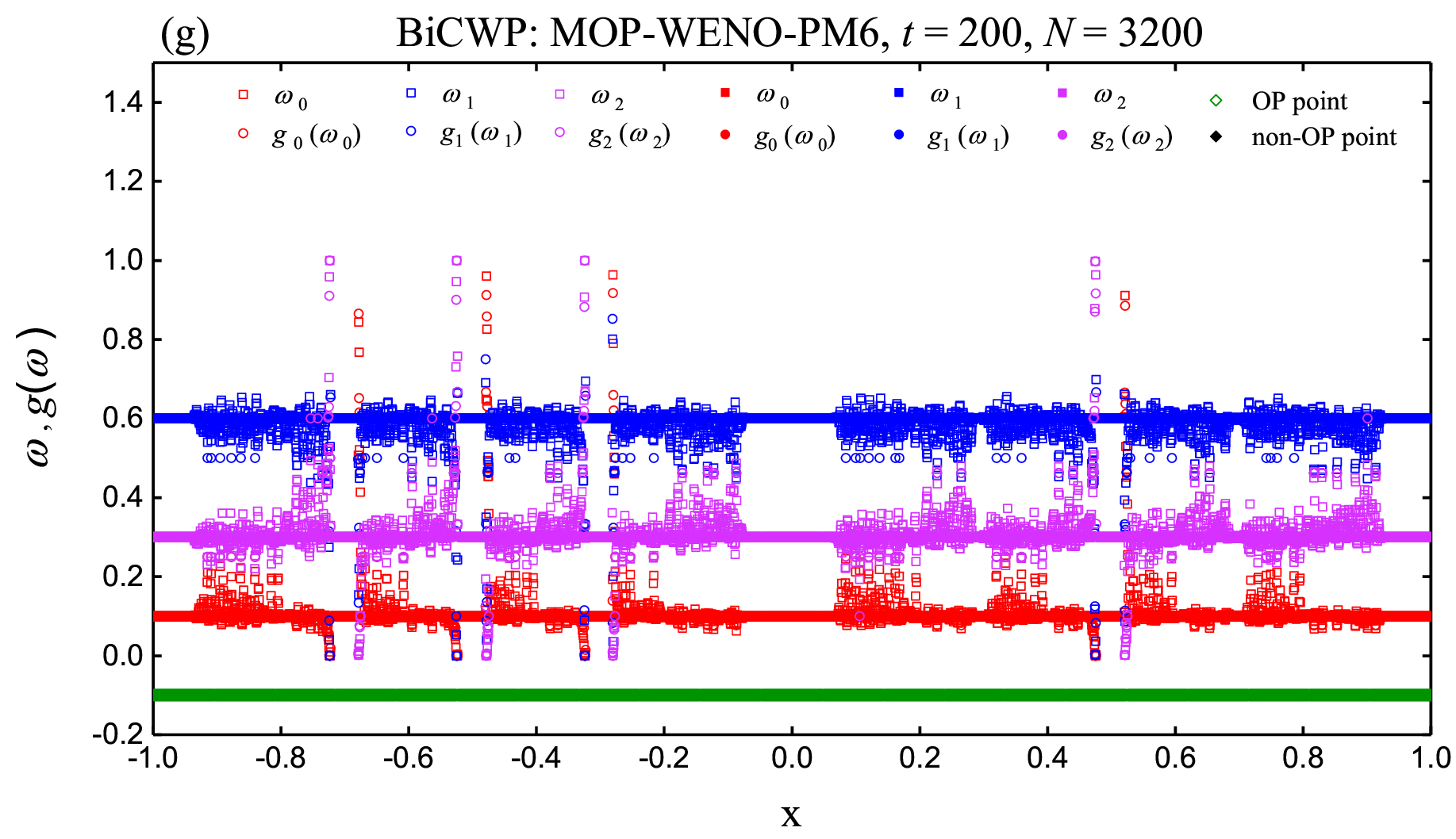}
\includegraphics[height=0.305\textwidth]
{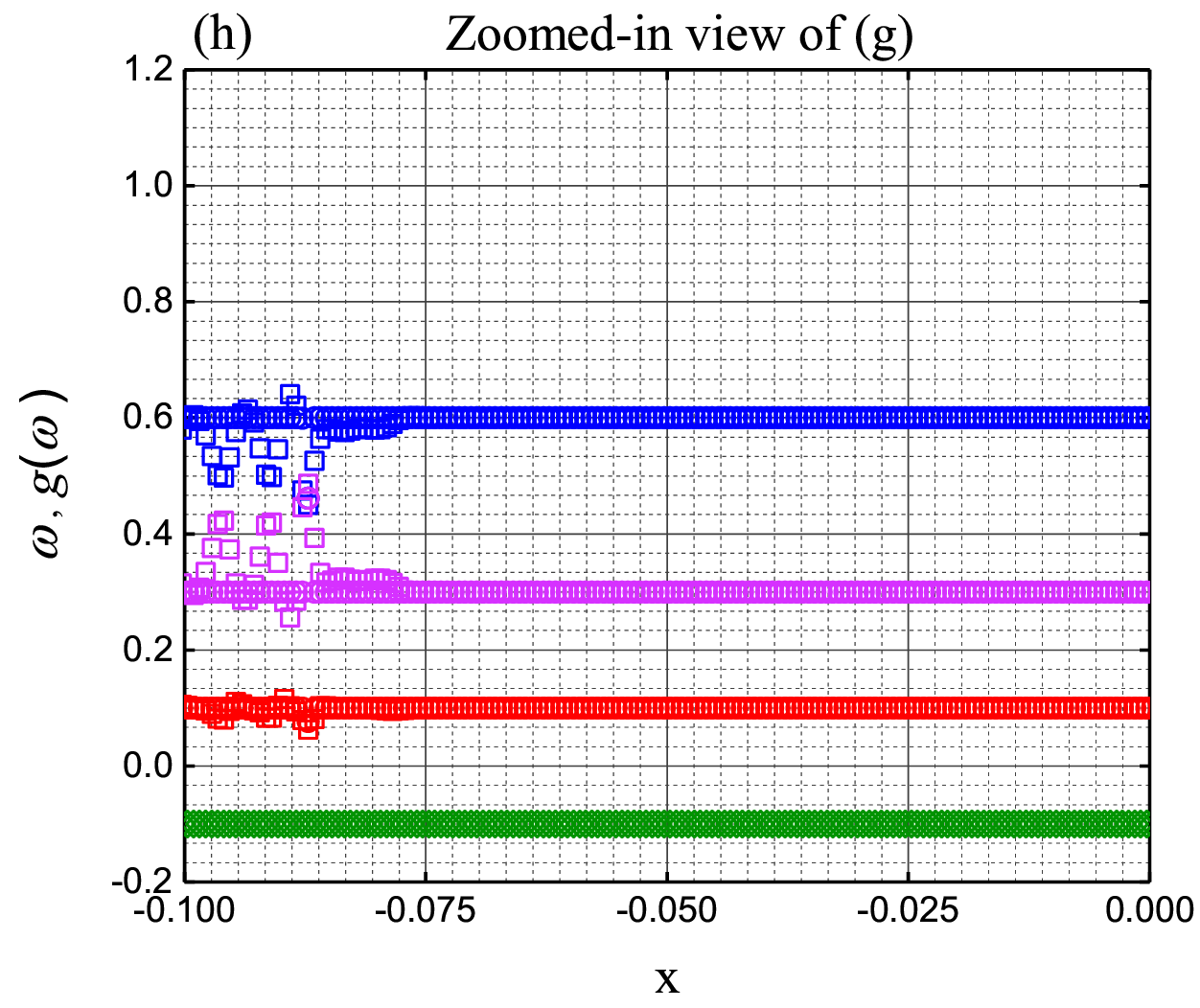}
\caption{The \textit{non-OP points} in the numerical solutions of 
SLP computed by the WENO-M and MOP-WENO-M schemes with $N=800,
t=2000$, and the \textit{non-OP points} in the numerical solutions 
of BiCWP computed by the WENO-PM6 and MOP-WENO-PM6 schemes with 
$N=3200, t= 200$.}
\label{fig:x-Omega:SLP_BiCWP}
\end{figure}

\begin{figure}[ht]
\centering
  \includegraphics[height=0.44\textwidth]
  {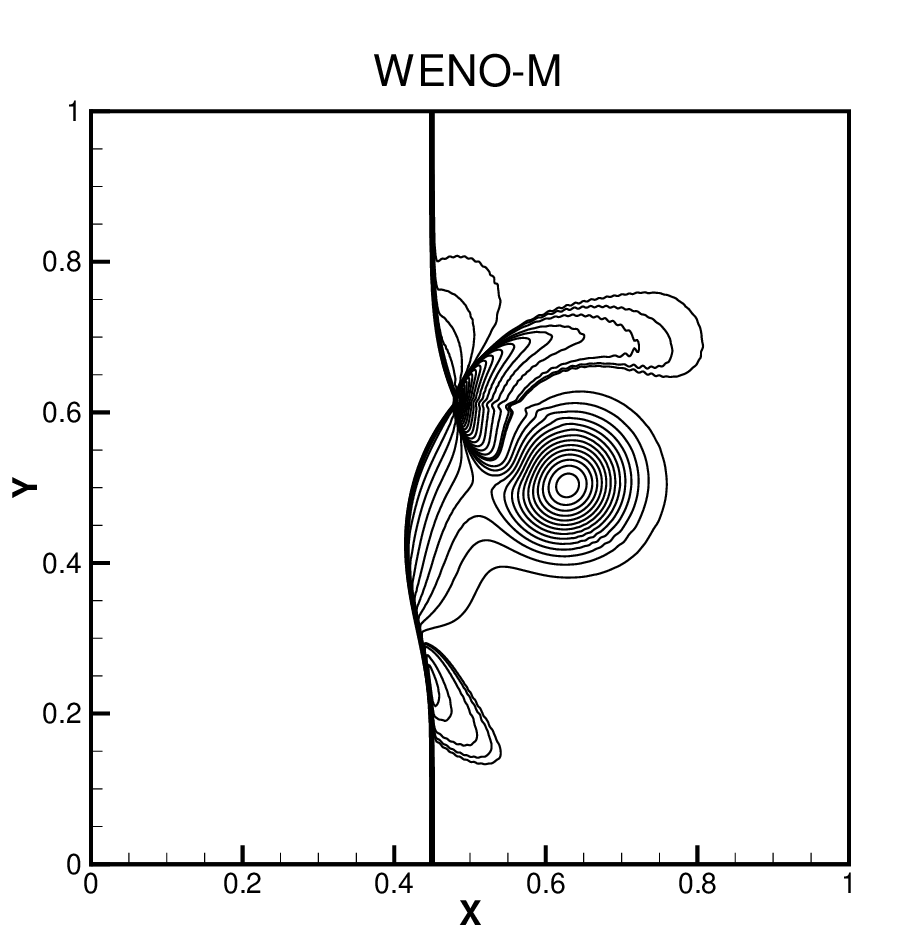}          \hspace{1.2ex}
  \includegraphics[height=0.44\textwidth]
  {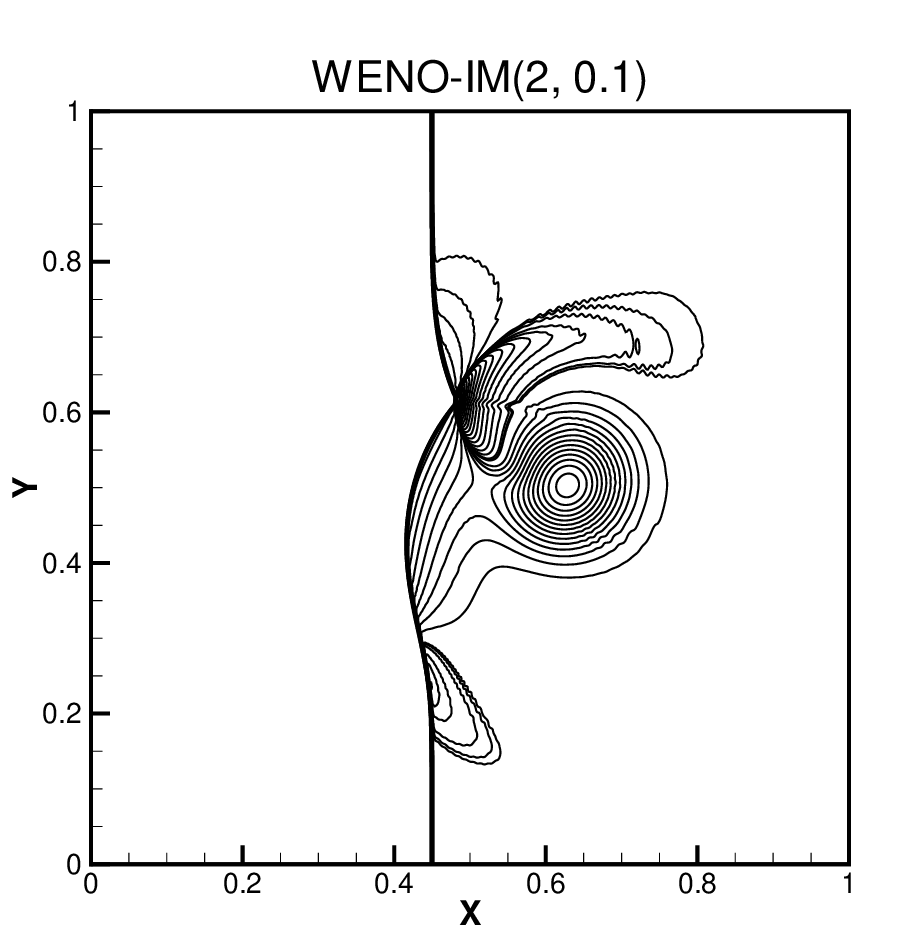}\\
  \includegraphics[height=0.44\textwidth]
  {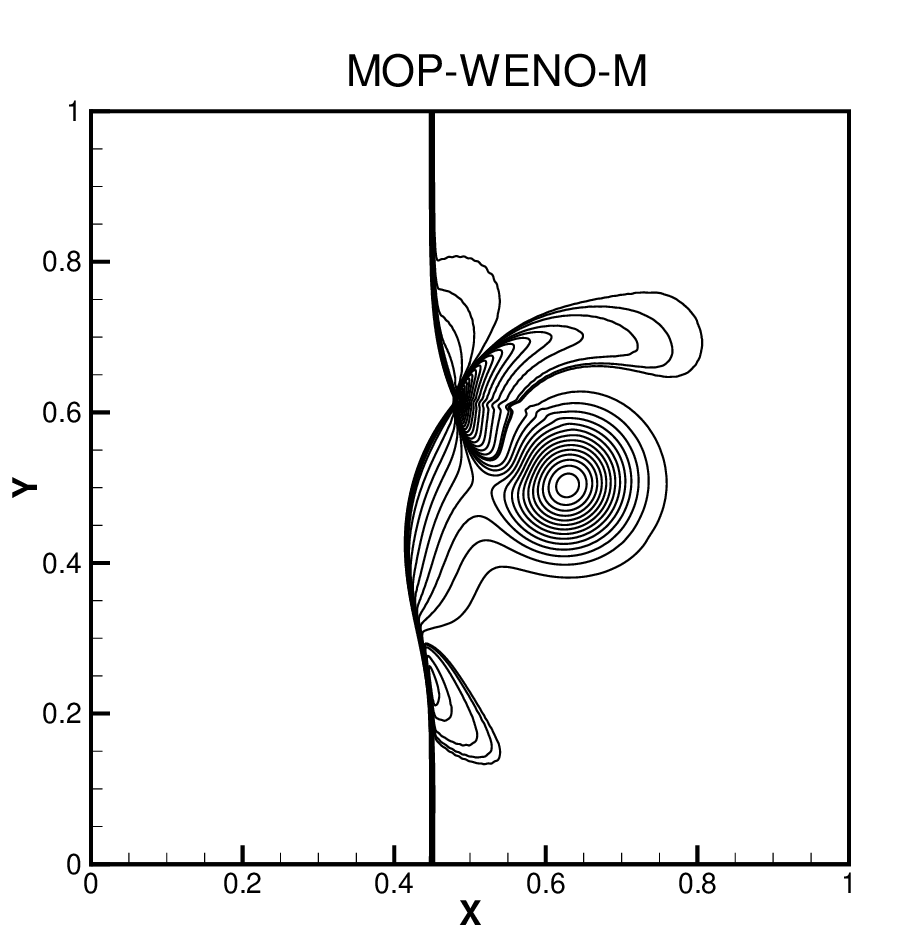}      \hspace{1.2ex}
  \includegraphics[height=0.44\textwidth]
  {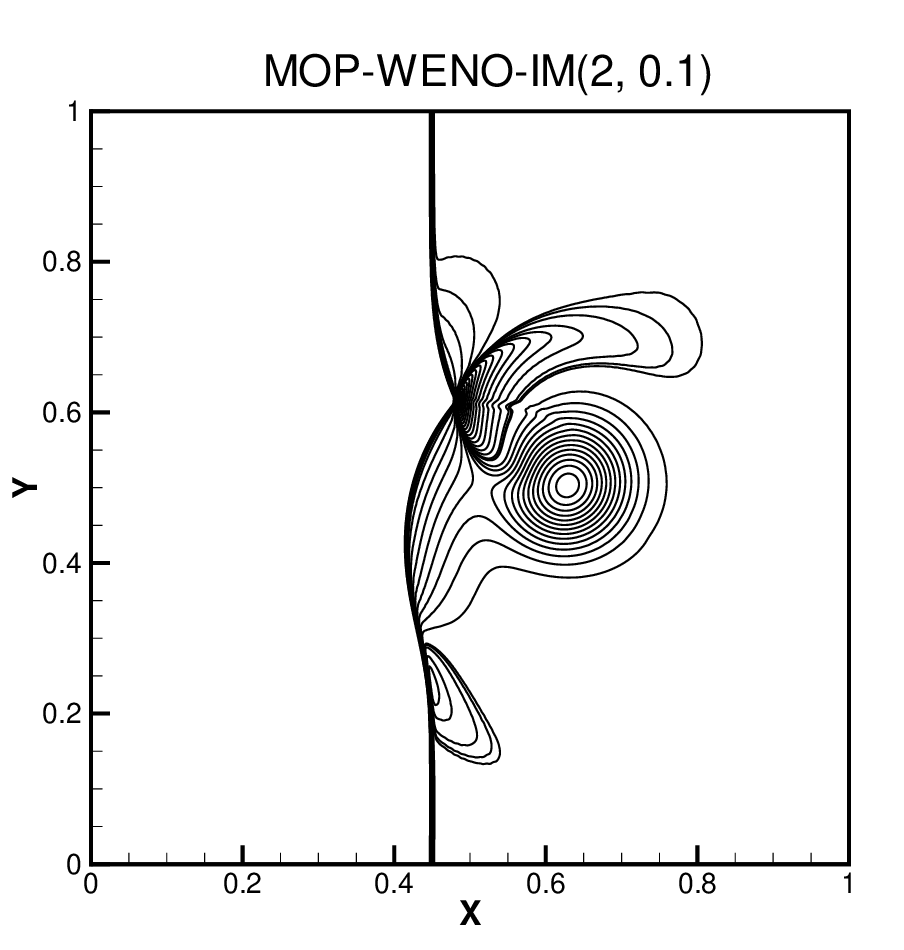}\\
\hspace{-3.5ex}  
  \includegraphics[height=0.345\textwidth]
  {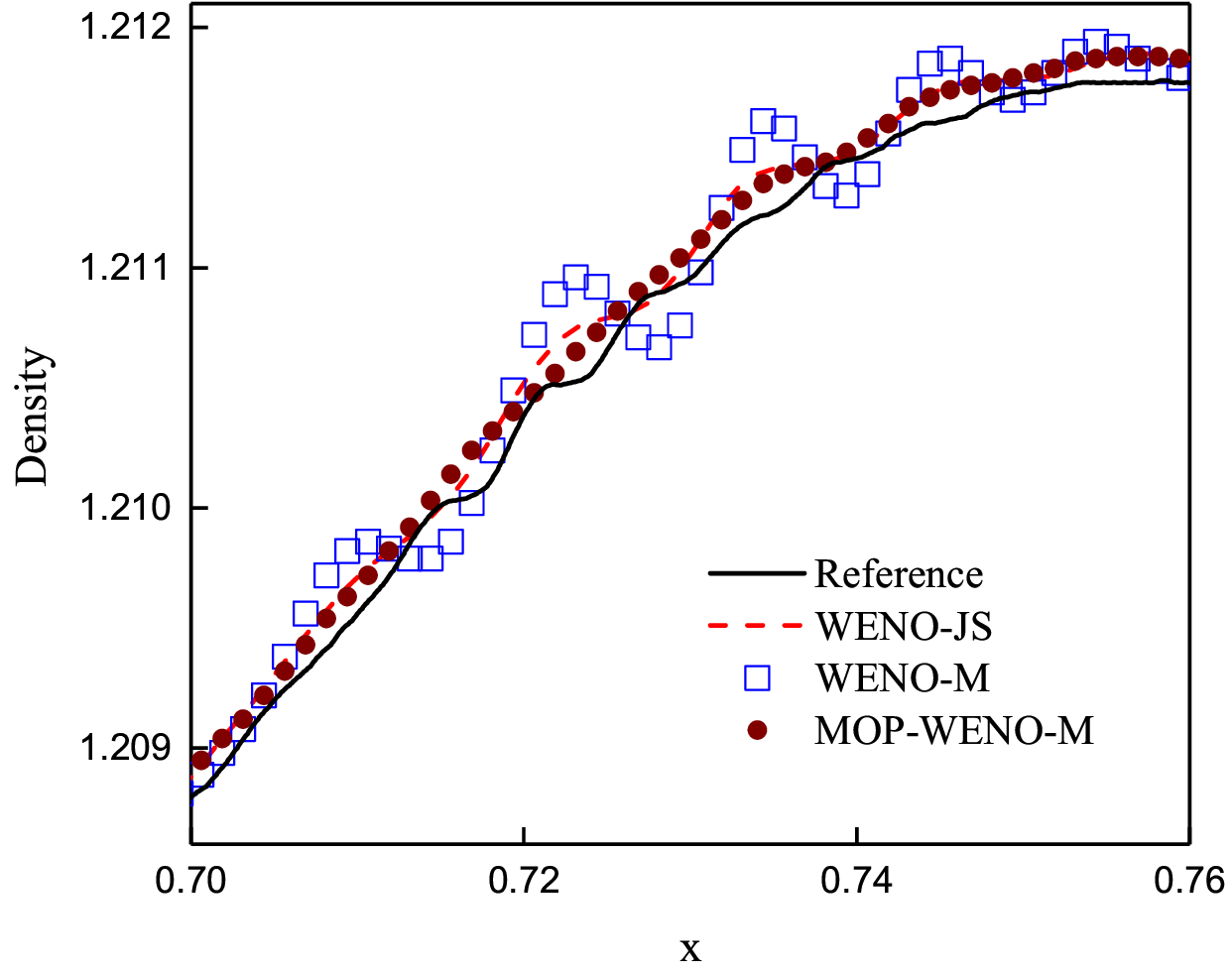}
  \includegraphics[height=0.345\textwidth]
  {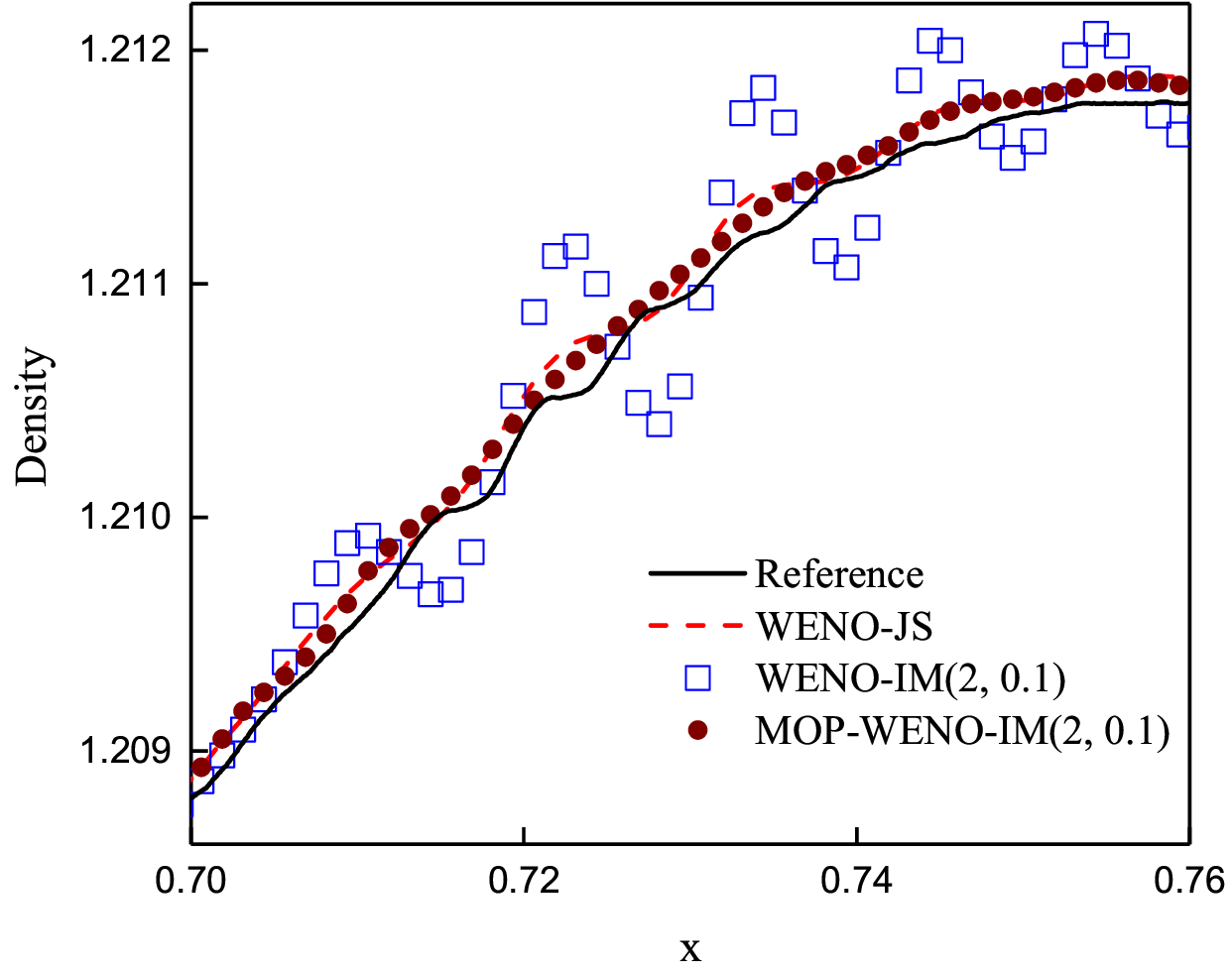}    
\caption{Density plots for the Shock-vortex interaction using $30$ 
contour lines with range from $0.9$ to $1.4$ (the first two rows) 
and the cross-sectional slices of density plot along the plane 
$y = 0.65$ where $x \in [0.70, 0.76]$ (the third row), computed 
using the WENO-M and MOP-WENO-M (left column), WENO-IM(2, 0.1) and 
MOP-WENO-IM(2, 0.1) (right column) schemes.}
\label{fig:ex:SVI:1}
\end{figure}

\begin{figure}[ht]
\centering
  \includegraphics[height=0.44\textwidth]
  {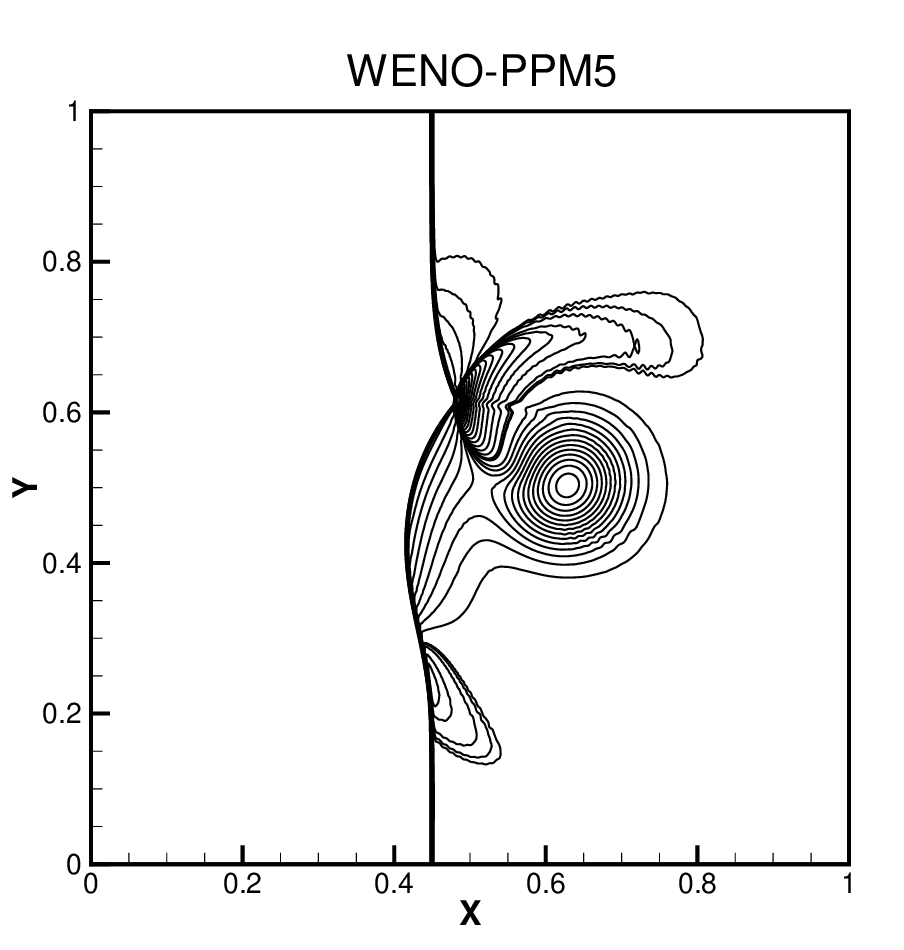}         \hspace{1.2ex}
  \includegraphics[height=0.44\textwidth]
  {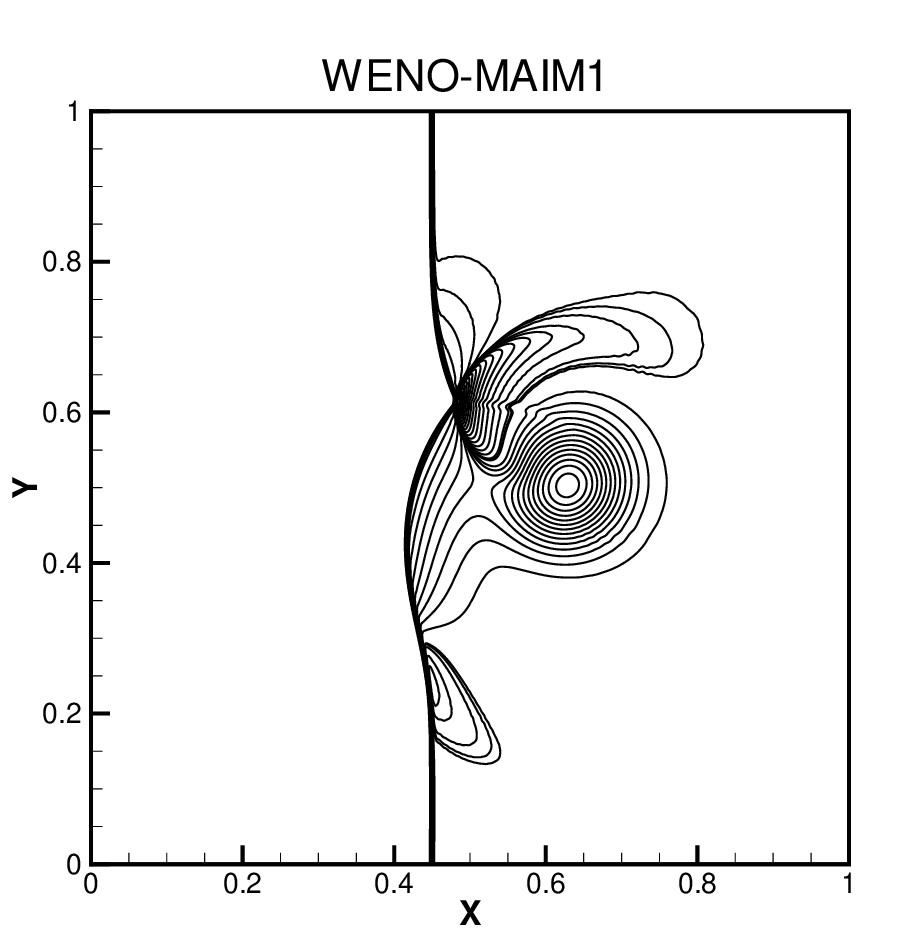}\\
  \includegraphics[height=0.44\textwidth]
  {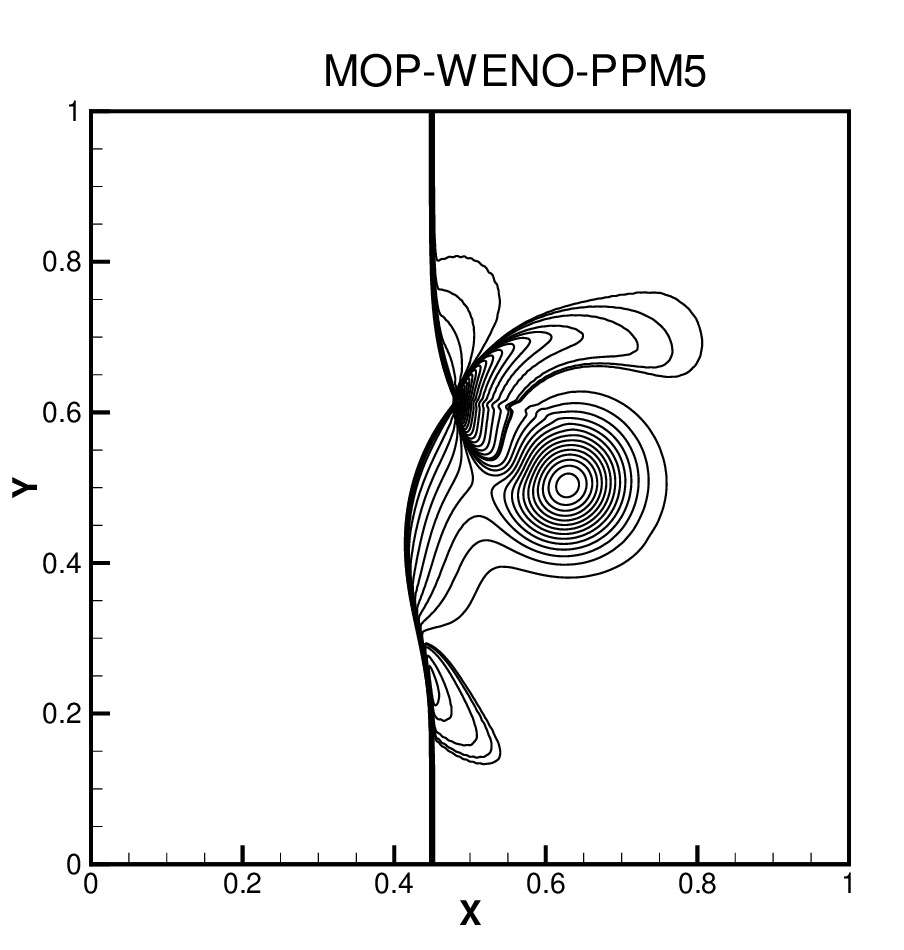}     \hspace{1.2ex}
  \includegraphics[height=0.44\textwidth]
  {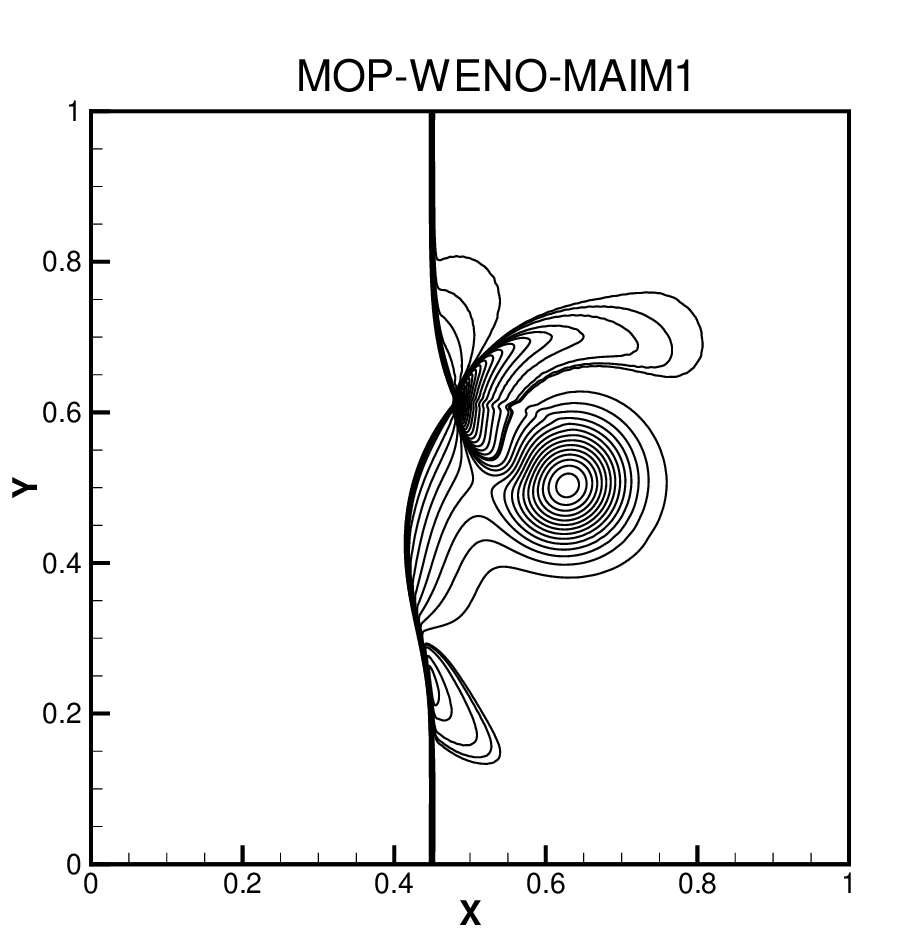}\\
\hspace{-3.5ex}
  \includegraphics[height=0.345\textwidth]
  {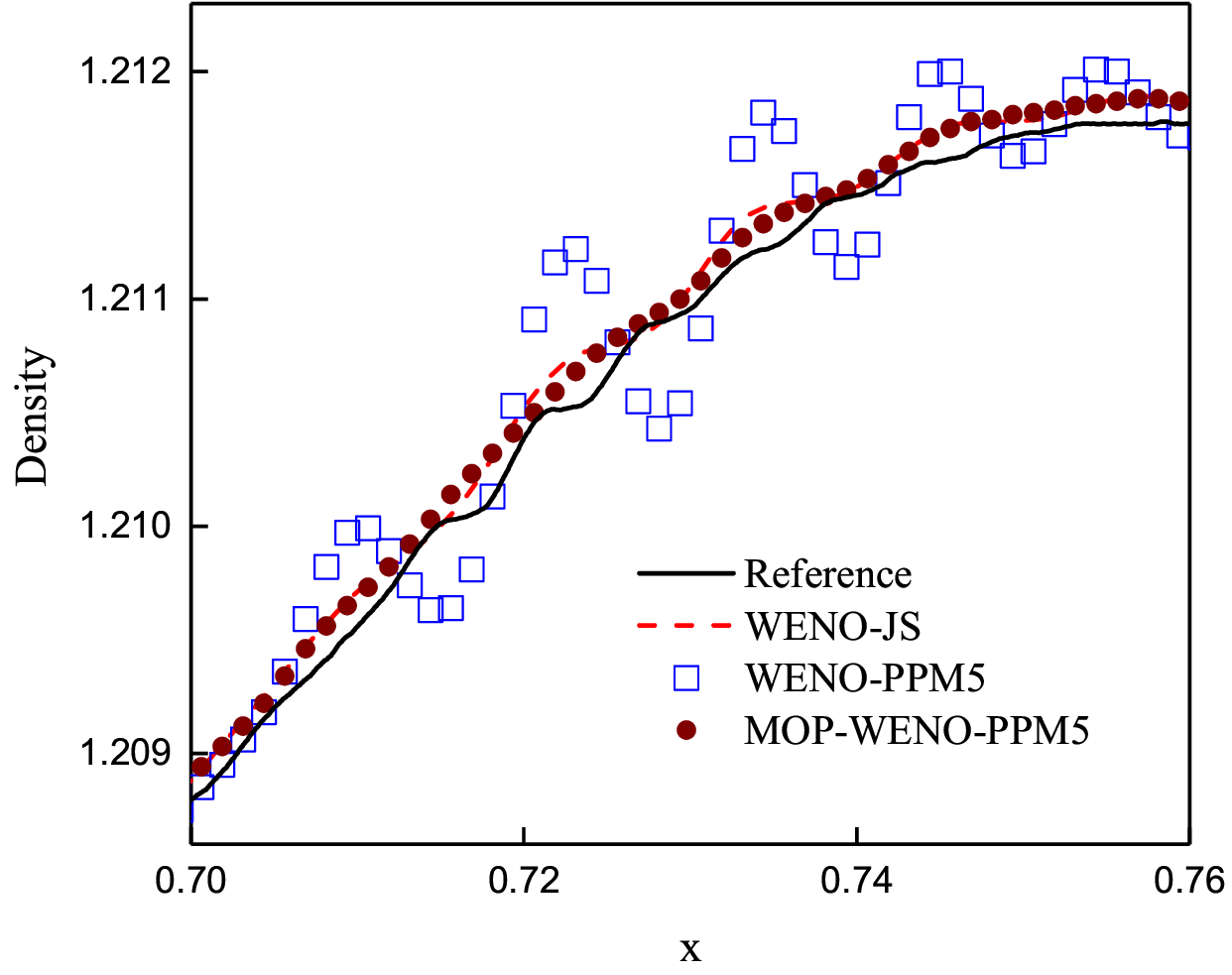}   
  \includegraphics[height=0.345\textwidth]
  {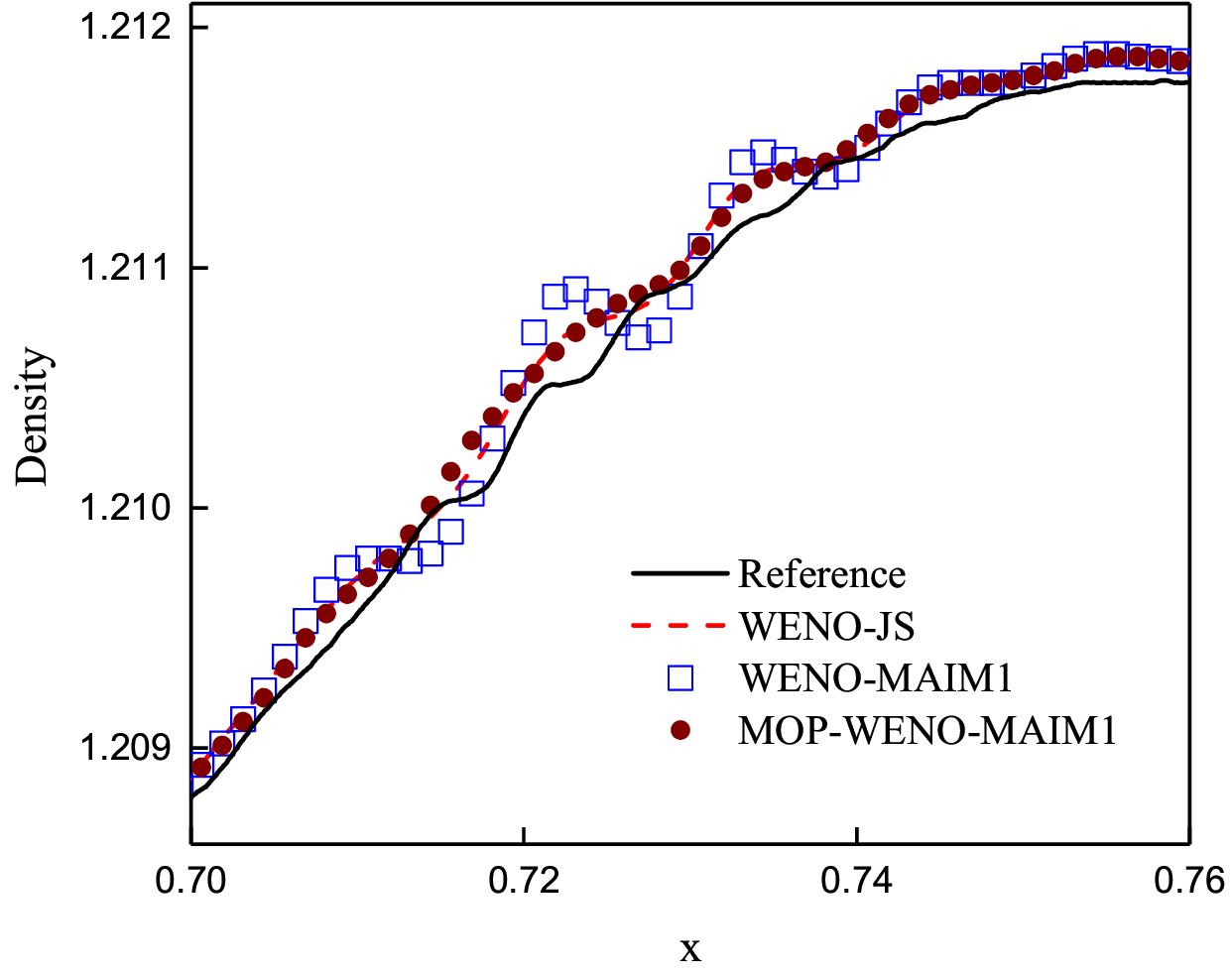}    
\caption{Density plots for the Shock-vortex interaction using $30$ 
contour lines with range from $0.9$ to $1.4$ (the first two rows) 
and the cross-sectional slices of density plot along the plane 
$y = 0.65$ where $x \in [0.70, 0.76]$ (the third row), computed 
using the WENO-PPM5 and MOP-WENO-PPM5 (left column), WENO-MAIM1 and 
MOP-WENO-MAIM1 (right column) schemes.}
\label{fig:ex:SVI:2}
\end{figure}

\begin{figure}[ht]
\centering
  \includegraphics[height=0.432\textwidth]
  {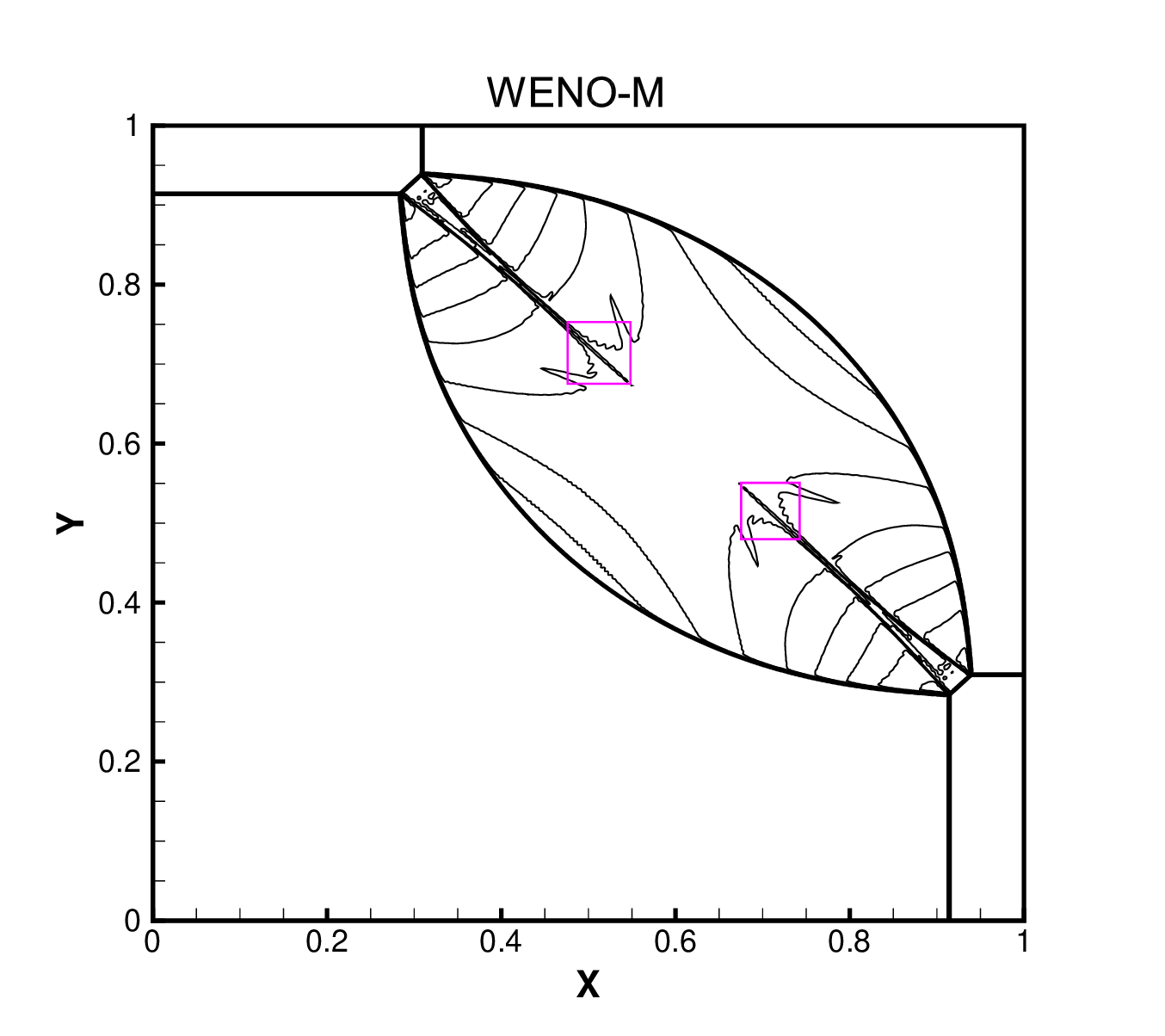}
  \includegraphics[height=0.432\textwidth]
  {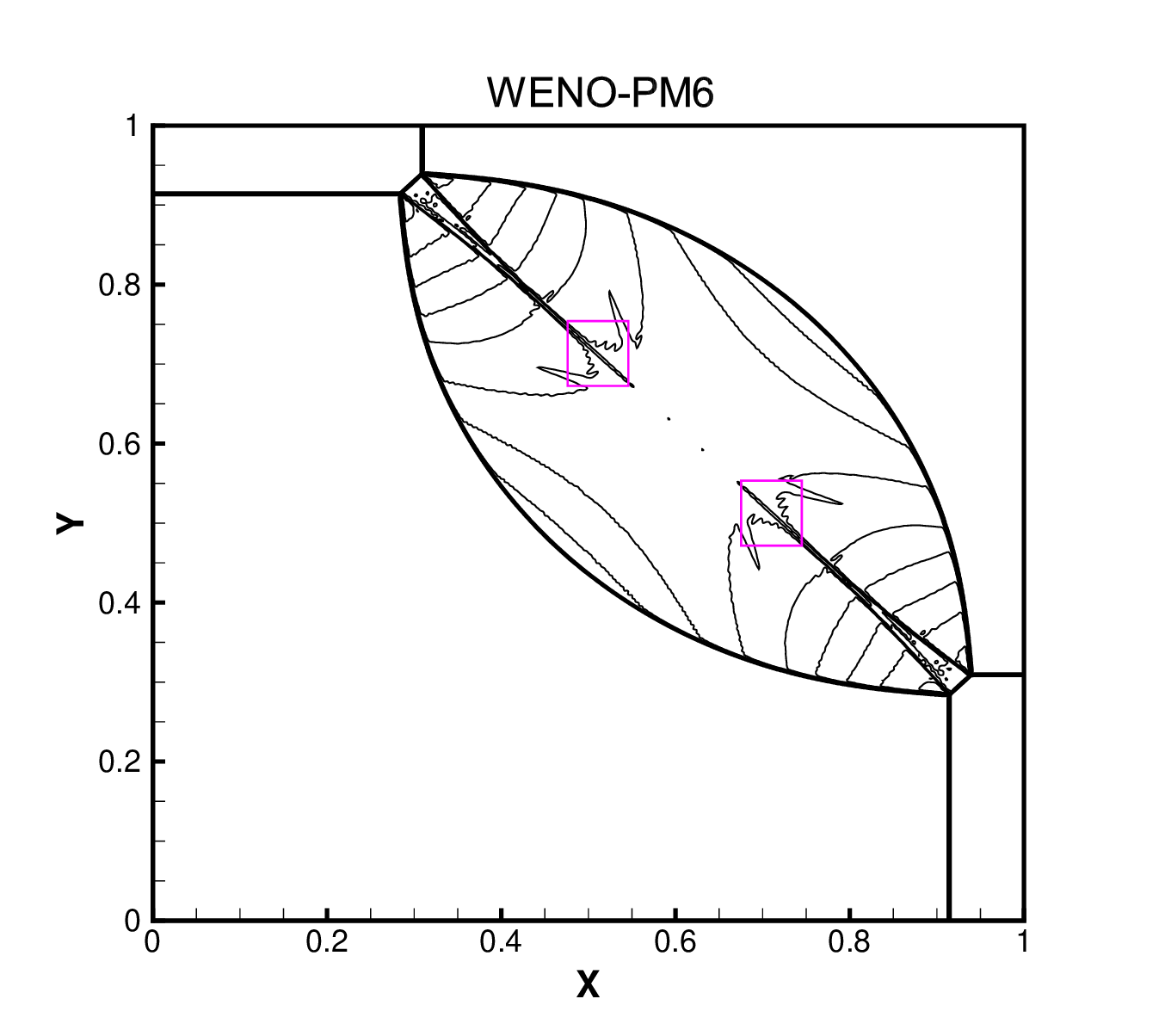}\\
  \includegraphics[height=0.432\textwidth]
  {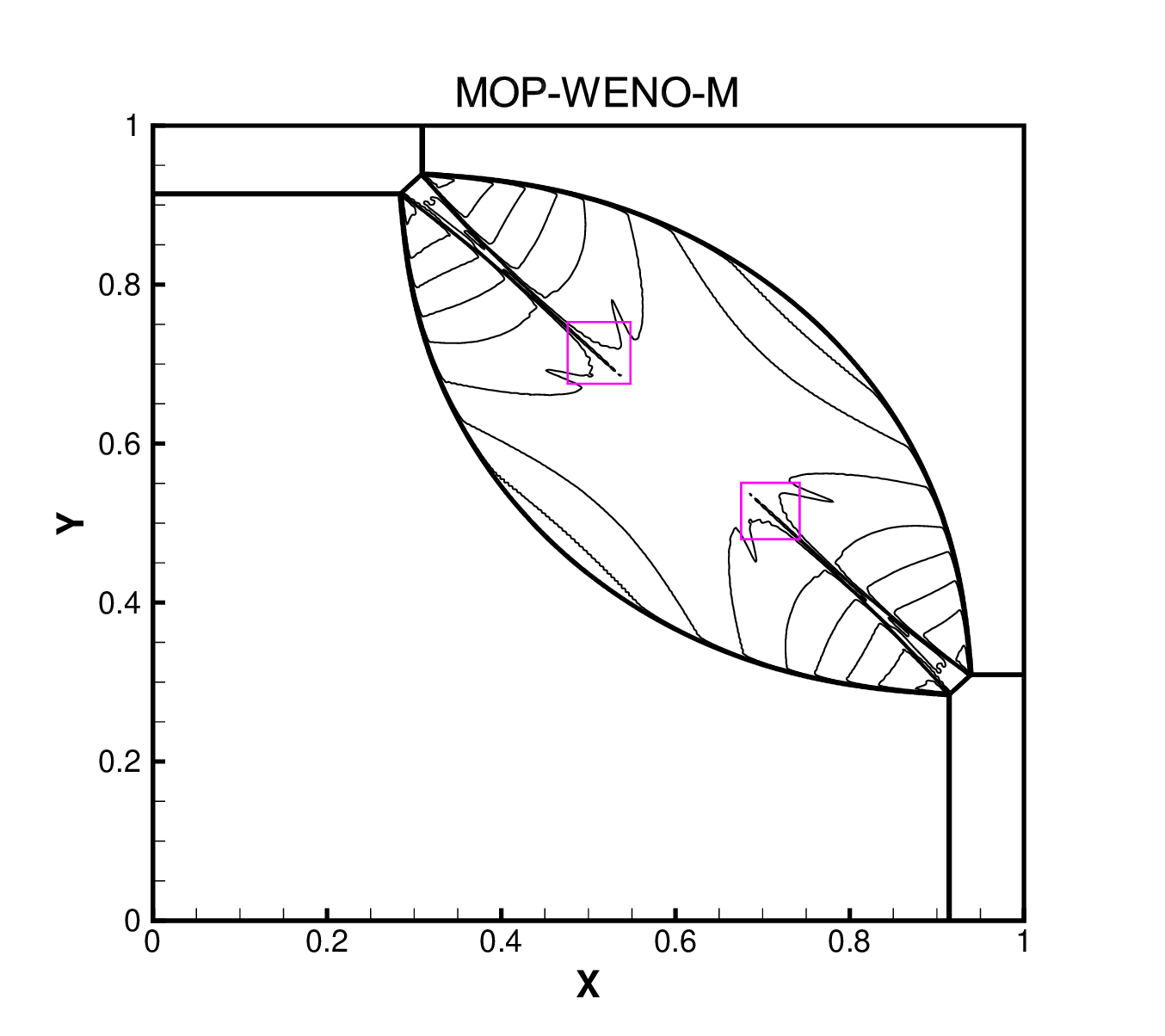}
  \includegraphics[height=0.432\textwidth]
  {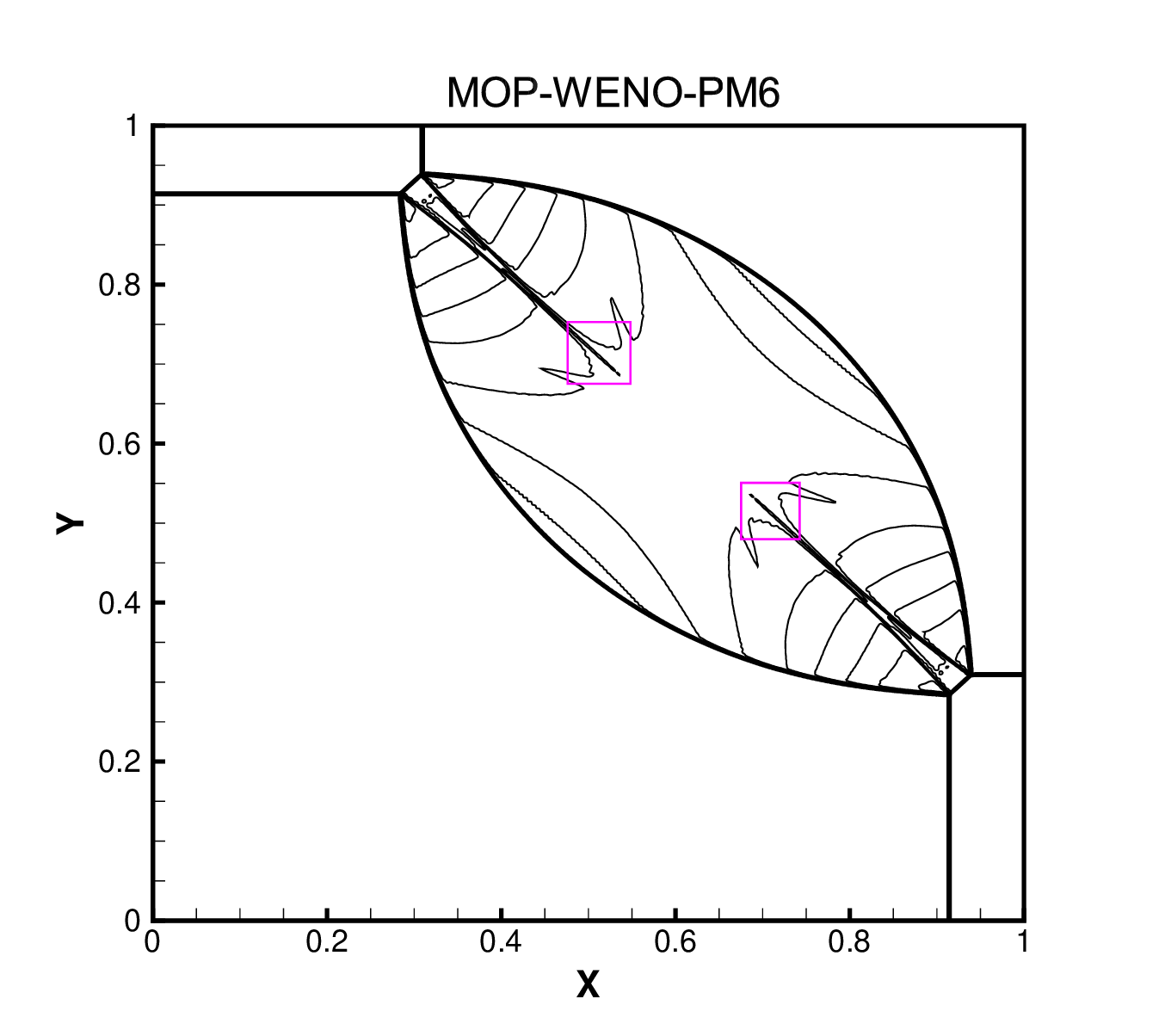}\\
\hspace{-6.2ex}  
  \includegraphics[height=0.355\textwidth]
  {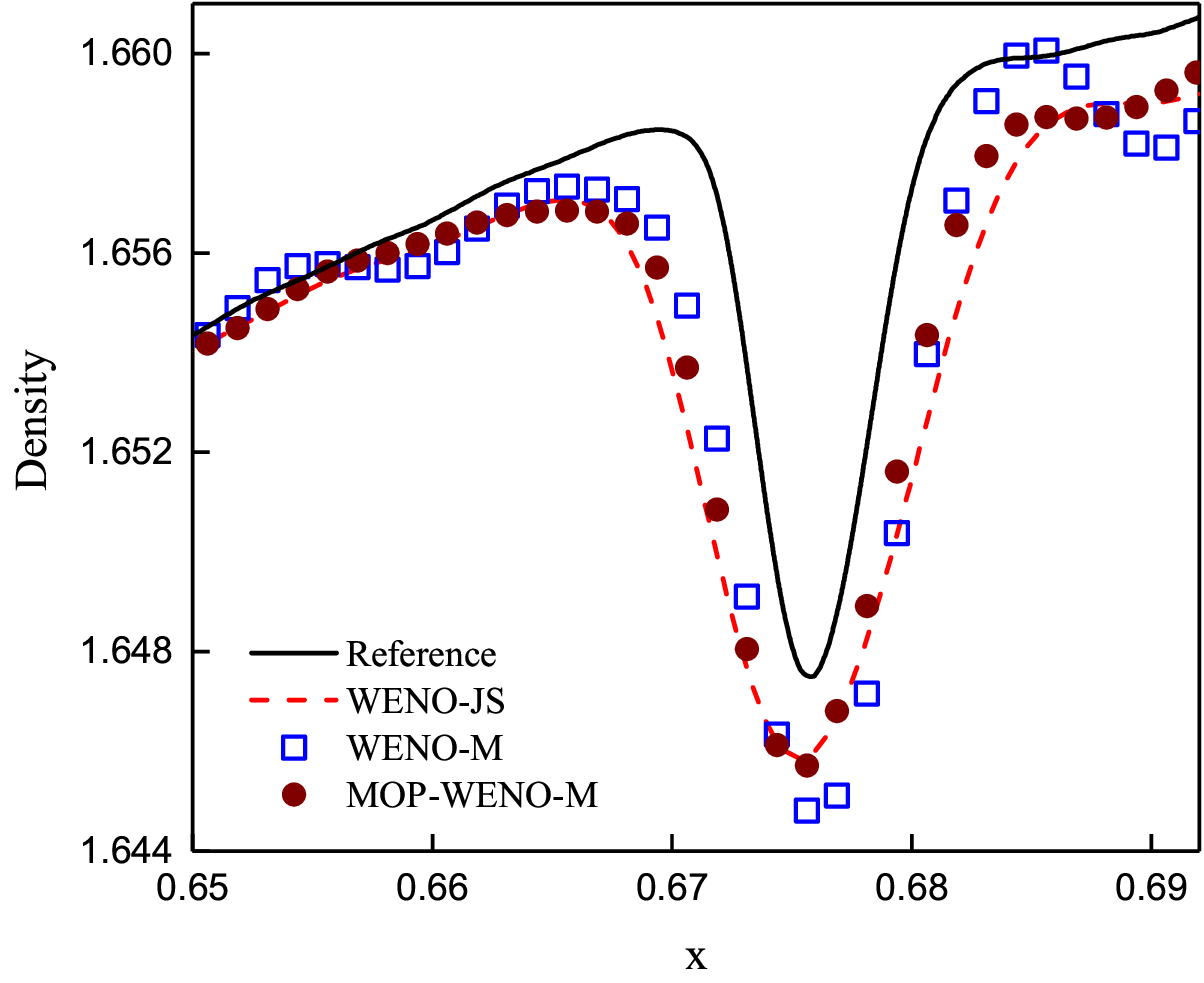}             \hspace{4.8ex}
  \includegraphics[height=0.355\textwidth]
  {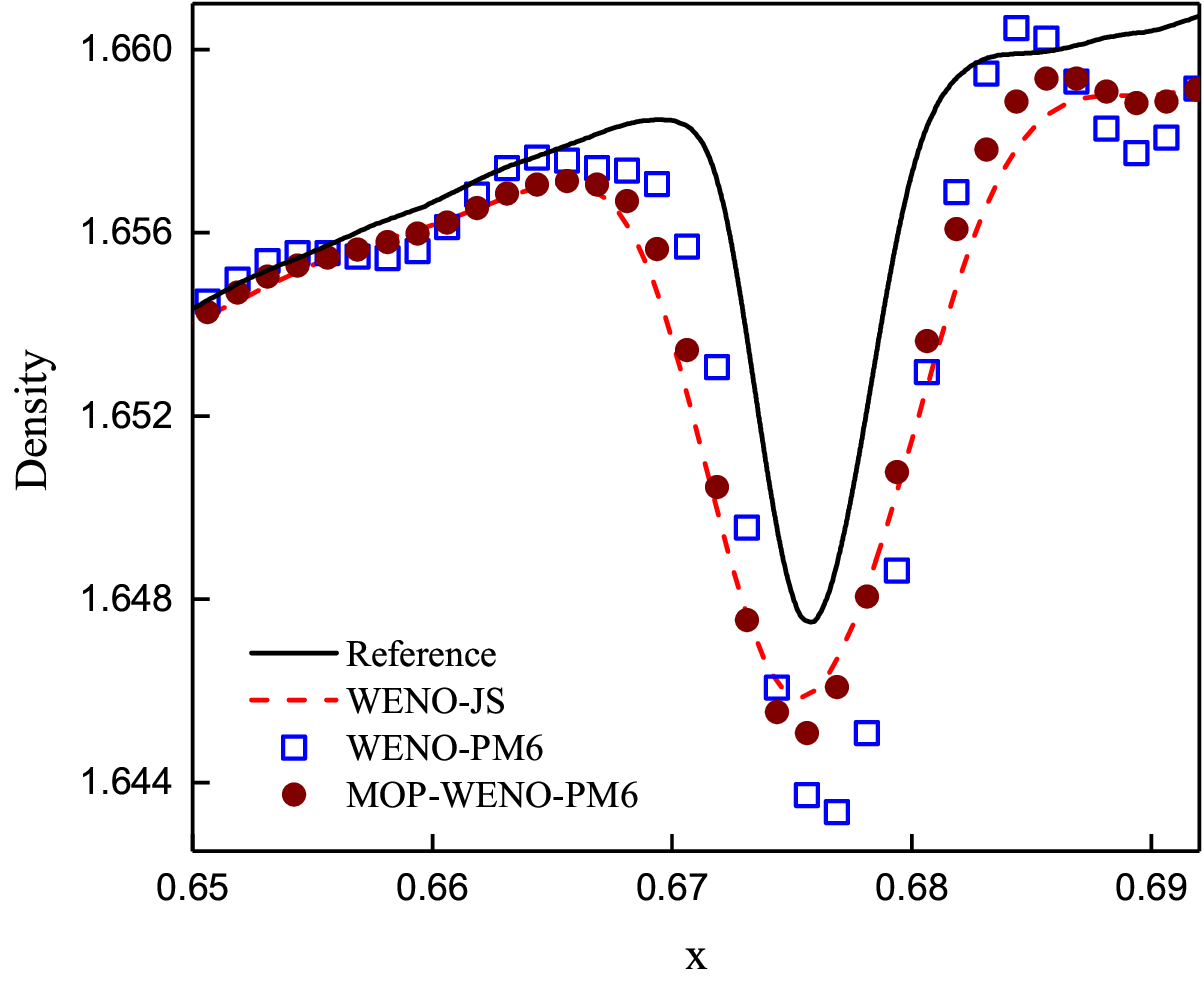}    
\caption{Density plots for the 2D Riemann problem using $30$ contour
lines with range from $0.5$ to $1.9$ (the first two rows) and the 
cross-sectional slices of density plot along the plane $y = 0.5$ 
where $x \in [0.65, 0.692]$ (the third row), computed using the 
WENO-M and MOP-WENO-M (left column), WENO-PM6 and MOP-WENO-PM6 
(right column) schemes.}
\label{fig:ex:Riemann2D:1}
\end{figure}

\begin{figure}[ht]
\centering
  \includegraphics[height=0.432\textwidth]
  {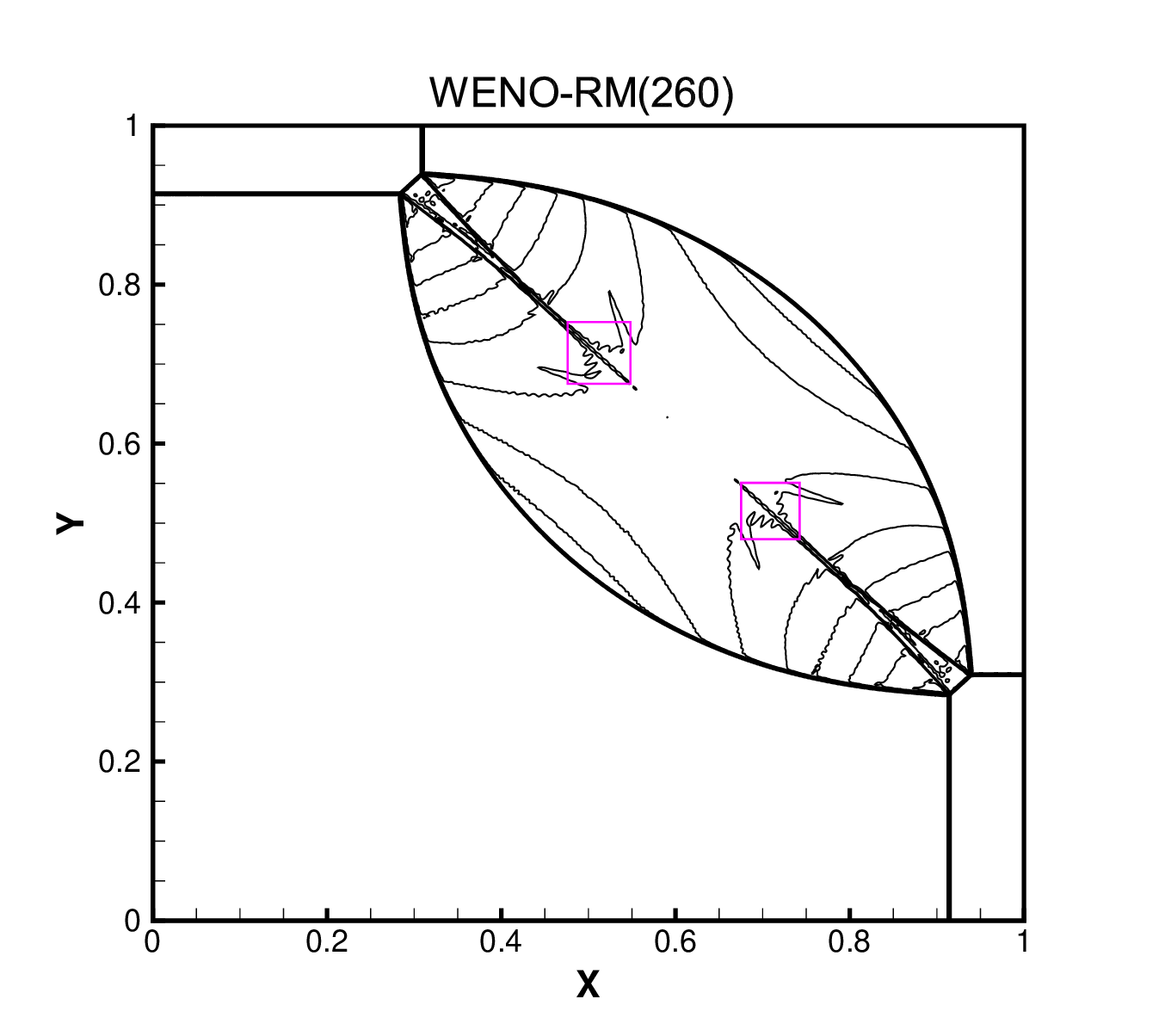}
  \includegraphics[height=0.432\textwidth]
  {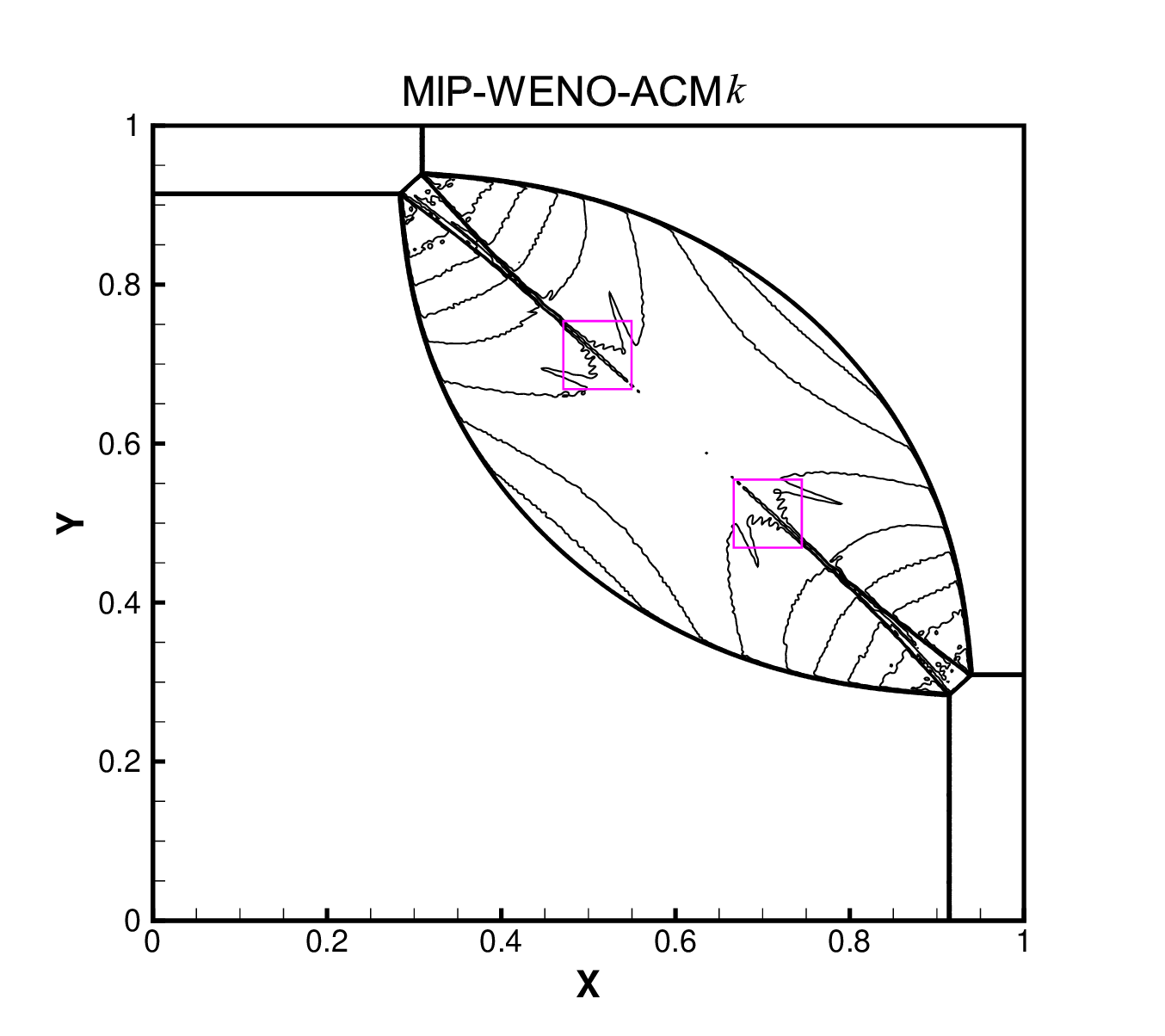}\\
  \includegraphics[height=0.432\textwidth]
  {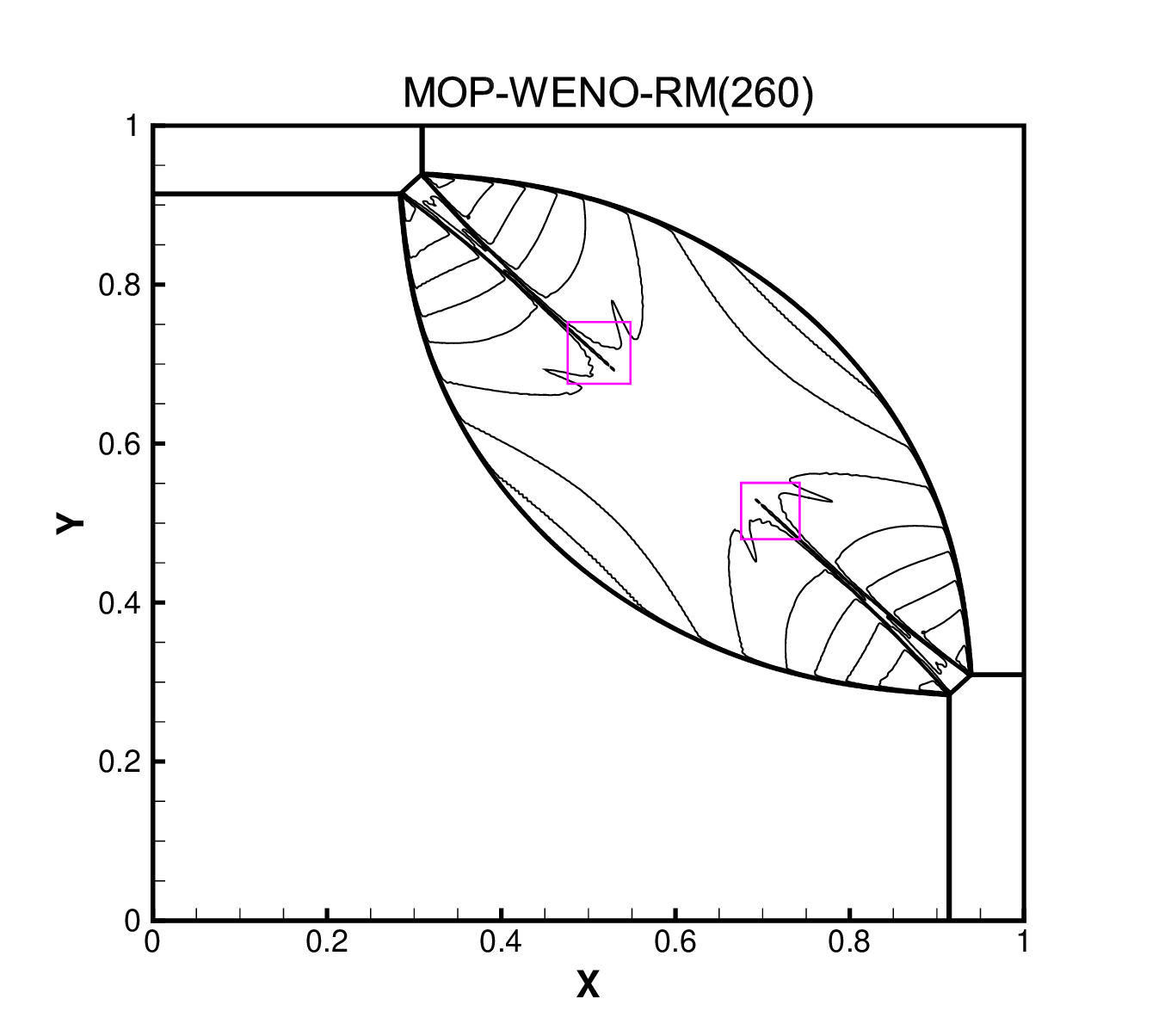}
  \includegraphics[height=0.432\textwidth]
  {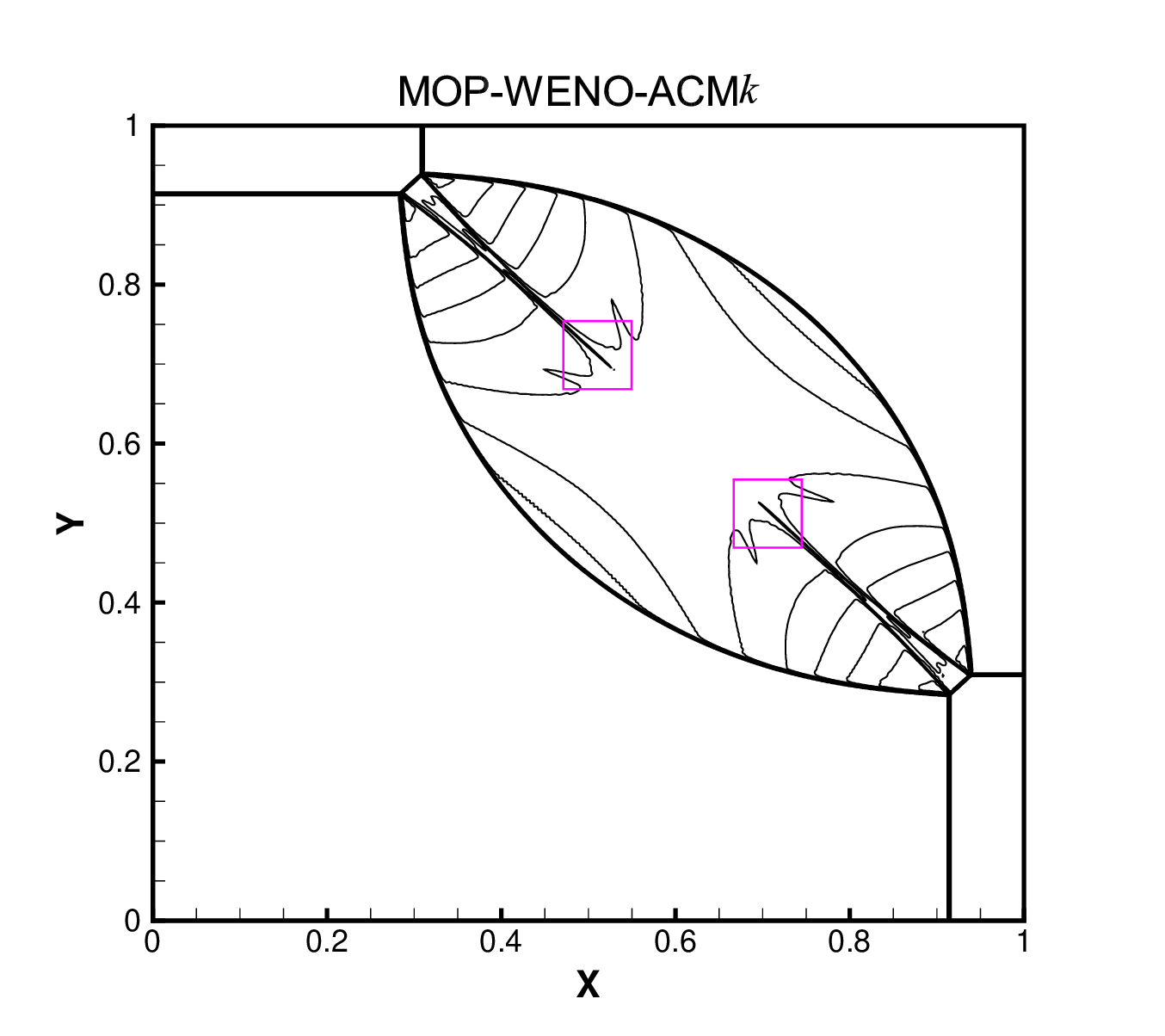}\\
\hspace{-6.2ex}    
  \includegraphics[height=0.355\textwidth]
  {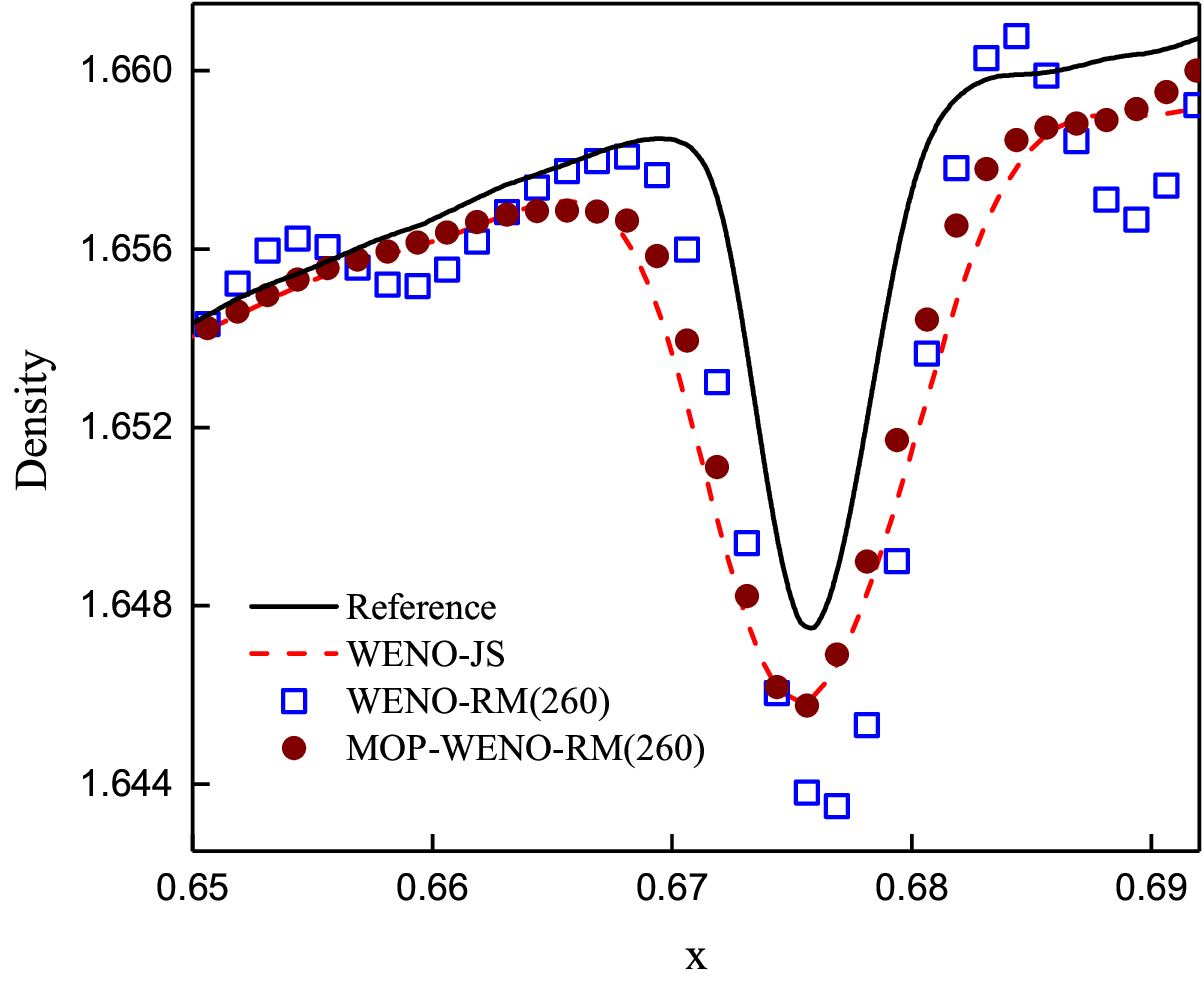}         \hspace{4.8ex}
  \includegraphics[height=0.355\textwidth]
  {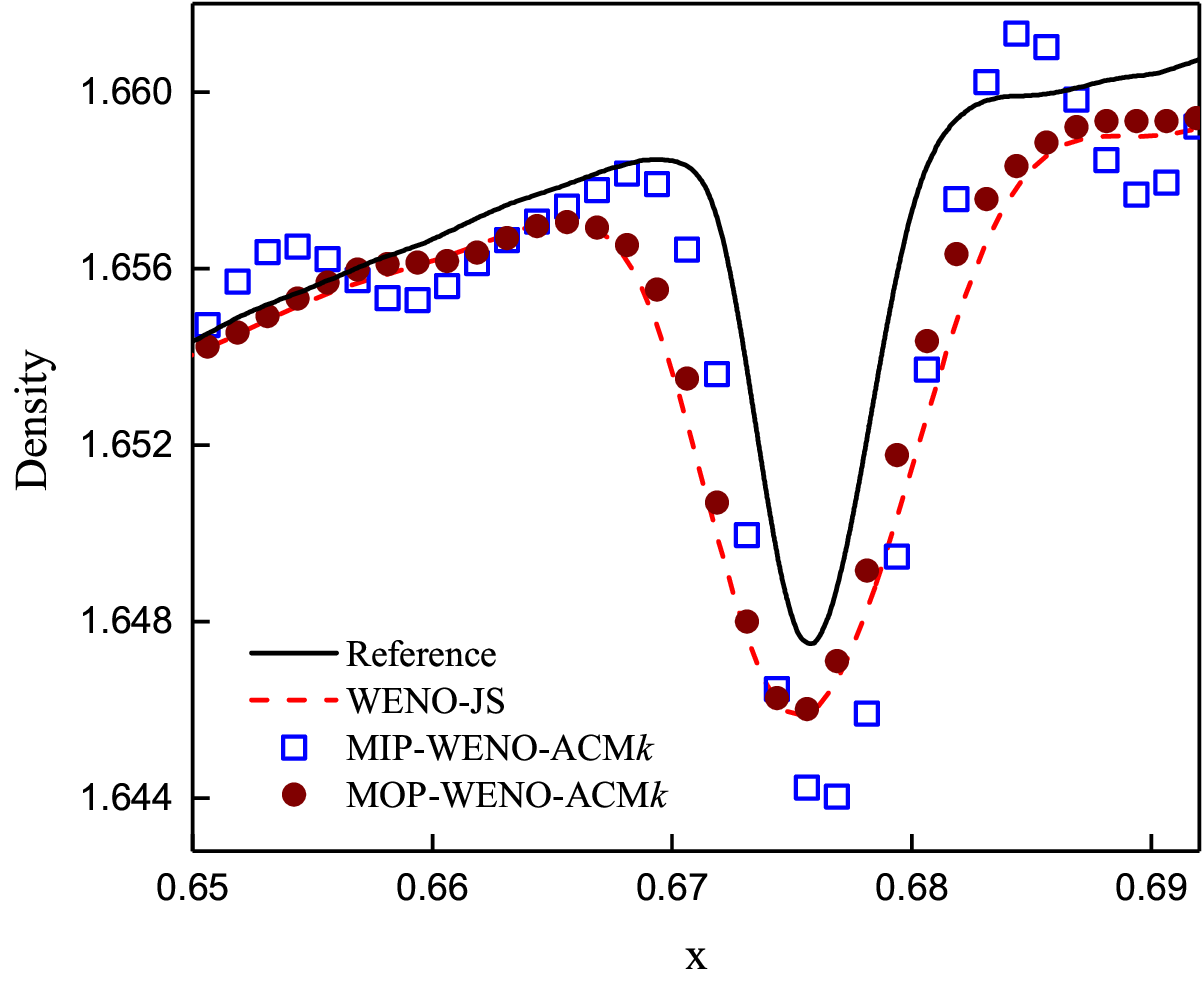}    
\caption{Density plots for the 2D Riemann problem using $30$ contour
lines with range from $0.5$ to $1.9$ (the first two rows) and the 
cross-sectional slices of density plot along the plane $y = 0.5$ 
where $x \in [0.65, 0.692]$ (the third row), computed using the 
WENO-RM(260) and MOP-WENO-RM(260) (left column), MIP-WENO-ACM$k$ and 
MOP-WENO-ACM$k$ (right column) schemes.}
\label{fig:ex:Riemann2D:2}
\end{figure}